\documentclass{amsart}
\usepackage{graphicx,amssymb,amscd}

\let\cal\mathcal

\def\xto#1{\xrightarrow[]{#1}}

\newtheorem{theorem}{Theorem}[section]
\newtheorem{proposition}{Proposition}[section]

\newtheorem{remark}{Remark}[section]
\newtheorem{definition}{Definition}[section]
\newtheorem{corollary}{Corollary}[section]
\newtheorem{example}{Example}[section]

\newcommand{\Sing}{\operatorname{Sing}}

\newcommand{\Id}{\operatorname{Id}}

\newcommand{\ox}{\otimes}
\newcommand{\x}{\times}
\newcommand{\uu}{\underline}
\newcommand{\oo}{\varnothing}
\newcommand{\sm}{\setminus}

\author{Tornike Kadeishvili}
\thanks{This research described in this publication was made possible in part by
Award No. GM1-2083 of the U.S. Civilian Research   and Development
Foundation for the Independent  States of the Former Soviet Union
(CRDF) and by Award  No. 99-00817 of INTAS}
\address{A. Razmadze Mathematical Institute\\
Georgian Academy of Sciences\\
M. Aleksidze st., 1\\
0193 Tbilisi, Georgia} \email{kade@@rmi.acnet.ge}

\author{Samson Saneblidze}
\address{A. Razmadze Mathematical Institute\\
Georgian Academy of Sciences\\
M. Aleksidze st., 1\\
0193 Tbilisi, Georgia} \email{sane@@rmi.acnet.ge}

\title{The twisted Cartesian model  for  the double path fibration}
\date{}

\begin{document}

\begin{abstract}
In the paper  the notion of  truncating twisting function from a
cubical set to a permutahedral set and the corresponding notion of
twisted Cartesian product of these sets are introduced. The latter
becomes a permutocubical set that models in particular the path
fibration on a loop space. The chain complex of this twisted
Cartesian product in fact is a comultiplicative twisted tensor
product of cubical chains of base and permutahedral chains of
fibre. This construction is formalized as a theory of twisted
tensor products for Hirsch algebras.
\end{abstract}

\maketitle

\vskip+5mm

\section{introduction}

The paper continues \cite{KS1} in which a combinatorial model for
a fibration was constructed based on the notion of a {\it
truncating twisting function} from a simplicial set to a cubical
set and on the corresponding notion of twisted Cartesian product
of these sets being a cubical set. Applying the cochain functor we
obtained  a {\it multiplicative} twisted tensor product modeling
the corresponding fibration.

There arises  a need to iterate this construction for fibrations
over loop or path spaces the bases of which are modeled  by
cubical sets. A cubical base naturally requires a permutahedral
fibre; this really agrees with the first usage of the permutahedra
(the Zilchgons) as modeling polytopes for loops on the standard
cube due to R.J. Milgram \cite{Milgram} (see also \cite{CM}).

For this  we proceed almost parallel to \cite{KS1}. Namely, let
$Q$ be a 1-reduced cubical set, ${\mathcal Z}$
 a monoidal permutahedral set, and ${\mathcal L}$ a permutahedral ${\mathcal Z}$-module, i.e.,
${\mathcal Z}$ and ${\mathcal L}$ are permutahedral sets with
given associative permutahedral maps ${\mathcal Z}\times {\mathcal
Z}\rightarrow {\mathcal Z}$ and ${\mathcal Z}\times {\mathcal
L}\rightarrow {\mathcal L}$ (see \cite{SU2} or Section 2 below).
We introduce the notion of {\it truncating twisting function}
$\vartheta :Q_*\to {\mathcal Z}_{*-1}$ from a  cubical set  to a
monoidal permutahedral set (the term {\it truncating} comes from
the universal example $\vartheta _{U}:I ^{n}\rightarrow P_{n}$ of
such functions obtained by the standard truncation procedure, see
Section \ref{trunc} below).
  Such a twisting function
 $\vartheta $ defines the
 twisted Cartesian product $Q\times _{\vartheta} {\mathcal L}$ as
 a
 {\it permutocubical set}. The {\it permutocube} is defined as a
polytope which is obtained from the standard cube by a certain
truncation procedure due to N. Berikashvili \cite{Berip}, see also
bellow. The permutocube can be thought of as a modeling polytope
for paths on the cube.

We construct a functor which assigns  to a cubical set $Q$ a
monoidal permutahedral set ${\bf \Omega} Q$ and present
 a truncating twisting function
  $\vartheta_U :Q\to {\bf \Omega} Q$   of degree $-1$ which is
universal in the following sense: Given an arbitrary truncating
function $\vartheta :Q_* \to {\mathcal Z}_{*-1},$ there is
 a monoidal
permutahedral map  $f_{\vartheta }: Q \to{\mathcal Z}$ such that
$\vartheta=f_{\vartheta}\vartheta_U.$ The twisted Cartesian
product ${\bf P }Q=Q\times_{\vartheta_U}{\bf \Omega} Q$ is a
permutocubical set functorially depending on $Q$. Note that ${\bf
\Omega} Q$ models the loop space $\Omega |Q|$ and ${\bf P }Q$
models the path fibration on $|Q|$.

The  chain complex $C^{\diamondsuit}_*({\bf \Omega} Q)$ coincides
with the cobar construction $\Omega C^{\Box}_*(Q)$. Similarly, the
chain complex $C^{\boxminus}_*(Q  \times_{\vartheta_{U} }{\bf
\Omega} Q)$ coincides with the acyclic cobar construction
$\Omega\left(C^{\Box}_*(Q)\,;C^{\Box}_*(Q)\right);$ furthermore,
$\vartheta _*= C_*(\vartheta ):C^{\Box}_*(Q)\to
C^{\diamondsuit}_{*-1}({\mathcal Z})$ is a twisting  cochain and
$C^{\boxminus}_*(Q\times _{\vartheta } {\mathcal L})$ coincides
with the twisted tensor product  $C^{\Box}_*(Q)\otimes
_{\vartheta_*}C^{\diamondsuit}_*({\mathcal L}).$

  We construct an explicit diagonal for the permutocube $B_{n}$
which agrees with that of $P_{n}$  \cite{SU2}
 by means of the natural embedding $P_{n}\to  B_{n}$. The
equalities $C^{\diamondsuit}_*({\bf \Omega} Q)=\Omega
C^{\Box}_*(Q)$ and $C^{\boxminus}_*(Q\times _{\vartheta }{\mathcal
L})=C^{\Box}_*(Q)\otimes
_{\vartheta_*}C^{\diamondsuit}_*({\mathcal L})$ allow us to
transport these diagonals \linebreak
 to the cobar construction $\Omega
C^{\Box}_*(Q)$ and the twisted tensor product
$C^{\Box}_*(Q)\otimes _{\vartheta_*}C^{\diamondsuit}_*({\mathcal
L})$ respectively. Dually, we immediately obtain a multiplication
on
 $
 C^*_{\Box}(Q)\otimes
_{\vartheta^*}C_{\diamondsuit}^*({\mathcal L}) \subset
C_{\boxminus}^*(Q\times _{\vartheta }{\mathcal L})$  (which is an
equality  if the graded sets are of finite type).
 Note that
this (co)multiplication is not strictly (co)associative but could
be extended to an $A_{\infty}$-(co)algebra structure.

Next we express the resulting comultiplication on
$C^{\Box}_*(Q)\otimes _{\vartheta_*}C^{\diamondsuit}_*({\mathcal
L})$ in terms of certain chain operations of degree $p+q-1:$
\[
\{E^{p,q}:C^{\Box}_{*}(Q)\to C^{\Box}_{*}(Q)^{\otimes p}\otimes
C^{\Box}_{*}(Q)^{\otimes q}\}_{p+q>0},
\]
which give $C^{\Box}_{*}(Q)$ a structure what we call a {\em
Hircsh coalgebra structure.} This structure is a consequence of
the permutahedral diagonal on $C^{\diamondsuit}_*({\bf \Omega}
Q)=\Omega C^{\Box}_*(Q):$ The permutahedral diagonal from
\cite{SU2} induces
 the diagonal $\Omega
C^{\Box}_*(Q)\to \Omega C^{\Box}_*(Q)\otimes \Omega C^{\Box}_*(Q)$
being a multiplicative map, thus it extends  a certain
homomorphism $C^{\Box}_*(Q)\to \Omega C^{\Box}_*(Q)\otimes \Omega
C^{\Box}_*(Q)$, which itself consists of components $ E^{p,q}:
C^{\Box}_*(Q)\to  C^{\Box}_*(Q)^{\ox p}\otimes C^{\Box}_*(Q)^{\ox
q},\,p,q\geq 0.$ The operation $E^{1,1}$ is dual to the cubical
version of Steenrod's $\smile _{1}$-cochain operation; thus when
$E^{1,1}=0$ a Hirsch coalgebra specializes to a cocommutative dg
coalgebra (and dually for {\em Hirsch algebras}).

Towards the end of the paper we develop the theory of
multiplicative twisted tensor products for  Hirsch algebras, which
provides a general algebraic framework for our multiplicative
model of a fibration.  A Hirsch algebra we define as an object
$(A,d,\cdot ,\{E_{p,q}:A^{\otimes p}\otimes A^{\otimes q}\to
A\}_{p+q>0}),$ i.e., $(A,d,\cdot )$ is an associative  dga and the
sequence of operations $\{E_{p,q}\}$ determines a product on the
bar construction $BA$ turning it into a dg Hopf algebra (this
multiplication can be viewed as a perturbation of the shuffle
product and is not necessarily associative). In particular
$E_{1,1}$ has properties similar to $\smile_1$ product, so that a
Hirsch algebra can be considered as to have a structure measuring
the lack of commutativity of $A$. Let $C$ be a dg Hopf algebra and
$M$ be a dga and a dg $C$-comodule simultaneously. We say that a
twisting cochain $\phi:C\to A$ is {\it multiplicative} if the
induced map $C\to BA$ is a dg Hopf algebra map.  We introduce on
$A\otimes_{\phi}M$ a twisted multiplication $\mu_{\phi }$ in terms
of $\phi $ and the Hirsch algebra structure of $A$ by the same
formulas as in the case $A=C_{\Box}^*(Q),$
$C=C_{\diamondsuit}^*({\mathcal Z})$
 and $M=C_{\diamondsuit}^{*}(
{\mathcal L});$ then $\phi=\vartheta^*: C_{\diamondsuit}^{*}(
{\mathcal Z})\to C_{\Box}^{*+1}(Q)$ provides a basic example of a
multiplicative twisting cochain.
\newpage

Applying  our machinery to a fibration $F\to E\to Y$ on a
1-connected space $Y$ and an associated  principal $G$-fibration
$G\to P\to Y$ with action $G\times F\to  F$ we obtain the
following permutocubical model (Theorem \ref{percubmodel}): Let
$Q={{\Sing}^1}^I Y\subset {{\Sing}}^I Y$ be the Eilenberg
1-subcomplex  generated by singular cubes that send the 1-skeleton
of the standard $n$-cube $I^n$ to the base point of $Y.$ Let
${\mathcal Z}={{\Sing}}^M G$ and $Y={{\Sing}}^M F$ be the singular
 multipermutahedral sets (see \cite{SU2} and Section 2). We
construct the Adams-Milgram map
$$\omega_*:\Omega
C^{\Box}_*(Q)\to C_*^{\diamondsuit}(\Omega Y)$$ which  is realized
by a monoidal multipermutahedral map $\omega:{\bf \Omega}
Q\rightarrow {\Sing}^M\Omega Y.$ Composing $\omega$ with a map of
monoidal multipermutahedral sets\,
 ${\Sing}^M\Omega Y  \rightarrow  \linebreak  {\Sing}^M G
= {\mathcal Z}$  induced by the canonical map $\Omega Y\to G$ of
monoids we immediately obtain a truncating twisting function
$\vartheta:Q\to {\mathcal Z}$. The resulting twisted Cartesian
product ${\Sing^1}^I Y\times_{\vartheta }{\Sing}^M F$ provides the
required permutocubical model of $E;$ and there exists a
permutocubical weak equivalence ${\Sing^1}^I Y\times_{\vartheta
}{\Sing}^M F\to {\Sing}^{B}E,$ where ${\Sing}^{B}$ denotes the
singular permutocubical complex of a space. Applying the cochain
functor we obtain a certain multiplicative twisted tensor product
for the fibration.

In particular, we can obtain a permutocubical model for the path
fibration $\Omega^2 Y'\to P\Omega Y' \to \Omega Y'$ in the
following way. Taking for the base $Y=\Omega Y'$ the cubical model
$Q={\bf\Omega} \Sing^2 Y' $ from \cite{KS1}  the above  machinery
yields the twisted Cartesian model ${\bf \Omega} \Sing^2
Y'\times_{\vartheta_U} {\bf\Omega}{\bf \Omega} \Sing^2 Y' $ being
a {\em permutocubical} set.

Consequently, we introduce the multiplication on the acyclic bar
construction $B\left(BC^*(Y);BC^*(Y)\right)$  whose restriction to
the double bar construction  $BB C^*(Y)$
 is just the one constructed in
   \cite{SU2}.

To summarize we observe the following. In \cite{KS1} it is
indicated the homotopy G-algebra structure on $C^*(Y)$ consisting
of cochain operations
$$
\{E_{k,1}:C^*(Y)^{\otimes k}\otimes C^*(Y)\to C^*(Y)\}_{k \geq 1
},
$$
defining a multiplication on $BC^*(Y)$. Here we extend this
multiplication to  the structure of Hirsch algebra on $BC^*(Y)$,
i.e., to  operations
$$
\{E_{p,q}:(BC^*(Y))^{\otimes p}\otimes (BC^*(Y))^{\otimes q}\to
BC^*(Y)\}_{p+q>0},
$$
which actually are cochain operations of type $C^*(Y)^{\otimes m
}\to C^*(Y)^{\otimes n}$. This two sets of operations including in
particular $\smile$, $ \smile_1$ and $\smile_2$ operations, allow
us to construct multiplicative models for $\Omega Y,\ \Omega^2Y$
and multiplicative twisted tensor products for path fibrations on
$Y$ and $\Omega Y$ as well as for fibrations associated with them.

As an example we present fibrations with the base  being  the loop
space on a double suspension (in this case the Hirsch algebra
structure consists just of $E_{1,1}=\smile_1$ and all other
operations $E_{p,q}$ are trivial) and for which the formula for
the multiplication in the twisted tensor product has a very simple
form. Moreover, in this case we present small multiplicative model
being the twisted tensor product of cohomologies of base and fiber
with the multiplicative structure purely defined by the $\smile,$
$\smile_1$ and  $\smile_2$ operations.

Finally, we mention that the geometric realization $|{\bf\Omega}
{\bf\Omega} \Sing ^2 Y|$ of ${\bf \Omega}{\bf\Omega} \Sing ^2 Y$
is homeomorphic to the cellular model for the double loop space
due to G. Carlsson and R. J. Milgram \cite{CM} and is
homotopically equivalent to the cellular model due to H.-J. Baues
\cite{Baues1}.

The paper is organized as follows. We adopt the notions and the
terminology from \cite{KS1}; note that here a (co)algebra need not
have a  (co)associative (co)mul\-ti\-pli\-ca\-ti\-on if it is not
specially emphasized.
 In Section 2 we  construct the  functor $\bf \Omega$
from the category of cubical sets to the category of permutahedral
sets; Section 3 introduces the permutocubes; in Section 4 we
introduce the notion of a permutocubical set; Section 5 introduces
the notion of a truncating twisting function and the resulting
twisted Cartesian product; in Section 6 we define an explicit
diagonal on the permutocubes; in Section 7 we build the
permutocubical set model for the double path fibration;
 in Section 8
a permutocubical model and the corresponding multiplicative
twisted tensor product for a fibration are constructed,
 and,
finally, in Section 9  the twisted tensor product theory for
Hirsch algebras is developed.

\section{The permutahedral set functor ${\bf \Omega }Q$ } \label{secpermut}

For completeness we first recall some basic facts about
permutahedral sets from \cite{SU2} (compare, \cite{KS2}).

\subsection{Permutahedral sets}\
This subsection introduces the notion of a permutahedral set
$\mathcal{Z}$, which is a combinatorial object generated by
permutahedra and equipped with appropriate face and degeneracy
operators. We construct the generating category $\mathbf{P\ }$and
show how to lift the diagonal on the permutahedra $P$ constructed
above to a diagonal on $\mathcal{Z}$. Naturally occurring examples
of permutahedral sets include the double cobar construction, i.e.,
Adams' cobar construction \cite{Adams} on the cobar with
coassociative coproduct \cite{Baues1}, \cite{CM}, \cite{KS1}.
Permutahedral sets are distinguished from simplicial or cubical
sets by their higher order structure relations. While our
construction of $\mathbf{P\ }$follows the analogous (but not
equivalent) construction for polyhedral sets given by D.W. Jones
in \cite{Jones}, there is no mention of structure relations in
\cite{Jones}.

Let $S_{n}$ be the symmetric group on $\underline{n}=\left\{
1,2,\ldots ,n\right\}  .$ Recall that the permutahedron $P_{n}$ is
the convex hull of $n!$ vertices $\left(
\sigma(1),\ldots,\sigma(n)\right)  \in\mathbb{R}^{n},$ $\sigma\in
S_{n}$ \cite{Coxeter}, \cite{Milgram}.$\ $As a cellular complex,
$P_{n}$ is an $\left(  n-1\right)  $-dimensional convex polytope
whose $\left(  n-p\right)  $-faces are indexed by (ordered)
partitions $U_{1}|\cdots|U_{p}$ of $\underline{n}$. We shall
define the permutahedra inductively as subdivisions of the
standard $n$-cube $I^{n}.$ With this representation the
combinatorial connection between faces and partitions is
immediately clear.

Assign the label $\underline{1}$ to the single point $P_{1}.$ If
$P_{n-1}$ has been constructed and $u=U_{1}|\cdots|U_{p}$ is one
of its faces, form the sequence $u_{\ast}=\left\{
u_{0}=0,u_{1},\ldots,u_{p-1},u_{p}=\infty\right\} $ where
$u_{j}=\#\left(  U_{p-j+1}\cup\cdots\cup U_{p}\right)  ,$ $1\leq
j\leq
p-1$ and $\#$ denotes cardinality. Define the \emph{subdivision of }%
$I$\emph{\ relative to }$u$ to be
\[
I/u_{\ast}=I_{1}\cup I_{2}\cup\cdots\cup I_{p},
\]
where $I_{j}=\left[
1-\frac{1}{2^{u_{j-1}}},1-\frac{1}{2^{u_{j}}}\right]  $ and
$\frac{1}{2^{\infty}}=0.$ Then
\[
P_{n}=\bigcup\limits_{u\in P_{n-1}}u\times I/u_{\ast}%
\]
with faces labeled as follows (see Figures 1 and 2)\vspace{0.1in}:%
\[%
\begin{tabular}
[c]{c|cc}%
\textbf{Face of }$\underset{\ }{u\times I/u_{\ast}}$ &
\textbf{Partition of }$\underline{n}$ & \\\hline
&  & \\
$u\times0$ & $U_{1}|\cdots|U_{p}|n$ & \\
&  & \\
$u\times(I_{j}\cap I_{j+1})$ &
$U_{1}|\cdots|U_{p-j}|n|U_{p-j+1}|\cdots
|U_{p},$ & $1\leq j\leq p-1$\\
&  & \\
$u\times1$ & $n|U_{1}|\cdots|U_{p}$ & \\
&  & \\
$u\times I_{j}$ & $U_{1}|\cdots|U_{p-j+1}\cup n|\cdots|U_{p}.$ &
\end{tabular}
\
\]

\vspace{0.1in}

A \emph{cubical }vertex of $P_{n}$ is a vertex common to both
$P_{n}$ and $I^{n-1}.$ Note that $u$ is a cubical\ vertex of
$P_{n-1}$ if and only if $u|n$ and $n|u$ are cubical\ vertices of
$P_{n}.$ Thus the cubical vertices of $P_{3}$ are $1|2|3,$
$2|1|3,$ $3|1|2$ and $3|2|1$ since $1|2$ and $2|1$ are cubical
vertices of $P_{2}.$ \vspace{0.1in}

\begin{center}
\setlength{\unitlength}{0.0004in}\begin{picture}
(2975,2685)(3126,-2038) \thicklines \put(3601,239){\line(
1,0){1800}} \put(5401,239){\line( 0,-1){1800}}
\put(5401,-1561){\line(-1, 0){1800}} \put(3601,-1561){\line(
0,1){1800}} \put(3601,239){\makebox(0,0){$\bullet$}}
\put(3601,-661){\makebox(0,0){$\bullet$}}
\put(3601,-1561){\makebox(0,0){$\bullet$}}
\put(5401,239){\makebox(0,0){$\bullet$}}
\put(5401,-661){\makebox(0,0){$\bullet$}}
\put(5401,-1561){\makebox(0,0){$\bullet$}}
\put(4500,-680){\makebox(0,0){$123$}}
\put(2980,-1861){\makebox(0,0){$1|2|3$}}
\put(2980,-699){\makebox(0,0){$1|3|2$}}
\put(2980,464){\makebox(0,0){$3|1|2$}}
\put(6000,-1861){\makebox(0,0){$2|1|3$}}
\put(6000,-699){\makebox(0,0){$2|3|1$}}
\put(6000,464){\makebox(0,0){$3|2|1$}}
\put(3040,-1260){\makebox(0,0){$1|23$}}
\put(4550,530){\makebox(0,0){$3|12$}}
\put(3040,-111){\makebox(0,0){$13|2$}}
\put(5960,-111){\makebox(0,0){$23|1$}}
\put(5960,-1260){\makebox(0,0){$2|13$}}
\put(4550,-1890){\makebox(0,0){$12|3$}}
\end{picture}\vspace{0.1in}

Figure 1: $P_{3}$ as a subdivision of $P_{2}\times I$.\
\vspace{0.4in}

\setlength{\unitlength}{0.00023in}\begin{picture} (7500,7500)
\thicklines
\put(3000,4800){\line( 0,-1){4800}}
\put(3000,4800){\makebox(0,0){$\bullet$}}
\put(3000,2400){\makebox(0,0){$\bullet$}} \put(3000,0){\makebox
(0,0){$\bullet$}} \put(3000,3600){\makebox(0,0){$\bullet$}}
\put(3000,0){\line( 1, 0){4800}}
\put(7800,0){\line( 0, 1){4800}}
\put(7800,4800){\line(-1, 0){4800}}
\put(7800,4800){\makebox(0,0){$\bullet$}}
\put(7800,3600){\makebox(0,0){$\bullet$}} \put(7800,2400){\makebox
(0,0){$\bullet$}} \put(7800,0){\makebox(0,0){$\bullet$}}
\put(3000,2400){\line( 1, 0){4800}}
\put(0,6800){\line( 0,-1){4800}}
\put(0,6800){\makebox(0,0){$\bullet$}}
\put(0,5600){\makebox(0,0){$\bullet$}}
\put(0,4400){\makebox(0,0){$\bullet$}}
\put(0,2000){\makebox(0,0){$\bullet$}}
\put(0,6800){\line( 1, 0){4800}}
\put(0,5600){\line( 1, 0){1200}} \put(2000,5600){\line( 1,
0){2800}}
\put(0,2000){\line( 1, 0){1200}} \put(1800,2000){\line( 1,
0){900}} \put(3300,2000){\line( 1, 0){1500}}
\put(4800,5000){\line( 0,1){1700}} \put(4800,2600){\line(
0,1){2050}} \put(4800,2000){\line( 0,1){200}}
\put(4800,6800){\makebox(0,0){$\bullet$}}
\put(4800,5600){\makebox(0,0){$\bullet$}} \put(4800,4400){\makebox
(0,0){$\bullet$}} \put(4800,2000){\makebox(0,0){$\bullet$}}
\put(0,2000){\line( 3,-2){3000}}
\put(3000,4800){\line(-3,2){3000}}
\put(1500,5800){\line( 0,-1){4800}}
\put(1500,5800){\makebox(0,0){$\bullet$}}
\put(1500,4600){\makebox(0,0){$\bullet$}} \put(1500,3400){\makebox
(0,0){$\bullet$}} \put(1500,1000){\makebox(0,0){$\bullet$}}
\put(3000,3600){\line(-3, 2){1500}}
\put(1500,3400){\line(-3, 2){1500}}
\put(6300,3400){\line(-3, 2){1500}}
\put(6300,5800){\line( 0,-1){800}} \put(6300,4600){\line(
0,-1){2000}} \put(6300,2180){\line( 0,-1){1200}}
\put(6300,5800){\makebox(0,0){$\bullet$ }}
\put(6300,3400){\makebox(0,0){$\bullet$}} \put(6300,1000){\makebox
(0,0){$\bullet$}} \put(6300,4600){\makebox(0,0){$\bullet$}}
\put(7800,3600){\line(-3, 2){1500}}
\put(4800,2000){\line(3, -2){3000}}
\put(4800,6800){\line( 3,-2){3000}}
\put(-1000,1000){\makebox(0,0){$(0,1,0)$}}
\put(5100,7600){\makebox(0,0){$(1,1,1)$}}
\put(8750,-1000){\makebox(0,0){$(1,0,0)$}} \put
(2700,-1000){\makebox(0,0){$(0,0,0)$}}
\end{picture}\vspace{0.4in}

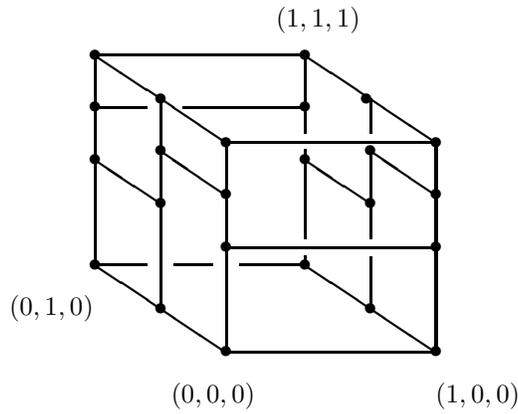
Figure 2a: $P_{4}$ as a subdivision of $P_{3}\times
I.\vspace{0.2in}$

\newpage

\setlength{\unitlength}{0.009in}\begin{picture} (500,-500)
\thicklines
\put(0,-120){\line( 0,-1){120}} \put(120,0){\line( 0,-1){360}}
\put(180,-120){\line( 0,-1){120}} \put(240,0){\line( 0,-1){360}}
\put(360,-120){\line( 0,-1){120}} \put(420,-120){\line(
0,-1){120}} \put(480,-120){\line( 0,-1){120}}
\put(120,0){\line( 1,0){120}} \put(0,-120){\line( 1,0){480}} \put
(120,-150){\line( 1,0){60}} \put(180,-180){\line( 1,0){60}} \put
(240,-150){\line( 1,0){120}} \put(420,-150){\line( 1,0){60}} \put
(0,-180){\line( 1,0){120}} \put(360,-180){\line( 1,0){60}}
\put(0,-240){\line( 1,0){480}} \put(120,-360){\line( 1,0){120}}
\put(120,0){\makebox(0,0){$\bullet$}}
\put(180,0){\makebox(0,0){$\bullet$}}
\put(240,0){\makebox(0,0){$\bullet$}}
\put(0,-120){\makebox(0,0){$\bullet$}}
\put(120,-120){\makebox(0,0){$\bullet$}} \put(180,-120){\makebox
(0,0){$\bullet$}} \put(240,-120){\makebox(0,0){$\bullet$}} \put
(360,-120){\makebox(0,0){$\bullet$}}
\put(420,-120){\makebox(0,0){$\bullet$}}
\put(480,-120){\makebox(0,0){$\bullet$}} \put(120,-150){\makebox
(0,0){$\bullet$}} \put(180,-150){\makebox(0,0){$\bullet$}} \put
(240,-150){\makebox(0,0){$\bullet$}}
\put(360,-150){\makebox(0,0){$\bullet$}}
\put(0,-180){\makebox(0,0){$\bullet$}}
\put(182.25,-180){\makebox(0,0){$\bullet$ }}
\put(240,-180){\makebox(0,0){$\bullet$}} \put(420,-180){\makebox
(0,0){$\bullet$}} \put(480,-180){\makebox(0,0){$\bullet$}} \put
(0,-150){\makebox(0,0){$\bullet$}}
\put(120,-240){\makebox(0,0){$\bullet$}}
\put(360,-240){\makebox(0,0){$\bullet$}} \put(420,-240){\makebox
(0,0){$\bullet$}} \put(480,-150){\makebox(0,0){$\bullet$}} \put
(0,-240){\makebox(0,0){$\bullet$}}
\put(120,-180){\makebox(0,0){$\bullet$}}
\put(180,-240){\makebox(0,0){$\bullet$}} \put(240,-240){\makebox
(0,0){$\bullet$}} \put(360,-180){\makebox(0,0){$\bullet$}} \put
(420,-150){\makebox(0,0){$\bullet$}}
\put(480,-240){\makebox(0,0){$\bullet$}}
\put(120,-360){\makebox(0,0){$\bullet$}} \put(180,-360){\makebox
(0,0){$\bullet$}} \put(240,-360){\makebox(0,0){$\bullet$}}
\put(275,-100){\makebox(0,0){$(1,1,1)$}}
\put(20,-260){\makebox(0,0){$(0,0,0)$ }}
\put(55,-152){\makebox(0,0){$124|3$}}
\put(55,-210){\makebox(0,0){$12|34$}}
\put(150,-134){\makebox(0,0){$24|13$}}
\put(150,-195){\makebox(0,0){$2|134$}}
\put(177,-55){\makebox(0,0){$4|123$}}
\put(177,-295){\makebox(0,0){$123|4$}}
\put(210,-152){\makebox(0,0){$234|1$}}
\put(210,-210){\makebox(0,0){$23|14$}}
\put(300,-134){\makebox(0,0){$34|12$}}
\put(300,-195){\makebox(0,0){$3|124$}}
\put(390,-152){\makebox(0,0){$134|2$}}
\put(390,-200){\makebox(0,0){$13|24$}}
\put(450,-134){\makebox(0,0){$14|23$}}
\put(450,-194){\makebox(0,0){$1|234$}}
\end{picture}\vspace*{3.4in}

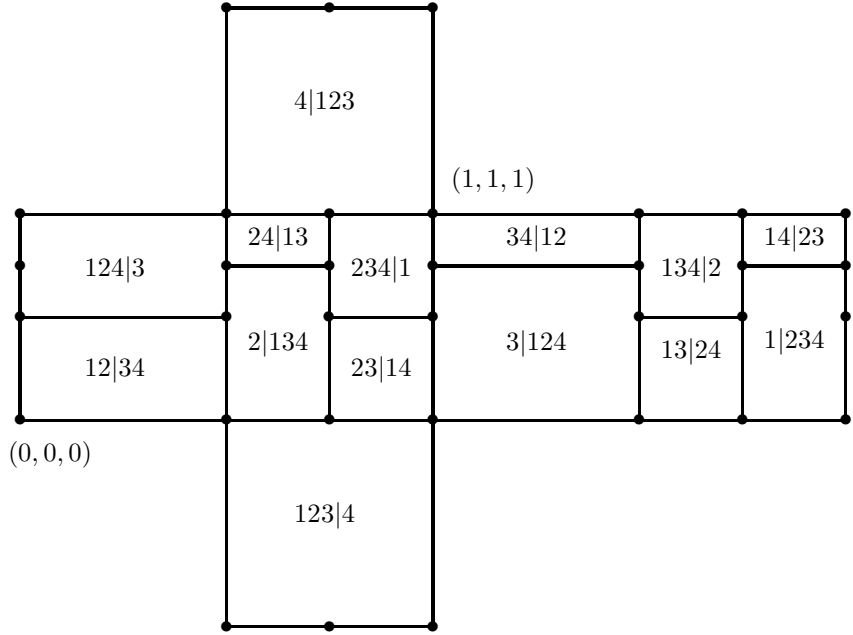
Figure 2b: The $2$-faces of $P_{4}.$
\end{center}


\subsection{Singular Permutahedral Sets}

By way of motivation we begin with constructions of two singular
permutahedral sets--our universal examples. Whereas the first
emphasizes coface and codegeneracy operators, the second
emphasizes cellular chains and is appropriate for homology theory.
We begin by constructing the various maps we need to define
singular coface and codegeneracy operators.

Fix a positive integer $n.$ For $0\leq p\leq n,$ let
\[
\underline{p}=\left\{
\begin{array}
[c]{cc}%
\varnothing, & p=0\\
\left\{  1,\ldots,p\right\}  , & 1\leq p\leq n
\end{array}
\right.  \text{ \ and \ }\overline{p}=\left\{
\begin{array}
[c]{cc}%
\varnothing, & p=0\\
\left\{  n-p+1,\ldots,n\right\}  , & 1\leq p\leq n;
\end{array}
\right.
\]
then $\underline{p}$ and $\overline{p}$ contain the first and last
$p$
elements of $\underline{n},$ respectively; note that $\underline{p}%
\cap\overline{q}=\left\{  p\right\}  $ whenever $p+q=n+1$. Given
integers $r,s\in\underline{n}$ such that $r+s=n+1,$ there is a
canonical projection $\Delta_{r,s}:P_{n}\rightarrow P_{r}\times
P_{s}$ whose restriction to a vertex $v=a_{1}|\cdots|a_{n}\in
P_{n}$ is given by
\[
\Delta_{r,s}(v)=b_{1}|\cdots|b_{r}\times c_{1}|\cdots|c_{s},
\]
where $\left(
b_{1},\ldots,b_{r};c_{1},\ldots,c_{k-1},c_{k+1},\ldots
,c_{s}\right)  \,$is the unshuffle of $\left(
a_{1},\ldots,a_{n}\right)  $ with $b_{i}\in\underline{r},$
$c_{j}\in\overline{s},$ $c_{k}=r.$ For example,
$\Delta_{2,3}(2|4|1|3)=2|1\times2|4|3$ and
$\Delta_{3,2}(2|4|1|3)=2|1|3\times 4|3.$ Since the image of the
vertices of a cell of $P_{n}$ uniquely determines a cell in
$P_{r}\times P_{s}$ the map $\Delta_{r,s}$ is well-defined and
cellular. Furthermore, the restriction of $\Delta_{r,s}$ to an
$(n-k)$-cell $A_{1}|\cdots|A_{k}\subset P_{n}$ is given by
\[
\Delta_{r,s}\left(  A_{1}|\cdots|A_{k}\right)  =\left\{
\begin{array}
[c]{ll}%
\underline{r}\times\left(
A_{1}|\cdots|A_{i}\setminus\underline{r-1}\text{
}|\cdots|A_{k}\right)  , & \hspace*{-0.3in}\text{if
}\underline{r}\subseteq
A_{i},\text{ some }i,\\
& \\
\left(  A_{1}|\cdots|A_{j}\setminus\overline{s-1}\text{
}|\cdots|A_{k}\right) \times\overline{s}, &
\hspace*{-0.3in}\text{if }\overline{s}\subseteq
A_{j},\text{ some }j,\\
& \\%
\begin{array}
[c]{l}%
\hspace*{-0.06in}\left(  A_{1}\setminus\overline{s-1}\text{
}|\cdots
|A_{k}\setminus\overline{s-1}\right) \\
\hspace*{0.5in}\times\left(  A_{1}\setminus\underline{r-1}\text{ }%
|\cdots|A_{k}\setminus\underline{r-1}\right)  ,
\end{array}
& \text{otherwise.}%
\end{array}
\right.
\]
Note that $\Delta_{r,s}$ acts homeomorphically in the first two
cases and degeneratively in the third when $1<k<n$. When $n=3$ for
example, $\Delta_{2,2}$ maps the edge $1|23$ onto the edge
$1|2\times23$ and the edge $13|2$ onto the vertex $1|2\times3|2$
(see Figure 3).

\hspace*{0.3in}\setlength{\unitlength}{0.0005in}\begin{picture}
(2975,1000) \thicklines \put(1,239){\line( 1,0){1800}}
\put(1801,239){\line( 0,-1){1800}} \put(1801,-1561){\line(-1,
0){1800}} \put(1,-1561){\line( 0,1){1800}}
\put(1,239){\makebox(0,0){$\bullet$}}
\put(1,-661){\makebox(0,0){$\bullet$}}
\put(1,-1561){\makebox(0,0){$\bullet$}}
\put(1801,239){\makebox(0,0){$\bullet$}}
\put(1801,-661){\makebox(0,0){$\bullet$}}
\put(1801,-1561){\makebox(0,0){$\bullet$}}
\put(900,-680){\makebox(0,0){$123$}}
\put(-490,-1260){\makebox(0,0){${1|23}$}}
\put(950,530){\makebox(0,0){${3|12}$}}
\put(-490,-111){\makebox(0,0){${13|2}$}}
\put(2300,-111){\makebox(0,0){${23|1}$}}
\put(2300,-1260){\makebox(0,0){${2|13}$}}
\put(950,-1890){\makebox(0,0){${12|3}$}}
\put(3200,-400){\makebox(0,0){$\Delta_{2,2}$}}
\put(2700,-661){\vector(1,0){1000}}
\put(4700,-661){\makebox(0,0){${1|2} \times 23$}}
\put(8325,-661){\makebox(0,0){${2|1} \times 23$}}
\put(6500,530){\makebox(0,0){$12 \times {3|2}$}}
\put(6500,-1890){\makebox(0,0){$12 \times {2|3}$}}
\put(5601,239){\line( 1,0){1800}} \put(7401,239){\line(
0,-1){1800}} \put(7401,-1561){\line(-1, 0){1800}}
\put(5601,-1561){\line( 0,1){1800}}
\put(5601,239){\makebox(0,0){$\bullet$}}
\put(5601,-1561){\makebox(0,0){$\bullet$}}
\put(7401,239){\makebox(0,0){$\bullet$}}
\put(7401,-1561){\makebox(0,0){$\bullet$}}
\put(6500,-661){\makebox(0,0){$12 \times 23$}}
\end{picture}$\vspace{1.1in}$

\begin{center}
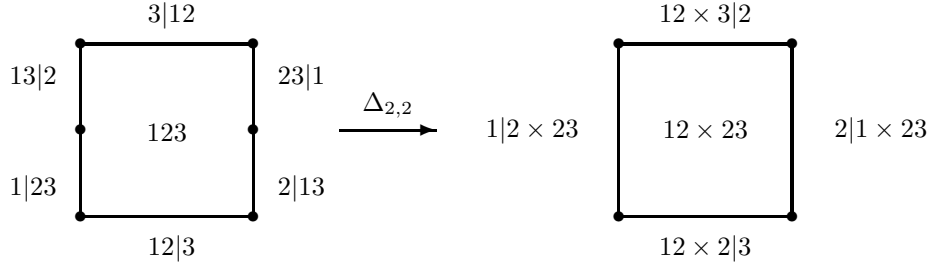
Figure 3: The projection $\Delta_{2,2}:P_{3}\rightarrow
I^{2}.\vspace{0.2in} $
\end{center}

Now identify the set $U=\left\{  u_{1}<\cdots<u_{n}\right\}  $
with $P_{n}$ and the ordered partitions of $U$ with the faces of
$P_{n}$ in the obvious way. Then $\left(
\Delta_{r,s}\times1\right)  \circ\Delta_{r+s-1,t}=\left(
1\times\Delta_{s,t}\right)  \circ\Delta_{r,s+t-1}$ whenever
$r+s+t=n+2$ so that $\Delta_{\ast,\ast}$ acts coassociatively with
respect to Cartesian product. It follows that each $k$-tuple
$\left(  n_{1},\ldots,n_{k}\right) \in\mathbb{N}^{k}$ with
$k\geq2$ and $n_{1}+\cdots+n_{k}=n+k-1$ uniquely determines a
cellular projection $\Delta_{n_{1}\cdots n_{k}}:P_{n}\rightarrow
P_{n_{1}}\times\cdots\times P_{n_{k}} $ given by the composition%
\[
\Delta_{n_{1}\cdots n_{k}}=\left(
\Delta_{n_{1},n_{2}}\times1^{\times k-2}\right)
\circ\cdots\circ\left(  \Delta_{n_{\left(  k-2\right)
}-k+3,n_{k-1}}\times1\right)  \circ\Delta_{n_{\left(  k-1\right)  }-k+2,n_{k}%
},
\]
where $n_{\left(  q\right)  }=n_{1}+\cdots+n_{q};$ and in
particular,
\begin{equation}
\Delta_{n_{1}\cdots n_{k}}\left(  \underline{n}\right)  =\underline{n_{1}%
}\times\underline{n_{\left(  2\right)  }-1}\setminus\underline{n_{1}-1}%
\times\cdots\times\underline{n_{\left(  k\right)  }-\left(
k-1\right) }\setminus\underline{n_{\left(  k-1\right)  }-\left(
k-1\right)  }.
\label{proj}%
\end{equation}
Note that formula \ref{proj} with $k=n-1$ and $n_{i}=2$ for all
$i$ defines a
projection $\rho_{n}:P_{n}\rightarrow I^{n-1}$%
\[
\rho_{n}\left(  \underline{n}\right)  =\Delta_{2\cdots2}\left(
\underline {n}\right)  =12\times23\times\cdots\times\left\{
n-1,n\right\}
\]
(see Figure 4) acting on a vertex $u=u_{1}|\cdots|u_{n}$ as
follows: For each $i\in\underline{n-1},$ let $\left\{
u_{j},u_{k}\text{ }|\text{ }j<k\right\} =\left\{
u_{1},\ldots,u_{n}\right\}  \cap\left\{  i,i+1\right\}  $ and set
$v_{i}=u_{j},$ $v_{i+1}=u_{k};$ then
$\rho_{n}(u)=v_{1}|v_{2}\times \cdots\times
v_{n-1}|v_{n}.$\vspace{0.2in}

\newpage

\hspace*{0.4in}\setlength{\unitlength}{0.00015in}\begin{picture}
(7500,7500) \thicklines
\put(6000,4800){\line( 0,-1){4800}}
\put(6000,4800){\makebox(0,0){$\bullet$}}
\put(6000,2400){\makebox(0,0){$\bullet$}}
\put(6000,0){\makebox(0,0){$\bullet$}}
\put(6000,3600){\makebox(0,0){$\bullet$}}
\put(6000,0){\line( 1, 0){4800}}
\put(10800,0){\line( 0, 1){4800}}
\put(10800,4800){\line(-1, 0){4800}}
\put(10800,4800){\makebox(0,0){$\bullet$}}
\put(10800,3600){\makebox(0,0){$\bullet$}}
\put(10800,2400){\makebox(0,0){$\bullet$}}
\put(10800,0){\makebox(0,0){$\bullet$}}
\put(6000,2400){\line( 1, 0){4800}}
\put(3000,6800){\line( 0,-1){4800}}
\put(3000,6800){\makebox(0,0){$\bullet$}}
\put(3000,5600){\makebox(0,0){$\bullet$}}
\put(3000,4400){\makebox(0,0){$\bullet$}}
\put(3000,2000){\makebox(0,0){$\bullet$}}
\put(3000,6800){\line( 1, 0){4800}}
\put(3000,5600){\line( 1, 0){1200}}
\put(5000,5600){\line(1,0){2800}}
\put(3000,2000){\line( 1, 0){1200}}
\put(4800,2000){\line(1,0){900}} \put(6300,2000){\line( 1,
0){1500}}
\put(7800,5000){\line( 0,1){1700}}
\put(7800,2600){\line(0,1){2050}} \put(7800,2000){\line(0,1){200}}
\put(7800,6800){\makebox(0,0){$\bullet$}}
\put(7800,5600){\makebox(0,0){$\bullet$}}
\put(7800,4400){\makebox(0,0){$\bullet$}}
\put(7800,2000){\makebox(0,0){$\bullet$}}
\put(3000,2000){\line( 3,-2){3000}}
\put(6000,4800){\line(-3,2){3000}}
\put(4500,5800){\line( 0,-1){4800}}
\put(4500,5800){\makebox(0,0){$\bullet$}}
\put(4500,4600){\makebox(0,0){$\bullet$}}
\put(4500,3400){\makebox(0,0){$\bullet$}}
\put(4500,1000){\makebox(0,0){$\bullet$}}
\put(6000,3600){\line(-3, 2){1500}}
\put(4500,3400){\line(-3, 2){1500}}
\put(9300,3400){\line(-3, 2){1500}}
\put(9300,5800){\line( 0,-1){800}}
\put(9300,4600){\line(0,-1){2000}}
\put(9300,2180){\line(0,-1){1200}}
\put(9370,5800){\makebox(0,0){$\bullet$ }}
\put(9300,3400){\makebox(0,0){$\bullet$}}
\put(9300,1000){\makebox(0,0){$\bullet$}}
\put(9300,4600){\makebox(0,0){$\bullet$}}
\put(10800,3600){\line(-3, 2){1500}}
\put(7800,2000){\line(3, -2){3000}}
\put(7800,6800){\line( 3,-2){3000}}
\put(17000,4800){\line( 0,-1){4800}}
\put(17000,4800){\makebox(0,0){$\bullet$}}
\put(17000,2400){\makebox(0,0){$\bullet$}}
\put(17000,0){\makebox(0,0){$\bullet$}}
\put(17000,0){\line( 1, 0){4800}}
\put(21800,0){\line( 0, 1){4800}}
\put(21800,4800){\line(-1, 0){4800}}
\put(21800,4800){\makebox(0,0){$\bullet$}}
\put(21800,2400){\makebox(0,0){$\bullet$}}
\put(21800,0){\makebox(0,0){$\bullet$}}
\put(17000,2400){\line( 1, 0){4800}}
\put(14000,6800){\line( 0,-1){4800}}
\put(14000,6800){\makebox(0,0){$\bullet$}}
\put(14000,4400){\makebox(0,0){$\bullet$}}
\put(14000,2000){\makebox(0,0){$\bullet$}}
\put(14000,6800){\line( 1, 0){4800}}
\put(14000,4400){\line( 1, 0){2600}} \put(17300,4400){\line(
1,0){1500}}
\put(14000,2000){\line( 1, 0){2600}} \put(17300,2000){\line(
1,0){1500}}
\put(18800,5000){\line( 0,1){1700}}
\put(18800,2600){\line(0,1){2050}} \put(18800,2000){\line(
0,1){200}} \put(18800,6800){\makebox(0,0){$\bullet$}}
\put(18800,4400){\makebox(0,0){$\bullet$}}
\put(18800,2000){\makebox(0,0){$\bullet$}}
\put(14000,2000){\line( 3,-2){3000}}
\put(17000,4800){\line(-3,2){3000}}
\put(18800,2000){\line(3, -2){3000}}
\put(18800,6800){\line( 3,-2){3000}}
\put(6000,-11200){\line( 0,-1){4800}}
\put(6000,-11200){\makebox(0,0){$\bullet$}}
\put(6000,-16000){\makebox(0,0){$\bullet$}}
\put(6000,-16000){\line( 1, 0){4800}}
\put(10800,-16000){\line( 0, 1){4800}}
\put(10800,-11200){\line(-1, 0){4800}}
\put(10800,-11200){\makebox(0,0){$\bullet$}}
\put(10800,-16000){\makebox(0,0){$\bullet$}}
\put(3000,-9200){\line( 0,-1){4800}}
\put(3000,-9200){\makebox(0,0){$\bullet$}}
\put(3000,-14000){\makebox(0,0){$\bullet$}}
\put(3000,-9200){\line( 1, 0){4800}}
\put(3000,-14000){\line( 1, 0){1200}}
\put(4800,-14000){\line(1,0){900}} \put(6300,-14000){\line( 1,
0){1500}}
\put(7800,-11000){\line( 0,1){1700}}
\put(7800,-14000){\line(0,1){2600}}
\put(7800,-9200){\makebox(0,0){$\bullet$}}
\put(7800,-14000){\makebox(0,0){$\bullet$}}
\put(3000,-14000){\line( 3,-2){3000}}
\put(6000,-11200){\line(-3,2){3000}}
\put(4500,-10200){\line( 0,-1){4800}}
\put(4500,-10200){\makebox(0,0){$\bullet$}}
\put(4500,-15050){\makebox(0,0){$\bullet$}}
\put(9300,-10200){\line( 0,-1){800}}
\put(9300,-11400){\line(0,-1){3600}}
\put(9370,-10200){\makebox(0,0){$\bullet$ }}
\put(9300,-15020){\makebox(0,0){$\bullet$}}
\put(7800,-14000){\line(3, -2){3000}}
\put(7800,-9200){\line( 3,-2){3000}}
\put(17000,-11200){\line( 0,-1){4800}}
\put(17000,-11200){\makebox(0,0){$\bullet$}}
\put(17000,-16000){\makebox(0,0){$\bullet$}}
\put(17000,-16000){\line( 1, 0){4800}}
\put(21800,-16000){\line( 0, 1){4800}}
\put(21800,-11200){\line(-1, 0){4800}}
\put(21800,-11200){\makebox(0,0){$\bullet$}}
\put(21800,-16000){\makebox(0,0){$\bullet$}}
\put(14000,-9200){\line( 0,-1){4800}}
\put(14000,-9200){\makebox(0,0){$\bullet$}}
\put(14000,-14000){\makebox(0,0){$\bullet$}}
\put(14000,-9200){\line( 1, 0){4800}}
\put(14000,-14000){\line( 1, 0){2600}} \put(17300,-14000){\line(
1,0){1500}}
\put(18800,-11000){\line( 0,1){1700}}
\put(18800,-14000){\line(0,1){2650}}
\put(18800,-9200){\makebox(0,0){$\bullet$}}
\put(18800,-14000){\makebox(0,0){$\bullet$}}
\put(14000,-14000){\line( 3,-2){3000}}
\put(14000,-9200){\line( 3,-2){3000}}
\put(18800,-14000){\line(3, -2){3000}}
\put(18800,-9200){\line( 3,-2){3000}}
\put(5750,-4300){\makebox(0,0){$\Delta_{3,2}$}}
\put(7700,-3400){\vector(0,-1){2000}}
\put(7700,-2700){\makebox(0,0){$1234$}}
\put(7700,-6300){\makebox(0,0){$123\times34$}}
\put(9500,-2700){\vector(1,0){4000}}
\put(9900,-6300){\vector(1,0){3200}}
\put(16000,-2700){\makebox(0,0){$12\times234$}}
\put(16000,-6300){\makebox(0,0){$12\times23\times34$}}
\put(16000,-3400){\vector(0,-1){2000}}
\put(18500,-4300){\makebox(0,0){$1\times\Delta_{2,2}$}}
\put(11500,-7600){\makebox(0,0){$\Delta_{2,2}\times1$}}
\put(11500,-1500){\makebox(0,0){$\Delta_{2,3}$}}
\end{picture}$\vspace{2.6in}$

\begin{center}
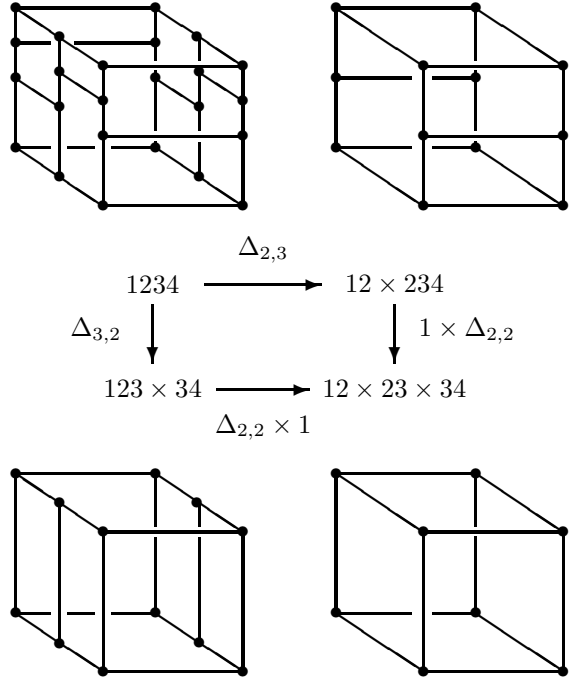
Figure 4: The projection $\rho_{4}:P_{4}\rightarrow
I^{3}.\vspace{0.2in}$
\end{center}

Now choose a (non-cellular) homeomorphism
$\gamma_{n}:I^{n-1}\rightarrow P_{n}$ whose restriction to a
vertex $v=v_{1}|v_{2}\times\cdots\times
v_{n-1}|v_{n}$ can be expressed inductively as follows: Set $A_{2}=v_{1}%
|v_{2};$ if $A_{k-1}$ has been obtained from
$v_{1}|v_{2}\times\cdots\times
v_{k-2}|v_{k-1},$ set%
\[
A_{k}=\left\{
\begin{array}
[c]{ll}%
A_{k-1}|k, & \text{if }v_{k}=k,\\
k|A_{k-1}, & \text{otherwise.}%
\end{array}
\right.
\]
For example, $\gamma_{4}\left(  2|1\times3|2\times3|4\right)
=3|2|1|4.$ Then $\gamma_{n}$ sends the vertices of $I^{n-1}$ to
cubical vertices of $P_{n}$ and the vertices of $P_{n}$ fixed by
$\gamma_{n}\rho_{n}$ are exactly its cubical vertices. Given a
codimension 1 face $A|B\subset P_{n},$ index the elements of $A$
and $B$ as follows: If $n\in A,$ write $A=\left\{
a_{1}<\cdots<a_{m}\right\}  $ and $B=\left\{
b_{1}<\cdots<b_{\ell}\right\} ;$ if $n\in B,$ write $A=\left\{
a_{1}<\cdots<a_{\ell}\right\}  $ and $B=\left\{
b_{1}<\cdots<b_{m}\right\}  .$ Then $A|B$ uniquely embeds in
$P_{n}$ as the subcomplex%
\[
P_{\ell}\times P_{m}=\left\{
\begin{array}
[c]{ll}%
a_{1}|\cdots|a_{m}|B\times A|b_{1}|\cdots|b_{\ell}, & \text{if }n\in A\\
A|b_{1}|\cdots|b_{m}\times a_{1}|\cdots|a_{\ell}|B, & \text{if
}n\in B.
\end{array}
\right.
\]
For example, $14|23$ embeds in $P_{4}$ as $1|4|23\times14|2|3$.
Let $\iota_{A|B}:A|B\hookrightarrow P_{\ell}\times P_{m}$ denote
this embedding and let $h_{A|B}=\iota_{A|B}^{-1};$ then
$h_{A|B}:P_{\ell}\times P_{m}\rightarrow A|B$ is an orientation
preserving homeomorphism. Also define the cellular projection
\[
\phi_{A|B}:P_{n}\rightarrow P_{\ell}\times P_{m}=\left\{
\begin{array}
[c]{ll}%
b_{1}\cdots b_{\ell}\times a_{1}\cdots a_{m}, & \text{if }n\in A\\
a_{1}\cdots a_{\ell}\times b_{1}\cdots b_{m}, & \text{if }n\in B
\end{array}
\right.
\]
on a vertex $c=c_{1}|\cdots|c_{n}$ by $\phi_{A|B}\left(  c\right)
=u_{1}|\cdots|u_{\ell}\times v_{1}|\cdots|v_{m},$ where
$(u_{1},\ldots
,u_{\ell};$ $v_{1},\ldots,v_{m})\,$is the unshuffle of $\left(  c_{1}%
,\ldots,c_{n}\right)  $ with $u_{i}\in B,$ $v_{j}\in A$ when $n\in
A$ or with $u_{i}\in A,$ $v_{j}\in B$ when $n\in B.$ Note that
unlike $\Delta_{r,s},$ the projection $\phi_{A|B}$ always
degenerates on the top cell; furthermore, $\phi_{A|B}\circ
h_{A|B}=\phi_{B|A}\circ h_{A|B}=1$. We note that when $A$ or $B$
is a singleton set, the projection $\phi_{A|B}$ was defined by
R.J. Milgram in \cite{Milgram}.

The \textit{singular codegeneracy operator associated with }$A|B$
is the map $\beta_{A|B}:P_{n}\rightarrow P_{n-1}$ given by the
composition
\[
P_{n}\overset{\phi_{A|B}}{\longrightarrow}P_{\ell}\times
P_{m}\overset
{\rho_{\ell}\times\rho_{m}}{\longrightarrow}I^{\ell-1}\times I^{m-1}%
=I^{n-2}\overset{\gamma_{n-1}}{\longrightarrow}P_{n-1};
\]
the \textit{singular coface operator associated with }$A|B$ is the
map
$\delta_{A|B}:P_{n-1}\rightarrow P_{n}$ given by the composition%
\[
P_{n-1}\overset{\rho_{n-1}}{\longrightarrow}I^{n-2}=I^{\ell-1}\times
I^{m-1}\overset{\gamma_{\ell}\times\gamma_{m}}{\longrightarrow}P_{\ell}\times
P_{m}\overset{h_{A|B}}{\longrightarrow}A|B\overset{i}{\hookrightarrow}P_{n}.
\]
Unlike the simplicial or cubical case, $\delta_{A|B}$ need not be
injective.
We shall often abuse notation and write $h_{A|B}:P_{\ell}\times P_{m}%
\rightarrow P_{n}$ when we mean $i\circ h_{A|B}.$

We are ready to define our first universal example. For future
reference and to emphasize the fact that our definition depends
only on positive integers, let $\left(  n_{1},\ldots,n_{k}\right)
\in\mathbb{N}^{k}$ such that
$n_{\left(  k\right)  }=n$ and denote%
\[
{\mathcal{P}}_{n_{1}\cdots n_{k}}\left(  n\right)
=\{\text{\textit{Partitions} }A_{1}|\cdots|A_{k}\text{ \textit{of }}%
\underline{n}\text{ }|\text{ }\#A_{i}=n_{i}\}.
\]

\begin{definition}
Let $Y$ be a topological space. The \underline{singular
permutahedral set of $Y$} consists of the singular set
\[
{\Sing}_{\ast}^{P}Y=\bigcup\limits_{n\geq1}\left[
{\Sing}_{n}^{P}Y=\left\{ \text{Continuous maps
}P_{n}{\rightarrow}Y\right\}  \right]
\]
together with singular face and degeneracy operators
\[
d_{A|B}:\Sing_{n}^{P}Y\rightarrow \Sing_{n-1}^{P}Y\text{ \ and \
}\varrho _{A|B}:\Sing_{n-1}^{P}Y\rightarrow \Sing_{n}^{P}Y
\]
defined respectively for each $n\geq2$ and
$A|B\in{\mathcal{P}}_{\ast\ast }\left(  n\right)  $ as the
pullback along $\delta_{A|B}$ and $\beta_{A|B}$, i.e., for $f\in
\Sing_{n}^{P}Y$ and $g\in \Sing_{n-1}^{P}Y,$
\[
d_{A|B}(f)=f\circ\delta_{A|B}\text{ \ and \
}\varrho_{A|B}(g)=g\circ \beta_{A|B}.
\]
\vspace*{0.1in}
\end{definition}

\hspace{1.3in}\setlength{\unitlength}{0.0002in}\begin{picture}
(0,0) \thicklines
\put(1500,-600){\vector(3,-1){8000}}
\put(11200,-1000){\vector(0,-1){1500}}
\put(4600,0){\makebox(0,0){$\delta_{A|B}:P_{n-1} \rightarrow
I^{n-1} \rightarrow P_{\ell} \times P_m \rightarrow A|B
\hookrightarrow P_n$}} \put(11200,-3500){\makebox(0,0){$Y$}}
\put(12000,-1500){\makebox(0,0){$f$}}
\put(4000,-2800){\makebox(0,0){$d_{A|B}(f)$}}
\end{picture}\vspace*{0.8in}

\begin{center}
Figure 5:\ The singular face operator associated with
$A|B$.\vspace{0.2in}
\end{center}

Although coface operators $\delta_{A|B}:P_{n-1}\rightarrow P_{n}$
need not be inclusions, the top cell of $P_{n-1}$ is always
non-degenerate; however, the top cell of $P_{n-2}$ may degenerate
under quadratic compositions
$\delta_{A|B}\delta_{C|D}:P_{n-2}\rightarrow P_{n}$ . For example,
$\delta_{12|34}\delta_{13|2}:P_{2}\rightarrow P_{4}$ is a constant
map, since $\delta_{12|34}:P_{3}\rightarrow P_{2}\times
P_{2}\hookrightarrow P_{4}$\ sends the edge\thinspace$13|2$ to the
vertex $1|2\times3|2$.

\newpage

\begin{definition}
\label{admissible} A quadratic composition of face operators
$d_{C|D}d_{A|B}$ \underline{acts on $P_{n}$} if the top cell of
$P_{n-2}$ is non-degenerate under the composition
 \[
\delta_{A|B}\delta_{C|D}:P_{n-2}\rightarrow P_{n}. \]

\end{definition}
 For comparison,
quadratic compositions of simplicial or cubical face operators
always act on the simplex or cube. When $d_{C|D}d_{A|B}$ acts on
$P_{n},$ we assign the label $d_{C|D}d_{A|B}$ to the codimension 2
face $\delta_{A|B}\delta_{C|D}\left(  \underline{n}\right) $. The
various paths of descent from the top cell to a cell in
codimension 2 gives rise to relations among compositions of face
and degeneracy operators (see Figure 6). \vspace{0.2in}

\setlength{\unitlength}{0.0005in}\begin{picture} (2975,1000)
\thicklines\put(3601,239){\line( 1, 0){1800}}
\put(5401,239){\line( 0,-1){1800}} \put(5401,-1561){\line(-1,
0){1800}} \put(3601,-1561){\line( 0, 1){1800}}
\put(3601,239){\makebox(0,0){$\bullet$}}
\put(3601,-661){\makebox(0,0){$\bullet$}}
\put(3601,-1561){\makebox(0,0){$\bullet$}}
\put(5401,239){\makebox(0,0){$\bullet$}} \put
(5401,-661){\makebox(0,0){$\bullet$}}
\put(5401,-1561){\makebox(0,0){$\bullet$}}
\put(4500,-680){\makebox(0,0){$123$}}
\put(2000,-1861){\makebox(0,0){$d_{1|2}d_{12|3}=d_{1|2}d_{1|23}$}}
\put(2000,-699){\makebox(0,0){$d_{1|2}d_{13|2}=d_{2|1}d_{1|23}$}}
\put(2000,464){\makebox(0,0){$d_{2|1}d_{13|2}=d_{1|2}d_{3|12}$}}
\put(7000,-1861){\makebox(0,0){$d_{1|2}d_{2|13}=d_{2|1}d_{12|3}$}}
\put(7000,-699){\makebox(0,0){$d_{2|1}d_{2|13}=d_{1|2}d_{23|1}$}}
\put(7000,464){\makebox(0,0){$d_{2|1}d_{23|1}=d_{2|1}d_{3|12}$}}
\put(3110,-1260){\makebox(0,0){$d_{1|23}$}}
\put(4550,530){\makebox(0,0){$d_{3|12}$}}
\put(3110,-111){\makebox(0,0){$d_{13|2}$}}
\put(5900,-111){\makebox(0,0){$d_{23|1}$}}
\put(5900,-1260){\makebox(0,0){$d_{2|13}$}}
\put(4550,-1890){\makebox(0,0){$d_{12|3}$}}
\end{picture}\vspace{1.1in}

\begin{center}
Figure 6: Quadratic relations on the vertices of
$P_{3}.$\vspace{0.2in}
\end{center}

It is interesting to note that singular permutahedral sets have
higher order structure relations, an example of which appears
below in Figure 7 (see also (\ref{hrelation})). This distinguishes
permutahedral sets from simplicial or cubical sets in which
relations are strictly quadratic. Our second universal example,
called a \textquotedblleft singular multipermutahedral
set,\textquotedblright\ specifies a singular permutahedral set by
restricting to maps $f=\bar{f}\circ\Delta_{n_{1}\cdots n_{k}}$ for
some continuous $\bar{f}:P_{n_{1}}\times\cdots\times
P_{n_{k}}\rightarrow Y$. Face and degeneracy operators satisfy
those relations above in which $\Delta _{n_{1}\cdots n_{k}}$ plays
no essential role.

\vspace*{0.7in} \hspace*{0.09in}%
\setlength{\unitlength}{0.000236in}\begin{picture} (0,0)
\thicklines
\put(800,0){\vector(1,0){1500}} \put(4200,0){\vector(1,0){6000}}
\put(12200,0){\vector(1,0){6000}} \put(0,800){\line(0,1){1000}}
\put(19200,1750){\vector(0,-1){1000}}
\put(0,1750){\line(1,0){19225}}
\put(19200,-600){\vector(0,-1){4000}}
\put(18400,-5500){\vector(-1,0){3000}}
\put(18800,-600){\vector(-1,-1){4000}}
\put(11400,-600){\vector(2,-3){2650}}
\put(4000,-500){\vector(2,-1){9200}}
\put(0,-5500){\vector(1,0){13200}} \put(0,-800){\line(0,-1){4720}}
\put(0,0){\makebox(0,0){$P_2$}} \put(3200,0){\makebox(0,0){$P_1$}}
\put(7400,900){\makebox(0,0){$\delta_{1|2}$}}
\put(11100,0){\makebox(0,0){$P_2$}}
\put(15300,900){\makebox(0,0){$\delta_{1|23}$}}
\put(19200,0){\makebox(0,0){$P_3$}}
\put(9621,2600){\makebox(0,0){$\delta_{13|2}$}}
\put(19200,0){\makebox(0,0){$P_3$}}
\put(18300,-3100){\makebox(0,0){$\delta_{12|34}$}}
\put(19200,-5500){\makebox(0,0){$P_4$}}
\put(14300,-5500){\makebox(0,0){$Y$}}
\put(1500,900){\makebox(0,0){$\beta_{2|1}$}}
\put(17000,-4800){\makebox(0,0){$f$}}
\put(10000,-1700){\makebox(0,0){$d_{1|23}d_{12|34}(f)$}}
\put(4800,-2700){\makebox(0,0){$d_{2|1}d_{1|23}d_{12|34}(f)$}}
\put(16000,-1700){\makebox(0,0){$d_{12|34}(f)$}}
\put(5500,-4800){\makebox(0,0){$\varrho_{2|1}d_{2|1}d_{1|23}d_{12|34}(f)=d_{13|2}d_{12|34}(f)$}}
\end{picture}\vspace*{1.4in}

\begin{center}
Figure 7: A quartic relation in
${\Sing}_{\ast}^{P}Y$.\vspace*{0.2in}
\end{center}

Once again, fix a positive integer $n,$ but this time consider
$\left( n_{1},\ldots,n_{k}\right)  \in\left(
\mathbb{N}\cup0\right)  ^{k}$ with
$n_{(k)}=n-1$ and the projection $\Delta_{n_{1}+1\cdots n_{k}+1}%
:P_{n}\rightarrow P_{n_{1}+1}\times\cdots\times P_{n_{k}+1}$ with
$\Delta _{n}:P_{n}\rightarrow P_{n}$ defined to be the identity.
Given a topological
space $Y,$ let%
\[
{\Sing}^{n_{1}\cdots n_{k}}Y=\left\{
\bar{f}\circ\Delta_{n_{1}+1\cdots n_{k}+1}:P_{n}\rightarrow
Y\text{ }|\text{\ }\bar{f}\ \text{\textit{is continuous}}\right\}
;
\]
define $f,f^{\prime}\in{\Sing}^{n_{1}\cdots n_{k}}Y$ to be
equivalent if there exists $g:P_{n_{1}+1}\times\cdots\times
P_{n_{i-1}+1}\times P_{1}\times P_{n_{i+1}+1}\times\cdots\times
P_{n_{k}+1}\rightarrow Y$ for some $i<k$ such that
\[
f=g\circ(1^{\times
i-1}\times\phi_{\underline{n_{i}+1}|n_{i}+1}\times1^{\times
k-i-1})\circ\Delta_{n_{1}+1\cdots n_{i-1}+1,n_{i}+2,n_{i+2}+1\cdots n_{k}+1}%
\]
and
\[
f^{\prime}=g\circ(1^{\times
i}\times\phi_{1|\underline{n_{i+2}+1}\setminus
1}\times1^{\times k-i-2})\circ\Delta_{n_{1}+1\cdots n_{i}+1,n_{i+2}%
+2,n_{i+3}+1\cdots n_{k}+1},
\]
in which case we write $f\sim f^{\prime}.$ The geometry of the
cube motivates this equivalence; the degeneracies in the product
of cubical sets implies the identification (c.f. \cite{Kan} or the
definition of the cubical set functor $\mathbf{\Omega}X$ in
\cite{KS1}).

Define the singular set
\[
{\Sing}_{n}^{M}Y=\bigcup_{\substack{\left(
n_{1},\ldots,n_{k}\right)
\in\left(  \mathbb{N}\cup0\right)  ^{k}\\n_{(k)}=n-1}}\left.  {\Sing}%
^{n_{1}\cdots n_{k}}Y\right/  \sim.
\]
Singular face and degeneracy operators
\[
d_{A|B}:\Sing_{n}^{M}Y\rightarrow \Sing_{n-1}^{M}Y\text{ \ and \
}\varrho _{A|B}:\Sing_{n-1}^{M}Y\rightarrow \Sing_{n}^{M}Y
\]
are defined piece-wise for each $n\geq2$ and
$A|B\in{\mathcal{P}}_{\ast,\ast }\left(  n\right)  ,$ depending on
the form of $A|B.$ More precisely, for each
pair of integers $\left(  p_{i},q_{i}\right)  ,$ $1\leq i\leq k,$ with%
\[
p_{i}=1+\sum_{j=1}^{i-1}n_{j}\text{ \ and \ }q_{i}=1+\sum_{j=i+1}^{k}%
n_{j},\text{ \ let}%
\]%
\[
\mathcal{Q}_{p_{i},q_{i}}\left(  n\right)  =\left\{  U|V\in{\mathcal{P}}%
_{\ast,\ast}\left(  n\right)  \text{ }|\text{ }\left(  \underline{p_{i}%
}\subseteq U\text{ or }\underline{p_{i}}\subseteq V\right)  \text{
and }\left(  \overline{q_{i}}\subseteq U\text{ or
}\overline{q_{i}}\subseteq
V\right)  \right\}  ;\vspace*{0.1in}%
\]
in particular, when $r+s=n+1,$ set $k=2,$ $p_{1}=q_{2}=1,$
$p_{2}=r$ and
$q_{1}=s,$ then%
\[
\mathcal{Q}_{r,1}\left(  n\right)  =\left\{
U|V\in{\mathcal{P}}_{\ast,\ast }\left(  n\right)  \text{ }|\text{
}\underline{r}\subseteq U\text{ or
}\underline{r}\subseteq V\right\}  \text{ and}%
\]%
\[
\mathcal{Q}_{1,s}\left(  n\right)  =\left\{
U|V\in{\mathcal{P}}_{\ast,\ast }\left(  n\right)  \text{ }|\text{
}\overline{s}\subseteq U\text{ or }\overline{s}\subseteq V\right\}
.
\]
Since we identify $\underline{r}|\overline{s}\subset P_{n+1}$ with
$P_{r}\times P_{s}=\Delta_{r,s}\left(  P_{n}\right)  ,$ it follows
that $A|B\in\mathcal{Q}_{p_{i},q_{i}}\left(  n\right)  $ for some
$i$ if and only if
$\delta_{A|B}\delta_{\underline{r}|\overline{s}}:P_{n-1}\rightarrow
P_{n+1}$ is non-degenerate; consequently we consider cases $A|B\in
\mathcal{Q}_{p_{i},q_{i}}\left(  n\right)  $ for some $i$ and
$A|B\notin \mathcal{Q}_{p_{i},q_{i}}\left(  n\right)  $ for all
$i.$

Since our definitions of $d_{A|B}$ and $\varrho_{A|B}$ are
independent in the first case and interdependent in the second, we
define both operators simultaneously. But first we need some
notation: Given an increasingly ordered set
$M=\{m_{1}<\cdots<m_{k}\}\subset{\mathbb{N}},$ let
$I_{M}:M\rightarrow \underline{\#M}$ denote the \emph{indexing
map} $m_{i}\mapsto i$ and let $M+z=\{m_{i}+z\}$ denote
\emph{translation by} $z\in\mathbb{Z}$. Of course, $M-z$ and $M+z$
are left and right translations when $z>0;$ we adopt the
convention that translation takes preference over set operations.

Assume $A|B\in{\mathcal{Q}}_{p_{i},q_{i}}(n)$ for some $i,$ and
let
\[
C_{i}=\left\{  p_{i},p_{i}+1,...,p_{i}+n_{i}\right\}  ;
\]%
\[
A_{i}=\left(  C_{i}\cap A\right)  -n_{(i-1)},\text{ \
}B_{i}=\left(  C_{i}\cap B\right)  -n_{(i-1)};
\]%
\begin{equation}
n_{i}^{\prime}=\#(A\cap C_{i})-1,\text{ \
}\,n_{i}^{\prime\prime}=\#(B\cap
C_{i})-1. \label{formal}%
\end{equation}
For example, $n=6,$ $n_{1}=3$ and $n_{2}=2$ determines the
projection
$\Delta_{4,3}:P_{6}\rightarrow1234\times456$ and pairs $\left(  p_{1}%
,q_{1}\right)  =\left(  1,3\right)  $ and $\left(
p_{2},q_{2}\right) =\left(  4,1\right)  .$ Thus
$A|B=1234|56\in{\mathcal{Q}}_{3,2}(6)$ and the composition
$\delta_{\underline{4}|\overline{3}}\delta_{A|B}:P_{5}\rightarrow
P_{7}$ is non-degenerate. Furthermore, $C_{2}=456,$ $A_{2}=\left(
456\cap1234\right)  -3=1,$ $B_{2}=23,$ $n_{i}^{\prime}=0,$ $n_{i}%
^{\prime\prime}=1$ and we may think of $d_{A|B}$ acting on
$1234\times456$ as $1\times d_{1|23}.$

\vspace*{0.2in} \hspace*{0.7in}%
\setlength{\unitlength}{0.00017in}\begin{picture} (0,0)
\thicklines
\put(2000,-500){\vector(3,-2){7000}}
\put(2900,-3700){\vector(3,-1){5800}}
\put(2000,-11900){\vector(3,2){7000}}
\put(800,-1000){\vector(0,-1){1500}}
\put(800,-4000){\vector(0,-1){4500}}
\put(2900,-8700){\vector(3,1){5800}}
\put(800,-11400){\vector(0,1){1500}}
\put(800,0){\makebox(0,0){$P_{n-1}$}}
\put(800,-12200){\makebox(0,0){$P_{n}$}}
\put(-2000,-1600){\makebox(0,0){$_{\Delta_{n_1^{\prime}+1,n_1^{\prime
\prime}+1,n_2+1}}$}}
\put(-950,-6200){\makebox(0,0){$_{h_{A_{1}|B_{1}}\times 1}$}}
\put(800,-3200){\makebox(0,0){$_{P_{n_1^{\prime}+1} \times
P_{n_1^{\prime \prime }+1}\times P_{n_2+1}}$}}
\put(800,-9200){\makebox(0,0){$_{ P_{n_1+1}\times P_{n_2+1}}$}}
\put(9600,-6200){\makebox(0,0){$Y$}}
\put(-1300,-10800){\makebox(0,0){$_{\Delta_{n_1+1,n_2+1}}$}}
\put(5600,-7200){\makebox(0,0){$\bar{f}$}}
\put(5600,-10300){\makebox(0,0){$f$}}
\put(6400,-1900){\makebox(0,0){$d_{A|B}(f)$}}
\put(5600,-5450){\makebox(0,0){$\tilde{f}$}}
\put(17200,-500){\vector(-3,-2){7000}}
\put(16300,-3700){\vector(-3,-1){5800}}
\put(17200,-11900){\vector(-3,2){7000}}
\put(18400,-1000){\vector(0,-1){1500}}
\put(18400,-4000){\vector(0,-1){4500}}
\put(16300,-8700){\vector(-3,1){5800}}
\put(18400,-11400){\vector(0,1){1500}}
\put(18400,0){\makebox(0,0){$P_{n}$}}
\put(18400,-12200){\makebox(0,0){$P_{n-1}$}}
\put(21600,-10800){\makebox(0,0){$_{\Delta_{n_1^{\prime}+1,n_1^{\prime\prime}+1,n_2+1}}$}}
\put(20200,-6200){\makebox(0,0){$_{\phi_{A_{1}|B_{1}} \times 1}$}}
\put(18400,-9200){\makebox(0,0){$_{P_{n_1^{\prime}+1} \times
P_{n_1^{\prime \prime }+1}\times P_{n_2+1}}$}}
\put(18400,-3200){\makebox(0,0){$_{ P_{n_1+1}\times P_{n_2+1}}$}}
\put(20700,-1600){\makebox(0,0){$_{\Delta_{n_1+1,n_2+1}}$}}
\put(13600,-5350){\makebox(0,0){$\tilde{g}$}}
\put(13600,-10300){\makebox(0,0){$g$}}
\put(12600,-1900){\makebox(0,0){$\varrho_{A|B}(g)$}}
\put(13600,-7000){\makebox(0,0){$\bar{g}$}}
\end{picture}\vspace*{2.2in}

\begin{center}
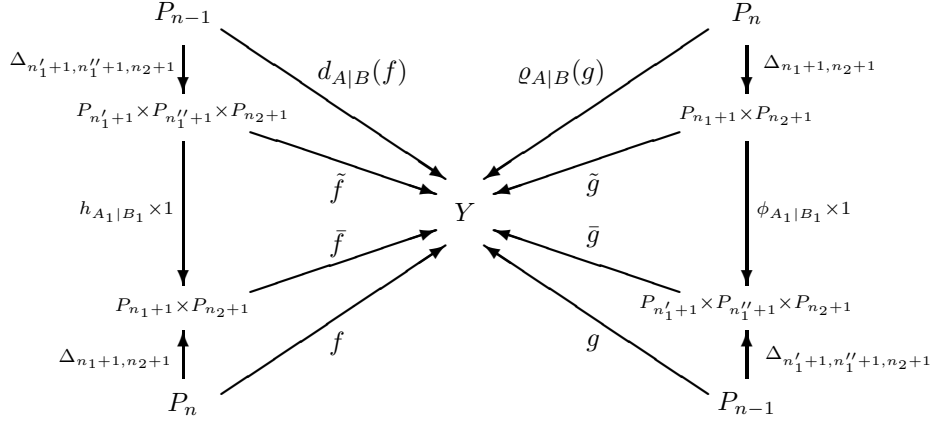
Figure 8:\ Face and degeneracy operators when $i=1$ and
$k=2.$\vspace*{0.2in}
\end{center}

For $f=\bar{f}\circ\Delta_{n_{1}+1\cdots n_{k}+1}\in
\Sing_{n}^{M}Y,$ let $\tilde{f}=\bar{f}\circ(1^{\times i-1}\times
h_{A_{i}|B_{i}}\times1^{\times k-i})$ and define
\[
d_{A|B}(f)=\tilde{f}\circ\Delta_{n_{1}+1\cdots n_{i}^{\prime}+1,n_{i}%
^{\prime\prime}+1\cdots n_{k}+1}.
\]
Dually, note that $n_{i}^{\prime}+n_{i}^{\prime\prime}=n_{i}-1$
implies the sum of coordinates
$(n_{1},\ldots,n_{i-1},n_{i}^{\prime},n_{i}^{\prime\prime },$
$n_{i+1},\ldots,n_{k})\in\left(  \mathbb{N}\cup0\right)  ^{k+1}$
is $n-2.$
So for $g=\bar{g}\circ\Delta_{n_{1}+1\cdots n_{i}^{\prime}+1,n_{i}%
^{\prime\prime}+1\cdots n_{k}+1}\in \Sing_{n-1}^{M}Y,$ let
$\tilde{g}=\bar {g}\circ(1^{\times
i-1}\times\phi_{A_{i}|B_{i}}\times1^{\times k-i})$ and
define%
\[
\varrho_{A|B}(g)=\tilde{g}\circ\Delta_{n_{1}+1\cdots n_{k}+1}%
\]
(see Figure 8).

On the other hand, assume that
$A|B\notin{\mathcal{Q}}_{p_{i},q_{i}}(n)$ for all $i$ and define
$d_{A|B}$ inductively as follows: When $k=2,$ set
$r=n_{1}+1,$ $s=n_{2}+1$ and let%
\[%
\begin{array}
[c]{ll}%
K|L= & \left\{
\begin{array}
[c]{ll}%
(\underline{r}\cap A)\cup\overline{s}\,|\,\underline{r}\cap B, & r\in A\\
\underline{r}\cap A\,|\,(\underline{r}\cap B)\cup\overline{s}, &
r\in B
\end{array}
\right. \\
& \\
M|N= & \left\{
\begin{array}
[c]{lll}%
(\overline{s}\cap
A)-1\,|\,\underline{n-1}\setminus(\overline{s}\cap A)-1, &
r\in B & \\
\underline{n-1}\setminus(\overline{s}\cap
B)-\#L\,|\,(\overline{s}\cap
B)-\#L, & r\in A, & n\in A\\
I_{\underline{n}\setminus L}(A)\,|\,\underline{n-1}\setminus
I_{\underline {n}\setminus L}(A), & r\in A, & n\in B
\end{array}
\right. \\
& \\
C|D= & \left\{
\begin{array}
[c]{lll}%
I_{\underline{n}\setminus B}(\underline{r}\cap
A)\,|\,\underline{n-1}\setminus
I_{\underline{n}\setminus B}(\underline{r}\cap A), & r\in B, & n\in B\\
I_{\underline{n}\setminus A}(\overline{s}\cap
B)\,|\,\underline{n-1}\setminus
I_{\underline{n}\setminus A}(\overline{s}\cap B), & r\in A, & n\in B\\
&  & \\
\underline{n-1}\setminus I_{\underline{n}\setminus
B}(\overline{s}\cap
A)\,|\,I_{\underline{n}\setminus B}(\overline{s}\cap A), & r\in B, & n\in A\\
\underline{n-1}\setminus I_{\underline{n}\setminus
A}(\underline{r}\cap B)\,|\,I_{\underline{n}\setminus
A}(\underline{r}\cap B), & r\in A, & n\in A.
\end{array}
\right.
\end{array}
\]
Then define
\begin{equation}
d_{A|B}=\varrho_{C|D}d_{M|N}d_{K|L}. \label{multip}%
\end{equation}

\begin{remark}
This definition makes sense since $K|L\in{\mathcal{Q}}_{p_{1},q_{1}}%
(n),$  $ M|N  \in \linebreak {\mathcal{Q}}_{p_{3},q_{3}}(n-1),$
 $C|D\in{\mathcal{Q}}%
_{p_{1},q_{1}}(n-1)\,$with either $r,n\in B$ or $r,n\in
A$\thinspace and $C|D\in{\mathcal{Q}}_{p_{3},q_{3}}(n-1)\,$with
either $r\in B,$ $n\in A$ or $r\in A,$ $n\in B$. Of course,
${\mathcal{Q}}_{\ast\ast}(n-1)$ is considered with respect to the
decomposition $n-2=m_{1}+m_{2}+m_{3}$ fixed after the action of
$d_{K|L}(\underline{r}\times\underline{s})$.
\end{remark}

If $k=3,$ consider the pair $(r,s)=(n_{1}+1,n-n_{1}),$ then $(r_{1}%
,s_{1})=(n_{2}+1,n-n_{1}-n_{2}-1)$ for $A_{1}|B_{1}=I_{\underline{n}%
\setminus\underline{r}}(\overline{s}\cap
A)|I_{\underline{n}\setminus \underline{r}}(\overline{s}\cap
B)\in{\mathcal{P}}_{p_{1},q_{1}}(n-r),$ and so on. Now dualize and
use the same formulas above to define the degeneracy operator
$\varrho_{A|B}$.

\begin{definition}\label{multisingular}
Let $Y$ be a topological space. The \underline{singular
multipermutahedral set} \underline{of $Y$} consists of the
singular set $\Sing_{\ast}^{M}Y$ together with the singular face
and degeneracy operators
\[
d_{A|B}:\Sing_{n}^{M}Y\rightarrow \Sing_{n-1}^{M}Y\text{ \ and \
}\varrho _{A|B}:\Sing_{n-1}^{M}Y\rightarrow \Sing_{n}^{M}Y
\]
defined respectively for each $n\geq2$ and
$A|B\in{\mathcal{P}}_{\ast\ast }\left(  n\right)  $.
\end{definition}

\begin{remark}
The operator $d_{A|B}$ defined in (\ref{multip}) applied to
$d_{U|V}$ for some$\ U|V\in{\mathcal{P}}_{r,s}(n+1)$ yields the
higher order structural relation
\begin{equation}
d_{A|B}d_{U|V}=\varrho_{C|D}d_{M|N}d_{K|L}d_{U|V} \label{hrelation}%
\end{equation}
discussed in our first universal example.
\end{remark}

Now $\Sing_{\ast}^{M}Y$ determines the singular (co)homology of a
space $Y$ in the following way: Let $R$ be a commutative ring with
identity. For $n\geq1,$ let $C_{n-1}(\Sing^{M}Y)$ denote the
$R$-module generated by $\Sing_{n}^{M}Y$
and form the \textquotedblleft chain complex\textquotedblright%
\[
(C_{\ast}(\Sing^{M}Y),d)=\bigoplus_{\substack{n_{(k)}=n-1 \\n\geq1}%
}(C_{n-1}(\Sing^{n_{1}\cdots n_{k}}Y),d_{n_{1}\cdots n_{k}}),
\]
where
\[
d_{n_{1}\cdots n_{k}}=\sum_{\substack{A|B\in\bigcup_{i=1}^{k}\mathcal{Q}%
_{p_{i},q_{i}}\left(  n\right)
}}\,-(-1)^{n_{(i-1)}+n_{i}^{\prime}}\text{ \textit{
shuff}}(C_{i}\cap A;C_{i}\cap B)\text{ }d_{A|B}.
\]
Refer to the example in Figure 7 and note that for $f\in
C_{4}(\Sing^{M}Y)$ with\linebreak\ $d_{13|2}d_{12|34}\left(
f\right)  \neq0$, the component $d_{13|2}d_{12|34}\left(  f\right)
$ of $d^{2}\left(  f\right)  \in C_{2}(\Sing^{M}Y)$ is not
cancelled and $d^{2}\neq0$. Hence $d$ is not a differential. To
remedy this, form the quotient
\[
C_{\ast}^{\diamondsuit}(Y)=C_{\ast}\left(  \Sing^{M}Y\right) /DGN,
\]
where $DGN$ is the submodule generated by the degeneracies, and
obtain the \emph{singular permutahedral chain complex }$\left(
C_{\ast}^{\diamondsuit }(Y),d\right)  $. Because the signs in $d$
are determined by the index $i$, which is missing in our first
universal example, we are unable to use our first example to
define a chain complex with signs. However, we could use it to
define a unoriented theory with $\mathbb{Z}_{2}$-coefficients.

The singular homology of $Y$ is recovered from the composition
\[
C_{\ast}(\Sing Y)\rightarrow C_{\ast}(\Sing^{I}Y)\rightarrow
C_{\ast} (\Sing^{M}Y)\rightarrow C_{\ast}^{\diamondsuit}(Y)
\]
arising from the canonical cellular projections
\[
P_{n+1}\rightarrow I^{n}\rightarrow\Delta^{n}.
\]
Since this composition is a chain map, there is a natural
isomorphism
\[
H_{\ast}(Y)\approx
H_{\ast}^{\diamondsuit}(Y)=H_{\ast}(C_{\ast}^{\diamondsuit
}(Y),d).
\]
The fact that our diagonal on $P$ and the A-W diagonal on
simplices commute with projections allows us to recover the
singular cohomology ring of $Y$ as well. Finally, we remark that a
cellular projection $f$ between polytopes induces a chain map
between corresponding singular chain complexes whenever chains on
the target are normalized. Here $C_{\ast}(\Sing Y)$ and $C_{\ast
}(\Sing^{I}Y)$ are non-normalized and the induced map $f^{\ast}$
is not a chain map; but fortunately $d^{2}=0$ does not depend
$df^{\ast}=f^{\ast}d$.

\subsection{Abstract Permutahedral Sets}

We begin by constructing a generating category $\mathbf{P}$ for
permutahedral sets similar to that of finite ordered sets and
monotonic maps for simplicial sets. The objects of $\mathbf{P}$
are the sets $n!=S_{n}$ of permutations of $\underline{n},$
$n\geq1.$ But before we can define the morphisms we need some
preliminaries. First note that when $P_{n}$ is identified with its
vertices $n!$, the maps $\rho_{n}$ and $\gamma_{n}$ defined above
become
\[
\rho_{n}:n!\rightarrow2!^{n-1}\ \ \text{and}\ \
\gamma_{n}:2!^{n-1}\rightarrow n!.
\]
Given a non-empty increasingly ordered set $M=\left\{ m_{1}<\cdots
<m_{k}\right\}  \subset\mathbb{N},$ let $M!$ denote the set of all
permutations of $M$ and let $J_{M}:M!\rightarrow k!$ be the map
defined for $a=\left(  m_{\sigma\left(  1\right)
},...,m_{\sigma\left(  k\right) }\right)  \in M!$ by
$J_{M}(a)=\sigma.$ For $n,m\in\mathbb{N}$ and partitions
$A_{1}|\cdots|A_{k}\in{\mathcal{P}}_{n_{1}\cdots n_{k}}(n)$ and $B_{1}%
|\cdots|B_{\ell}\in{\mathcal{P}}_{m_{1}\cdots m_{\ell}}(m)$ with
$n-k=m-\ell=\varkappa,$ define the morphism
\[
f_{A_{1}|\cdots|A_{k}}^{B_{1}|\cdots|B_{\ell}}:m!\rightarrow n!
\]
by the composition
\begin{multline*}
m!\overset{sh_{B}}{\longrightarrow}\prod_{j=1}^{\ell}B_{j}\overset{\sigma_{\max}%
}{\longrightarrow}\prod_{r=1}^{\ell}B_{j_{r}}\overset{J_{B}}{\longrightarrow}%
\prod_{j=r}^{\ell}m_{j_{r}}!\overset{\rho_{\ast}}{\longrightarrow}2!^{\varkappa
}\overset{\gamma_{\ast}}{\longrightarrow}
\\
\prod_{s=1}^{k}n_{i_{s}}!\overset
{J_{A}^{-1}}{\longrightarrow}\prod_{s=1}^{k}A_{i_{s}}\overset{\sigma_{\max}^{-1}
}{\longrightarrow}\prod_{i=1}^{k}A_{i}
\overset{\iota_{A}}{\longrightarrow}n!
\end{multline*}
where $sh_{B}$ is a surjection defined for
$b=\{b_{1},...,b_{m}\}\in m!$ by
\[
sh_{B}(b)=(b_{1,1},..,b_{m_{1},1};...;b_{1,\ell},..,b_{m_{\ell},\ell}),\,
\]
in which the right-hand side is the unshuffle of $b$ with
$b_{r,t}\in B_{t},\,1\leq r\leq m_{t},\,1\leq t\leq\ell;$
$\sigma_{\max}\in S_{\ell}$ is a permutation defined by
$j_{r}=\sigma_{\max}(r),\,$ $\max B_{j_{r}}=\max (B_{1}\cup
B_{2}\cup\cdots\cup B_{j_{r}});$ $J_{B}=\prod_{r=1}^{\ell
}J_{B_{j_{r}}};$ $\rho_{\ast}=\prod_{r=1}^{\ell}\rho_{j_{r}}$ and
$\gamma_{\ast}=\prod_{s=1}^{k}\gamma_{i_{s}};$ finally,
$\iota_{A}$ is the inclusion. It is easy to see that
\[
f_{A_{1}|\cdots|A_{k}}^{B_{1}|\cdots|B_{\ell}}=f_{A_{1}|\cdots|A_{k}%
}^{\underline{\varkappa+1}}\circ f_{\underline{\varkappa+1}}^{B_{1}%
|\cdots|B_{\ell}}\ \ \text{and}\ \ f_{\underline{n}}^{\underline{n}}%
=\gamma_{n}\circ\rho_{n}.
\]
In particular, the maps
$f_{A|B}^{\underline{n-1}}:(n-1)!\rightarrow n!$ and
$f_{\underline{n-1}}^{A|B}:n!\rightarrow\left(  n-1\right)  !$ are
generator morphisms denoted by $\delta_{A|B}$ and $\beta_{A|B},$
respectively (see Theorem \ref{highrelations} below, the statement
of which requires some new set operations).

\begin{definition}
Given non-empty disjoint subsets $A,B,U\subset\underline{n+1}$
with $A\cup B\subseteq U,$ define the \underline{lower and upper
disjoint unions} (with
respect to $U$) by\vspace{0.1in}\newline%
\begin{tabular}
[c]{lll}%
$\hspace*{0.3in}$ & $A\underline{\sqcup}B=$ & $\left\{
\begin{array}
[c]{ll}%
I_{U\diagdown A}\left(  B\right)  +\#A-1, & \text{if }\min
B>\min\left(
U\diagdown A\right) \\
I_{U\diagdown A}\left(  B\right)  +\#A-1\cup\underline{\#A}, &
\text{if }\min B=\min\left(  U\diagdown A\right)
\end{array}
\right.  $\\
and &  & \\
& $A\overline{\sqcup}B=$ & $\left\{
\begin{array}
[c]{ll}%
I_{U\diagdown B}\left(  A\right)  , & \text{if }\max A<\max\left(
U\diagdown
B\right) \\
I_{U\diagdown B}\left(  A\right)  \cup\overline{\#B}-1, & \text{if
}\max
A=\max\left(  U\diagdown B\right)  \text{.}%
\end{array}
\right.  $%
\end{tabular}
\vspace{0.1in}\newline If either $A$ or $B$ is empty, define
$A\underline {\sqcup}B=A\overline{\sqcup}B=A\cup B$. Furthermore,
given non-empty disjoint subsets
$A,B_{1},\ldots,B_{k}\subset\underline{n+1}$ with $k\geq1,$ set
$U=A\cup B_{1}\cup\cdots\cup B_{k}$ and define
\[
A\square(B_{1}|\cdots|B_{k})=(B_{1}|\cdots|B_{k})\square A=\left\{
\begin{array}
[c]{ll}%
A\underline{\sqcup}B_{1}|\cdots|A\underline{\sqcup}B_{k}, &
\text{if }\max
A<\max U\\
B_{1}\overline{\sqcup}A|\cdots|B_{k}\overline{\sqcup}A, & \text{if
}\max A=\max U.
\end{array}
\right.
\]

\end{definition}

\noindent Note that if $A|B$ is a partition of $\underline{n+1}$,
then
\[
A\underline{\sqcup}B=A\overline{\sqcup}B=\underline{n}.
\]
Given a partition $A_{1}|\cdots|A_{k+1}$ of $\underline{n},$ define $A_{1}%
^{1}|\cdots|A_{k+1}^{1}=A_{1}^{1}|\cdots|A_{1}^{k+1}=$\linebreak$A_{1}%
|\cdots|A_{k+1};$ inductively, given
$A_{1}^{i}|\cdots|A_{k-i+2}^{i}$ the partition of
$\underline{n-i+1},\,1\leq i<k,$ let
\[
A_{1}^{i+1}|\cdots|A_{k-i+1}^{i+1}=A_{1}^{i}\square(A_{2}^{i}|\cdots
|A_{k-i+2}^{i})
\]
be the partition of $\underline{n-i};$ and given $A_{i}^{1}|\cdots
|A_{i}^{k-i+2}$ the partition of $\underline{n-i+1},\,1\leq i<k,$
let
\[
A_{i+1}^{1}|\cdots|A_{i+1}^{k-i+1}=(A_{i}^{1}|\cdots|A_{i}^{k-i+1})\square
A_{i}^{k-i+2}%
\]
be the partition of $\underline{n-i}.$

\begin{theorem}
\label{highrelations}For
$A_{1}|\cdots|A_{k+1}\in{\mathcal{P}}_{n_{1}\cdots
n_{k+1}}(n),$ $2\leq k\leq n,$ the map $f_{A_{1}|\cdots|A_{k+1}}%
^{\underline{n-k}}:(n-k)!\rightarrow n!\,$ can be expressed as a
composition of $\delta$'s two ways:
\[
f_{A_{1}|\cdots|A_{k+1}}^{\underline{n-k}}=\delta_{A_{1}^{1}\,|\,A_{2}^{1}%
\cup\cdots\cup A_{k+1}^{1}}\cdots\delta_{A_{1}^{k}\,|\,A_{2}^{k}}%
=\delta_{A_{1}^{1}\cup\cdots\cup
A_{1}^{k}\,|\,A_{1}^{k+1}}\cdots\delta _{A_{k}^{1}\,|\,A_{k}^{2}}.
\]

\end{theorem}

\begin{proof}
The proof is straightforward and omitted.
\end{proof}

\noindent There is also the dual set of relations among the
$\beta$'s.

\begin{example}
Theorem \ref{highrelations} defines structure relations among the
$\delta$'s, the first of which is
\begin{equation}
\delta_{A|B\cup C}\,\delta_{A\square(B|C)}=\delta_{A\cup
B|C}\,\delta
_{(A|B)\square C} \label{coquadrel}%
\end{equation}
when $k=2.$ In particular, let $A|B|C=12|345|678$. Since
$A\underline{\sqcup }B=\{1234\},$ $A\underline{\sqcup}C=\left\{
567\right\}  ,$ $A\overline {\sqcup}C=\{12\}$ and
$B\overline{\sqcup}C=\left\{  34567\right\}  ,$ we obtain the
following quadratic relation on $12|345|678$:
\[
\delta_{12|345678}\delta_{1234|567}=\delta_{12345|678}\delta_{12|34567};
\]
similarly, on $345|12|678$ we have
\[
\delta_{345|12678}\delta_{1234|567}=\delta_{12345|678}\delta_{34567|12}.
\]

\end{example}

\begin{definition}
Let $\mathcal{C}$ be the category of sets. A
\underline{permutahedral set} is a contravariant functor
\[
\mathcal{Z}:\mathbf{{P}\rightarrow\mathcal{C}}.
\]

\end{definition}

Thus a permutahedral set $\mathcal{Z}$ is a graded set $\mathcal{Z}%
=\{\mathcal{Z}_{n}\}_{n\geq1}\,$endowed with face and degeneracy
operators
\[
d_{A|B}=\mathcal{Z}(\delta_{A|B}):\mathcal{Z}_{n}\rightarrow\mathcal{Z}%
_{n-1}\text{ \ and \ }\varrho_{M|N}=\mathcal{Z}(\beta_{M|N}):\mathcal{Z}%
_{n}\rightarrow\mathcal{Z}_{n+1}%
\]
satisfying an appropriate set of relations, which includes
quadratic relations such as
\begin{equation}
d_{A\square(B|C)}d_{A|B\cup C}\ =d_{(A|B)\square C}d_{A\cup B|C}
\label{quadrel}%
\end{equation}
induced by (\ref{coquadrel}) and higher order relations such as
\[
d_{A|B}d_{U|V}=\varrho_{C|D}d_{M|N}d_{K|L}d_{U|V}%
\]
discussed in (\ref{hrelation}).

Let us define the abstract analog of a singular multipermutahedral
set, which leads to a singular chain complex with arbitrary
coefficients.

\begin{definition}
\label{derformal}For $n\geq1,\mathbb{\ }$let
$X_{n}=\bigcup_{n_{\left( k\right)  }=n-1,\text{
}n_{k}\geq0\,}X^{n_{1}\cdots n_{k}}$ and $X_{n-1}=$
\linebreak$\bigcup_{m_{\left(  \ell\right)  }=n-2,\text{
}m_{\ell}\geq
0\,}X^{m_{1}\cdots m_{\ell}}$ be filtered sets; let $A|B\in{\mathcal{Q}%
}_{p_{i},q_{i}}(n)$ for some $i.$ A map $g:X_{n}\rightarrow
X_{n-1}$ acts as
an \underline{$A|B$-formal derivation} if $g|_{X^{n_{1}\cdots n_{k}}}%
:X^{n_{1}\cdots n_{k}}\rightarrow X^{n_{1}\cdots n_{i}^{\prime},n_{i}%
^{\prime\prime}\cdots n_{k}},$ where
$(n_{i}^{\prime},n_{i}^{\prime\prime})$ is given by
(\ref{formal}).
\end{definition}

Let $\mathcal{C}_{M}$ denote the category whose objects are
positively graded sets $X_{\ast}$ filtered by subsets
$X_{n}=\bigcup_{n_{\left(  k\right) }=n-1,\text{
}n_{k}\geq0}X^{n_{1}\cdots n_{k}}$ and whose morphisms are
filtration preserving set maps.

\begin{definition}
A \underline{multipermutahedral set} is a contravariant functor $\mathcal{Z}%
:\mathbf{{P}\rightarrow\mathcal{C}}_{M}$ such that
\[
\mathcal{Z}(\delta_{A|B}):\mathcal{Z}(n!)\rightarrow\mathcal{Z}((n-1)!)
\]
acts as an $A|B$-formal derivation for each $A|B\in{\mathcal{Q}}_{p_{i},q_{i}%
},$ all $i\geq1.$
\end{definition}

Thus a multipermutahedral set $\mathcal{Z}$ is a graded set
$\left\{ \mathcal{Z}_{n}\right\}  _{n\geq1}$ with
\[
\mathcal{Z}_{n}=\bigcup_{\substack{n_{\left(  k\right)  }=n-1
\\n_{k}\geq 0}}\mathcal{Z}^{n_{1}\cdots n_{k}},
\]
$\,$together with face and degeneracy operators
\[
d_{A|B}=\mathcal{Z}(\delta_{A|B}):\mathcal{Z}_{n}\rightarrow\mathcal{Z}%
_{n-1}\text{ \ and \ }\varrho_{M|N}=\mathcal{Z}(\beta_{M|N}):\mathcal{Z}%
_{n}\rightarrow\mathcal{Z}_{n+1}%
\]
satisfying the relations of a permutahedral set and the additional
requirement that $d_{A|B}$ respect underlying multigrading. This
later condition allows us to form the chain complex of
$\mathcal{Z}$ with signs mimicking the cellular chain complex of
permutahedra (see below). Note that the chain complex of a
permutahedral set is only defined with $\mathbb{Z}_{2}
$-coefficients in general.

\subsection{The Cartesian product of permutahedral sets}

The objects and morphisms in the category $\mathbf{P\times P}$ are
the sets and maps
\[
n!!=\bigcup_{r+s=n}r!\times s!\text{ \ and \ }\bigcup_{f,g\in\mathbf{P}%
}f\times g:m!!\rightarrow n!!
\]
all $m,n\geq1.$ There is a functor $\Delta:\mathbf{P\rightarrow
P\times P}$ defined as follows. If
${A|B}\in\mathcal{Q}_{r,1}(n)\cup\mathcal{Q}_{1,s}(n),$ define
$\Delta_{r,s}(A|B)=A_{1}|B_{1}\times A_{2}|B_{2}\in r!\times s!\
$and define $\delta_{A|B}:(n-1)!\rightarrow n!$ by
\[
\Delta(\delta_{A|B})=\delta_{A_{1}|B_{1}}\times\delta_{A_{2}|B_{2}},
\]
where $\delta_{A_{i}|B_{i}}=1$ for either $i=1$ or $i=2.$ Define
$\Delta
(\beta_{A|B})$ similarly$.$ On the other hand, if $A|B\notin\mathcal{Q}%
_{r,1}(n)\cup\mathcal{Q}_{1,s}(n),$ define
\[
\Delta(\delta_{A|B})=\Delta(\delta_{K|L})\Delta(\delta_{M|N})\Delta
(\beta_{C|D}),
\]
where $K|L,\,M|N,\,C|D$ are given by the formulas in
(\ref{multip}). Dually, define $\Delta(\beta_{M|N}).$ It is easy
to check that $\Delta$ is well defined.

Given multipermutahedral sets
$\mathcal{Z}^{\prime},\mathcal{Z}^{\prime\prime
}:\mathbf{P}\rightarrow\mathcal{C}_{M}$, first define a functor
\[
\mathcal{Z}^{\prime}\tilde{\times}\mathcal{Z}^{\prime\prime}:\mathbf{P}%
\times\mathbf{P}\rightarrow\mathcal{C}_{M}%
\]
on an object $n!!$ by
\[
(\mathcal{Z}^{\prime}\tilde{\times}\mathcal{Z}^{\prime\prime})(n!!)=\bigcup
_{r+s=n}\mathcal{Z}^{\prime}(r!)\times\mathcal{Z}^{\prime\prime}
(s!)\diagup{ \sim},
\]
where $\left(  a,b\right)  \sim\left(  c,e\right)  $ if and only
if $a=\varrho_{\underline{r}|r+1}^{\prime}(c)$ and
$e=\varrho_{1|\underline
{s+1}\setminus\underline{1}}^{\prime\prime}(b).$ On a map
$h=\bigcup(f\times g):m!!\rightarrow n!!,$
\[
(\mathcal{Z}^{\prime}\tilde{\times}\mathcal{Z}^{\prime\prime})(h):(\mathcal{Z}%
^{\prime}\tilde{\times}\mathcal{Z}^{\prime\prime})(n!!)\rightarrow
(\mathcal{Z}^{\prime}\tilde{\times}\mathcal{Z}^{\prime\prime})(m!!)
\]
is the map induced by $\bigcup(\mathcal{Z}^{\prime}(f)\times\mathcal{Z}%
^{\prime\prime}(g)).$ Now define the product $\mathcal{Z}^{\prime}%
\times\mathcal{Z}^{\prime\prime}$ to be the composition of
functors
\[
\mathcal{Z}^{\prime}\times\mathcal{Z}^{\prime\prime}=\mathcal{Z}^{\prime
}\tilde{\times}\mathcal{Z}^{\prime\prime}\circ\Delta:\mathbf{P\rightarrow
\mathcal{C}}_{M}.
\]
The face operator $d_{A|B}$ on $\mathcal{Z}^{\prime}\times\mathcal{Z}%
^{\prime\prime}$ is given by
\begin{equation}
d_{A|B}(a\times b)=\left\{
\begin{array}
[c]{ll}%
d_{\underline{r}\cap A|\underline{r}\cap B}^{\prime}\left(
a\right)  \times
b, & \text{\textit{if }}A|B\in\mathcal{Q}_{1,s}\left(  n\right)  ,\\
a\times d_{(\overline{s}\cap A)-r+1\,|\,(\overline{s}\cap B)-r+1}%
^{\prime\prime}(b), & \text{\textit{if
}}A|B\in\mathcal{Q}_{r,1}\left(
n\right)  ,\\
\varrho_{C|D}d_{M|N}d_{K|L}\left(  a\times b\right)  , &
\text{\textit{otherwise}},
\end{array}
\right.  \label{productd1}%
\end{equation}
with $M|N,\,K|L,\,C|D$ given by the formulas in (\ref{multip}).

\begin{example}\label{monoidal}
The canonical map $\iota:\Sing^{P}X\times \Sing^{P}Y\rightarrow
\Sing^{P}(X\times Y)$ defined for $(f,g)\in \Sing_{r}^{P}X\times
\Sing_{s}^{P}Y$ by
\[
\iota(f,g)=(f\times g)\circ\Delta_{r,s}
\]
is a map of permutahedral sets. Consequently, if $X$ is a
topological monoid, the singular permutahedral complex $Sing^{P}X$
inherits a canonical monoidal structure.
\end{example}


\begin{definition}

A  monoidal permutahedral  set is a permutahedral set $ {\mathcal
Z}$ with a map $ \mu :{\mathcal Z}\times {\mathcal Z}\rightarrow
{\mathcal Z} $ of permutahedral sets which is associative and has
the unit $e\in \mathcal Z_1.$
\end{definition}
 Clearly, for a
monoidal multipermutahedral set ${\mathcal Z},$ its chain complex
$(C^{\diamondsuit}_*({\mathcal Z};R),\\ d)$ is a dg Hopf algebra.

Given a monoidal  multipermutahedral  set ${\mathcal Z},$  a
${\mathcal Z}$-{\em module} is a multipermutahedral set ${\mathcal
L}$
 together with associative action $
{\mathcal Z}\times {\mathcal L}\rightarrow {\mathcal L}$ with the
unit of ${\mathcal Z}$ acting as identity. In this case
$C_{\diamondsuit}^*({\mathcal L};R)$ is a dga comodule over the dg
Hopf algebra $(C_{\diamondsuit}^*({\mathcal Z};R),d)$.

\subsection{The permutahedral set functor ${\bf \Omega} Q$}\
Let $Q=(Q_n,d^0_i,d^1_i,\eta_i)_{n\geq 0}$ be a cubical set.
Recall that the diagonal
$$\Delta : C^{\Box}_*(Q)\to   C^{\Box}_*(Q)\otimes  C^{\Box}_*(Q)$$
of  $Q$ is defined on $a\in Q_n$ by
$$
\Delta (a) =\sum \text{{\it shuff}}\,(A;B)\,    d^0_{B}
  (a)\otimes d^1_{A}(a),
$$
where $d^0_{B}=d^0_{j_1}...d^0_{j_{q}},$\
$d^1_{A}=d^1_{i_1}...d^1_{i_{p}},$
 the summation is over all shuffles $\{A,
B\}=\{i_1<...<i_q, j_1<...<j_p \}$ of the set $\uu {n}.$ The
primitive components  of the diagonal are given by the extreme
cases $A=\oo$ and $B=\oo.$

Assume $Q$ is 1-reduced. Let $\bar{Q}=s^{-1}(Q_{>0})$ denote the
desuspension of $Q,$ let $\mathbf{\Omega}^{\prime\prime}Q$ be the
free graded monoid generated by $\bar{Q}$ with the unit
$e\in\bar{Q}_{1}\subset\mathbf{\Omega}^{\prime\prime }Q$ and let
$\Upsilon$ be the set of formal expressions
\[
\Upsilon=\{\varrho_{M_{k}|N_{k}}((\cdots\varrho_{M_{2}|N_{2}}(\varrho
_{M_{1}|N_{1}}(\bar{a}_{1}\cdot\bar{a}_{2})\cdot\bar{a}_{3})\cdots)\cdot
\bar{a}_{k+1})|\,a_{i}\in Q_{r_{i}}\}_{r_{i}\geq1;k\geq2},
\]
$M_{i}|N_{i}\in{\mathcal{P}}_{r_{(i)},r_{i+1}}(r_{(i+1)})$ or $M_{i}|N_{i}%
\in{\mathcal{P}}_{r_{i+1},r_{(i)}}(r_{(i+1)}),$
$r_{(i)}=r_1+\cdots+r_i,$ $\,1\leq i\leq k. $ Note that one or
more of the $a_{i}$'s can be the unit $e$. Adjoin the elements of
$\Upsilon$ to $\mathbf{\Omega}^{\prime\prime}Q$ and obtain the
graded monoid $\mathbf{\Omega}^{\prime}Q\ $and let
$\mathbf{\Omega}Q$ be the monoid
\[
\mathbf{\Omega}Q=\mathbf{\Omega}^{\prime}Q\diagup{\sim}\text{,}%
\]
where
$\varrho_{M|N}(\bar{a}\cdot\bar{b})\sim\varrho_{N|M}(\bar{a}\cdot\bar
{b}),\,$ $\varrho_{j|\underline{n}\setminus j}(e\cdot\bar{a})\sim
\varrho_{\underline{n}\setminus j|j}(\bar{a}\cdot e)\sim\overline{\eta_{j}%
(a)},\,a,b\in Q_{>0},$ and $\bar{a}_{1}\cdots{\varrho_{\underline{r_{i}}%
|r_{i}+1}(\bar{a}_{i}\cdot
e)}\cdot\bar{a}_{i+2}\cdots\bar{a}_{k+1}\sim\bar
{a}_{1}\cdots\bar{a}_{i}\cdot\varrho_{1|\underline{r_{i+2}+1}\setminus
1}(e\cdot\bar{a}_{i+2})\cdots\bar{a}_{k+1}$ for $a_{i}\in Q_{r_{i}}%
,a_{i+1}=e,\,1\leq i\leq k.$ Then $\mathbf{\Omega}Q$ is
canonically a multipermutahedral set in the following way: First,
define the face operator $d_{A|B}$ on a monoidal generator
$\bar{a}\in\bar{Q}_{n}$ by
\[
d_{A|B}(\bar{a})=\overline{d_{B}^{0}(a)}\cdot\overline{d_{A}^{1}%
(a)},\ \ \ A|B\in\mathcal{P}_{\ast,\ast}{(n)}.
\]
Next, use the formulas in the definition of a singular
multipermutahedral set (\ref{multip}) to define $d_{A|B}$ and
$\varrho_{M|N}$ on decomposables. In
particular, the following identities hold for $1\leq i\leq n$:%
\[
d_{i|\underline{n+1}\setminus i}\left(  \overline{a}\right)
=\overline {d_{i}^{1}(a)}\text{ \ and \
}d_{\underline{n+1}\setminus i|i}\left( \overline{a}\right)
=\overline{d_{i}^{0}(a)}.
\]
It is easy to see that
$(\mathbf{\Omega}Q,\,d_{A|B},\,\varrho_{M|N})$ is a
multipermutahedral set that depends functorially on $Q.$

\begin{remark}
\label{lift}The fact that the definition of $\mathbf{\Omega}Q$
uses all cubical degeneracies is justified geometrically by the
fact that a degenerate singular $n$-cube in the base of a path
fibration lifts to a singular $(n-1)$-permutahedron in the fibre,
which is degenerate with respect to Milgram's projections
\cite{Milgram} (c.f., the definition of the cubical set
$\mathbf{\Omega}X$ on a simplicial set $X$ \cite{KS1}).
\end{remark}

\section{The permutocubes}\label{secpermutocube}

The pertmutocube $B_n$ is an n-dimensional  polytope discovered by
N. Berikashvili  which can be thought of as a "twisted Cartesian
product" of the cube and the permutahedron. Originally the
permutocube $B_n$ has been obtained   from $I^n$ by the following
truncation procedure: First the n-cube is truncated at the minimal
vertex $a_0=(0,...,0),$ then it is truncated along those
$(n-1)$-faces that contained $a_0,$ and continuing so the last
truncation is along those 1-faces (edges) of the n-cube
 that contained  $a_0$. Hence, $B_2$ is a pentagon (Figure 9), for $B_3$ see
Figure 10. In particular, this truncation procedure fixes the
permutahedron $P_{n}$
 at the vertex $a_0.$
 So that we get
 the  natural cellular
embedding (see Figures 11 and 12)
\begin{equation}\label{pb}
\delta_{0]\uu{n}}: P_{n}\to B_{n}.
\end{equation}

The notation for the above inclusion map is motivated by the
following combinatorial description of $B_n.$  We have that  the
faces of $B_n$ are in one-to-one correspondence with partitions
$C_0] C_1|...|C_p$ of all (non-empty) subsets  of the set
 \[ \uu{n_0}=\{0,1,...,n\}\] with $0\in C_0.$  More precisely,
an $(n-k-p)$-face $u$=$ C_0] C_1|...|C_p$ of $B_n$ can be labelled
by the composition of face operators, denoted by
$d_u$=$d_{A_0]A_1|...|A_p;(i_k,...,i_1)}$ for
$\{i_k<\cdots<i_1\}\subset \uu{n},$
\begin{multline*}
u=d_{A_0]A_1|...|A_p;(i_k,...,i_1)}(\uu{n_0})=\\
d_{A_{0}]A_{1}}d_{A_{0}\cup A_{1}]A_2}\cdots d_{A_0\cup A_{1}\cup
\cdots \cup A_{p-1}]A_{p}} {d}_{i_k}\dotsb {d}_{i_1}(\uu{n_0}),
\end{multline*}
where
 $d_i$ acts by deleting  the $(i+1)^{th}$ integer,
while $d_{A]M}$ forms the partition $A]M$ of the domain, and
$A_i=I_{\uu{n_0}\setminus (i_k,...,i_1)}(C_i),\,0\leq i\leq p.$
 In
particular, $d_{0] \uu n}(\uu{n_0})$ just corresponds to the
single $(n-1)$-permutahedral face  $\delta_{0]\uu n}(P_n)\subset
B_n.$ On the other hand, a face of the form ${d}_{i_k}\dotsb
{d}_{i_1}(\uu{n_0})$ may be identified with $B_{n-k},\,0\leq k<
n.$

 For example, for $u=038]14|6|79\subset B_9,$  we have
\[u=d_{026]13|4|57;\,(2,5)}(\uu{9_0})=d_{024]
13}d_{1235]4}d_{12346]57}{d}_2{d}_5(\uu{9_0}).\]

We have that $B_n$ also admits the realization   as a subdivision
of the standard $n$-cube $I^{n}$ compatible with the inclusion
$P_n\subset B_n$ (see Figures 9,10). Indeed, let $B_{0}=I^0;$ for
$1\leq i\leq n,$ let $e_{i,\epsilon}^{n-1}$ denote the $\left(
n-1\right)  $-face $(x_{1},\ldots,x_{i-1},\epsilon,x_{i+1}
,\ldots,x_{n})\subset I^{n}$ and label the endpoints of
$B_{1}=[0,1] $ via $e_{1,0}^0\leftrightarrow d_{ 0]1 }$ and
$e_{1,1}^0\leftrightarrow {d}_{1 }.$ Inductively, if $B_{n-1}$ has
been constructed, obtain $B_{n}$ as a subdivision of
$B_{n-1}\times I$ in the following way:

\vspace{0.2in}

\begin{tabular}
[c]{l|ll} Face of $B_{n}$ & $\text{Label}$\\\hline
& \\
$e^{n-1}_{n,0}$ & $d_{ \uu{(n-1)_0}\,]\,n  }$\\
& \\
$e^{n-1}_{i,1}$ & ${d}_{i},$   &  $  i \in \uu n$\\
& \\
$d_{ A]M   }\times I_{0, \# M}$ & $d_{  A] M\cup n  }$ \\
& \\
$d_{   A] M}\times I_{\# M,\infty}$ & $d_{  A\cup n] M }.$
\end{tabular}

\vspace{0.2in}

Thus,  proper cells of $B_n$ are represented as the Cartesian
product of the permutocube and permutahedra. In particular, on a
proper cell $e=e_1\x e_2\subset B_n$  a permutahedral face
operator $d_{M_1|M_2}$ acts as $d_{M_1|M_2}(e)=e_1\x
d_{M_1|M_2}(e_2).$

 \vspace{0.2in}
\unitlength=0.8mm
 \special{em:linewidth 0.4pt}
\linethickness{0.4pt}
\begin{picture}(125.00,32.67)
\put(120.33,0.67){\line(-1,0){29.67}}
\put(90.66,0.67){\line(0,1){27.00}}
\put(90.66,27.67){\line(1,0){29.33}}
\put(120.00,27.67){\line(0,-1){27.00}}
\put(90.66,13.00){\circle*{1.33}} \put(90.66,0.67){\circle*{1.33}}
\put(90.66,27.67){\circle*{1.33}}
\put(120.00,27.67){\circle*{1.33}}
\put(120.00,0.67){\circle*{1.33}}
\put(105.33,13.33){\makebox(0,0)[cc]{$B_2$}}
\put(30.33,0.67){\line(1,0){31.00}}
\put(61.33,0.67){\line(0,0){0.00}}
\put(30.66,0.67){\circle*{1.33}} \put(61.00,0.67){\circle*{1.33}}
\put(85.66,6.33){\makebox(0,0)[cc]{$d_{0]12}$}}
\put(85.66,19.67){\makebox(0,0)[cc]{$d_{02]1}$}}
\put(104.66,32.67){\makebox(0,0)[cc]{${d}_2$}}
\put(125.00,13.33){\makebox(0,0)[cc]{${d}_1$}}
\put(104.66,-4.33){\makebox(0,0)[cc]{$d_{01]2}$}}
\put(30.33,-4.33){\makebox(0,0)[cc]{$d_{0]1}$}}
\put(61.00,-3.67){\makebox(0,0)[cc]{${d}_1$}}
\put(46.00,6.00){\makebox(0,0)[cc]{$B_1$}}
\end{picture}

\vspace{0.4in}
\begin{center}
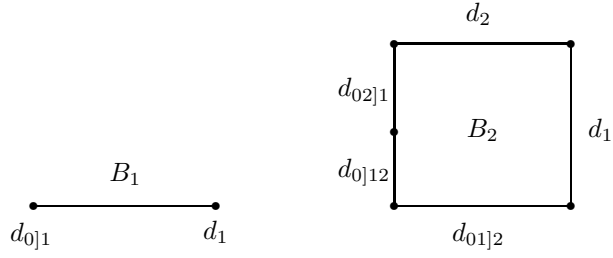
Figure 9: $B_{2}$ as a subdivision of $B_{1}\times I.$
\end{center}
\vspace{0.2in}

\unitlength=0.80mm \special{em:linewidth 0.4pt}
\linethickness{0.4pt}
\begin{picture}(133.67,39.66)
\put(118.00,1.33){\line(-1,0){28.67}}
\put(89.33,1.33){\line(-5,3){20.67}}
\put(68.66,13.67){\line(0,1){25.67}}
\put(68.66,39.33){\line(5,-3){21.00}}
\put(89.66,26.67){\line(1,0){28.00}}
\put(117.66,26.67){\line(-3,2){18.33}}
\put(99.33,39.00){\line(-1,0){30.67}}
\put(117.33,26.67){\line(0,-1){25.33}}
\put(89.33,1.33){\line(0,1){25.33}}
\put(99.66,39.00){\line(0,-1){11.33}}
\put(90.66,13.67){\line(1,0){9.00}}
\put(68.66,39.00){\circle*{1.33}}
\put(99.66,39.00){\circle*{1.33}}
\put(117.33,26.67){\circle*{1.33}}
\put(89.33,26.67){\circle*{1.33}}
\put(117.33,1.33){\circle*{1.33}} \put(89.33,1.33){\circle*{1.33}}
\put(68.66,13.67){\circle*{1.33}}
\put(99.66,13.67){\circle*{1.33}}
\put(79.66,19.67){\line(-5,3){11.00}}
\put(79.66,25.67){\line(5,-3){9.67}}
\put(68.66,13.67){\line(1,0){10.00}}
\put(80.66,13.67){\line(1,0){7.67}}
\put(79.66,6.67){\circle*{1.33}} \put(79.66,20.00){\circle*{1.33}}
\put(79.99,25.66){\circle*{1.33}}
\put(79.66,32.67){\circle*{0.00}}
\put(79.99,32.34){\circle*{1.33}}
\put(68.66,26.00){\circle*{1.33}}
\put(89.33,20.33){\circle*{1.33}}
\put(80.00,32.33){\line(0,-1){25.67}}
\put(133.67,13.00){\makebox(0,0)[cc]{$B_3$}}
\put(84.97,26.33){\makebox(0,0)[cc]{$_{d_{03]12}}$}}
\put(84.87,10.33){\makebox(0,0)[cc]{$_{d_{0]123}}$}}
\put(74.33,16.67){\makebox(0,0)[cc]{$_{d_{02]13}}$}}
\put(75.00,29.67){\makebox(0,0)[cc]{$_{d_{023]1}}$}}
\put(105.33,22.00){\makebox(0,0)[cc]{$_{d_{013]2}}$}}
\put(105.67,5.67){\makebox(0,0)[cc]{$_{d_{01]23}}$}}
\put(96.67,9.47){\makebox(0,0)[cc]{$_{d_{012]3}}$}}
\put(94.33,32.67){\makebox(0,0)[cc]{$_{{d}_{3}}$}}
\put(83.33,36.67){\makebox(0,0)[cc]{$_{{d}_{2}}$}}
\put(105.33,31.67){\makebox(0,0)[cc]{$_{{d}_{1}}$}}
\put(89.33,17.00){\line(1,0){28.33}}


\put(99.67,17.67){\line(0,1){8.33}}
\put(99.67,13.33){\line(0,1){3.33}}
\put(89.33,17.00){\circle*{1.33}}
\put(117.33,17.00){\circle*{1.33}}
\put(117.33,1.33){\line(-3,2){18.00}}
\end{picture}

\vspace{0.2in}

\begin{center}
Figure 10: $B_{3}$ as a subdivision of $B_{2}\times I.$
\end{center}
\vspace{0.2in}

The above face operators are connected to each other by the
following combinatorial  relations:  The relations between
$d_{A]M}$ and $d_{M_1|M_2}$ reflect the associativity of the
partition
 procedure, while the relations
  between   ${d}_i$ and  either $d_{A]M}$ or  $d_{M_1|M_2}$
 reflect the commutativity of the deleting and the partition
 procedures.
These relations together with those involving degeneracies are
incorporated in the singular permutocubes (see Example
\ref{Bsing}) which motivate structural relations for a
 {\em permutocubical set}  in the next section.

\section{Permutocubical sets}\label{trunc}
Here we give the formal definition of a {\it permutocubical set}.
The original motivation of that definition is
  the "twisted Cartesian product" of  cubical and
(multi) \linebreak permutahedral
 sets (see Section \ref{twistcart}). Note that
for applications here structural relations may be assumed modulo
degeneracies, so that we give the relations explicitly only for
face operators; the other relations involving degeneracies can be
written by examination of the universal example-the {\em singular
permutocubical set} if necessary.

Let $\mathcal{P}_{\ast\ast}^{0}(n)=\mathcal{P}_{\ast\ast}(n)\cup
\oo|\uu{n}$ with $\oo|\uu{n}\in \mathcal{P}_{0,n}^{0}(n).$

\begin{definition}
\label{permutocubical} A permutocubical set is a bigraded set
\[
{\mathcal B}=\{   {\mathcal B}_{p,q} \}_{p\geq 0,q\geq 1} \]
together with face and degeneracy operators

$$
\begin{array}{rlll}
{d}_i  &  : {\mathcal B}_{p,q}\rightarrow   {\mathcal B}_{p-1,q},
&
   $\newline$\vspace{0.1in}\\

 d_{A] M }
& : {\mathcal B}_{p,q}\rightarrow {\mathcal B}_{p-r, q+r-1},
    &  A\sm 0]M\in   {\mathcal P}^{0}_{p-r,r }(  {p}),

$\newline$\vspace{0.1in}
\\
d_{M_1|M_2 } & : {\mathcal B}_{p,q}\rightarrow {\mathcal
B}_{p,q-1}, & $\newline$\vspace{0.1in}
\\

\eta_{i}&
:{\mathcal B}_{p-1,q}\rightarrow {\mathcal B}_{p,q}, &  i\in \uu {p}, $\newline$\vspace{0.1in}\\

\varrho_{M_1|M_2}& :{\mathcal B}_{p,q-1}\rightarrow {\mathcal
B}_{p,q}, &
M_1|M_2\in   {\cal  P}_{*,*}({q}) , $\newline$\vspace{0.1in}\\
\end{array}
$$
such that
 for each $p\geq 0,$ the graded set
$$ \{   {\mathcal B}_{p,q};\,
 d_{M_1|M_2},\,\varrho_{M_1|M_2} \}_{q\geq 1} $$ is a  multipermutahedral
 set
 and the relations among face operators are:
\begin{equation*}
 \begin{array}{rlllll}
{d}_{i}{d}_{j}  &=&{d}_{j-1}{d}_{i},       \ \ \   \ \ \    i<j,$\newline$\vspace{0.1in}\\
{d}_{i}{d}_{A]M}&=&{d}_{A\sm  j ]M}{d}_{j}, \ \ \
i=I_{\uu{p}\setminus A}(j), \,\,\, i\in \uu
{p-r},$\newline$\vspace{0.1in}
\\
{d}_{i}{d}_{M_1|M_2}&=&{d}_{M_1|M_2}{d}_{i}, $\newline$\vspace{0.1in}\\

d_{K|L}d_{A]M}  &=& \left\{\begin{array}{lll}

d_{A]K\cap \uu r}d _{A\cup (K\cap \uu r)]L\cap \uu r},
        &       K|L \in
{\cal Q}_{1,q}({q+r}),\\

d_{A]M}d_{K+1-r| L+1-r} ,& K|L\in {\cal Q}_{r,1}({q+r}).
\end{array}
\right.
\end{array}
\end{equation*}
\end{definition}

\vspace{0.2in}

\begin{example}
\label{Bsing} For   a topological space $Y,$ define the
\uu{singular permutocubical} \linebreak \uu{complex} $\Sing^{B}Y$
as follows: Let
\[
{\Sing}^{B}_{p,q}Y=\{\text{continuous maps }B_{p}\times
P_{q}\rightarrow Y\}_{p\geq 0;q\geq 1},
\]
$B_p\times  P_{q}$ is a Cartesian product of the permutocube $B_p$
and the permutohedron $P_{q}$. Let
$$
\begin{array}{llll}
 {\delta}_{i}\times 1
&:
  B_{p-1}\times
 P_{q} \to  B_{p}\times
P_{q},\hfill  i\in \uu{ p},

      $\newline$\vspace{1mm} \\
\bar{\delta}_{A] M }&: B_{p-r}\times P_{q+r-1}
\xto{1\times{\Delta}_{r,q}}
 B_{p-r}\times P_{r}\times P_{q}\xto{\delta_{A] M}\times 1} B_{p}\times  P_{q},
          $\newline$\vspace{1mm}

\end{array}
$$
 be the  maps
in which ${\delta}_{i}$ and $\delta_{A] M}$ are the canonical
inclusions. Consider also
 the map
$$
\begin{array}{ll}
 \varsigma _i \x 1:B_{p}\times P_{q}
\to
 B_{p-1}\times P_{q}, &
i\in \uu{ p},

\end{array}
$$
where $\varsigma_i: B_{p}\to B_{p-1}$  is the projection that
identifies the faces $d_{\uu {p+1}\sm i] i}$ and $ {d}_i.$

Then for $f\in \Sing^{B}_{p,q}Y,$ define
$$
\begin{array}{rll}
{d}_i   & :  \Sing ^B _{p,q} Y\rightarrow \Sing ^B _{p -1, q} Y,
$\newline$\vspace{1mm}\\

d_{A] M}&: \Sing ^B _{p,q} Y\rightarrow
  \Sing ^B
_{p-r,q+r-1} Y,$\newline$\vspace{1mm}
\\
  $\newline$\vspace{1mm}

 \eta_{i}&: \Sing ^B _{p,q} Y\rightarrow  \Sing ^B
_{p+1,q}Y$\newline$\vspace{1mm}\\
\end{array}
$$
as compositions
$$
\begin{array}{rlll}
{d}_i(f) &= & f\circ  ( {\delta} _ i \times 1), $\newline$\vspace{1mm}\\

d_{  A]M}(f)&= &f\circ \bar{\delta}_{A]M},  $\newline$\vspace{1mm}\\

\eta_{i}(f)&=& f\circ (\varsigma_i\x 1),
\end{array}
$$
while define
\[
 d_{M_1|M_2} :
\Sing ^B _{p,q} Y\rightarrow  \Sing ^B _{p,q-1} Y \ \ \text{and} \
\
 \varrho_{M_1|M_2} :
\Sing ^B _{p,q} Y\rightarrow  \Sing ^B _{p,q+1} Y
\]
as $ d_{ M_1|M_2}=1\times d_{ M_1|M_2}$ and $\varrho_{
M_1|M_2}=1\times \varrho_{ M_1|M_2} $ to obtain  the singular
permutocubical set $(\Sing^{B}Y,\ {d}_{i},\, d_{A]M},d_{M_1|M_2},\
\eta_{i},\,\varrho_{M_1|M_2}).$

The singular permutocubical complex\, $\Sing^{B}Y$ determines the
singular (co)homo-logy of $Y$ in the following way: Form the
"chain complex"
\[
(C_{\ast}(\Sing^{B}Y),d)=\bigoplus_{\substack{p+q=n-1 \\n\geq1}
}(C_{n-1}(\Sing^{B}_{p,q}Y),d_{p,q}),
\]
where
$$d_{p,q}=\sum_{\substack{A\sm 0\,|M\in {\mathcal P}^{0}_{**}(p)\\ 1\leq i\leq p }}
 \left((-1)^{i+1}{d}_i+ (-1)^{\# A}\text{{\it shuff}}\,(A;M)\,d_{A]M}\right)+
(-1)^{p}d_q,$$
 $d_q:C_{q-1}(\Sing ^B_{p,q}Y)\rightarrow
C_{q-2}(\Sing ^B_{p,q-1}Y),$ the singular permutahedral
differential. Let
\[
C_{\ast}^{\boxminus}(Y)=C_{\ast}(\Sing^{B}Y)/DGN,
\]
where $DGN$ is the submodule of $C_{\ast}(\Sing^{B}Y)$ generated
by the degenerate elements of $\Sing^{B}Y,$ and obtain the {\it
singular permutocubical chain complex }
$(C_{\ast}^{\boxminus}(Y),d)$ of $Y.$

Now let $\varphi: B_n \to I^n$ be the  cellular projection defined
by the property that it maps homeomorphically    the faces $d_{\uu
n\sm i]i }(B_n)$ and $d_i(B_n)$ onto the faces $d^0_i(I^n)$ and
$d^1_i(I^n)$ respectively, $1\leq i\leq n$. Then the composition
of maps
\[
B_{p}\times P_{q}\xto{\varphi \x \rho } I^{p}\x
I^{q-1}=I^{p+q-1}\xto{\psi}\Delta^{p+q-1},
\]
clearly induces a composition of maps of graded sets
\[
 \Sing  Y\xto{\psi_{\ast}} \Sing^{I}Y\xto{(\varphi \x \rho)_{\ast}}
\Sing^{B}Y.
\]
 After the passage on the non-normalized chains (unless the last one) one gets a sequence of homomorphisms
\[
C_{\ast}( \Sing  Y)\rightarrow C_{\ast}(\Sing^{I}Y)\rightarrow
C_{\ast} ( \Sing^{B}Y)\rightarrow C_{\ast}^{\boxminus}(Y),
\]
whose composition is a chain map inducing  a natural isomorphism
$$H_{\ast}(Y)\approx H_{\ast}^{\boxminus}(Y)= H_{\ast}(C_{\ast
}^{\boxminus}(Y),d).$$
 Since the
diagonal on the permutocube
 constructed in
Section \ref{dsection}
 is compatible  with the A-W diagonal on the standard simplex under the above cellular
 projections,
  $H_{\ast}^{\boxminus}(Y)$
 determines the
singular cohomology ring of $Y$ as well.
\end{example}

Basic examples of a permutocubical set are provided in the next
section.

\section{Truncating twisting functions and twisted Cartesian
products}\label{twistcart}

As in \cite{KS1}
 in this
section we introduce the notion of a twisting function between
graded sets in which the domain and the target have face and
degeneracy operators of different types; but this time  we define
a {\it truncating twisting function} $\vartheta :Q\rightarrow
{\mathcal Z}$ from a cubical set $Q$ to a monoidal permutahedral
set ${\mathcal Z}$. For a permutahedral ${\mathcal Z}$-module with
action ${\mathcal Z}\times {\mathcal L}\rightarrow {\mathcal L},$
such $\vartheta $ defines a {\it twisted Cartesian product}
$Q\times _{\vartheta }{\mathcal L}$ as a {\it permutocubical set}.

These notions are motivated by the permutocubical set ${\bf P}Q$
(see below) which can be viewed as a twisted Cartesian product
determined by the canonical inclusion $ \vartheta :Q\rightarrow
{\bf \Omega} Q,$\ $x\mapsto \bar{x}$ of degree $-1,$ referred to
as the {\em universal truncating twisting function.}

The  geometrical interpretation of $\vartheta_U$
 answers to the truncation procedure that converts
 $I^n$ into $B_n$ mentioned in Section \ref{secpermutocube}. By this the
permutocube is thought of as
 a "twisted Cartesian product" of the cube and the permutohedron
(see Figures 11 and 12).

\vspace{0.2in}

\unitlength=0.700mm \special{em:linewidth 0.4pt}
\linethickness{0.4pt}
\begin{picture}(118.99,88.00)
\put(118.67,6.33){\line(-1,0){29.67}}
\put(89.00,6.33){\line(0,1){27.00}}
\put(89.00,33.33){\line(1,0){29.33}}
\put(118.33,33.33){\line(0,-1){27.00}}
\put(118.00,56.67){\line(-1,0){29.67}}
\put(88.33,56.67){\line(0,1){27.00}}
\put(88.33,83.67){\line(1,0){29.33}}
\put(117.67,83.67){\line(0,-1){27.00}}
\put(88.33,69.00){\circle*{1.33}}
\put(88.33,56.67){\circle*{1.33}}
\put(42.33,57.34){\line(0,1){12.33}}
\put(42.33,69.67){\circle*{1.33}}
\put(42.33,57.00){\circle*{1.33}}
\put(88.33,83.67){\circle*{1.33}}
\put(117.67,83.67){\circle*{1.33}}
\put(117.67,56.67){\circle*{1.33}}
\put(118.33,33.33){\circle*{1.33}}
\put(89.00,33.33){\circle*{1.33}} \put(89.00,6.33){\circle*{1.33}}
\put(118.33,6.33){\circle*{1.33}}
\put(48.00,63.67){\vector(1,0){21.67}}
\put(102.33,50.67){\vector(0,-1){13.33}}
\put(83.99,53.00){\makebox(0,0)[cc]{$_{0 ]1|2}$}}
\put(117.67,53.00){\makebox(0,0)[cc]{$_{0 ] 2}$}}
\put(117.67,87.33){\makebox(0,0)[cc]{$_{0 ]}$}}
\put(88.33,88.00){\makebox(0,0)[cc]{$_{0 ]1}$}}
\put(84.00,68.67){\makebox(0,0)[cc]{$_{0 ]2|1}$}}
\put(42.33,53.67){\makebox(0,0)[cc]{$_{1|2}$}}
\put(42.33,74.00){\makebox(0,0)[cc]{$_{2|1}$}}
\put(103.00,69.33){\makebox(0,0)[cc]{$B_2$}}
\put(103.00,18.00){\makebox(0,0)[cc]{$I^2$}}
\put(34.33,62.00){\makebox(0,0)[cc]{$P_2$}}
\put(88.66,2.00){\makebox(0,0)[cc]{$_{0 ]1|2}$}}
\put(118.00,2.00){\makebox(0,0)[cc]{$_{0 ]2}$}}
\put(88.66,36.33){\makebox(0,0)[cc]{$_{0 ]1}$}}
\put(117.67,36.00){\makebox(0,0)[cc]{$_{0 ]}$}}
\put(60.67,46.67){\makebox(0,0)[cc]{$\vartheta_U$}}
\put(107.67,43.67){\makebox(0,0)[cc]{$\varphi$}}
\put(78.67,23.00){\vector(-1,1){31.00}}
\put(58.00,67.00){\makebox(0,0)[cc]{$\delta_{0]\uu 2}$}}
\end{picture}

\vspace{0.2in}
\begin{center}
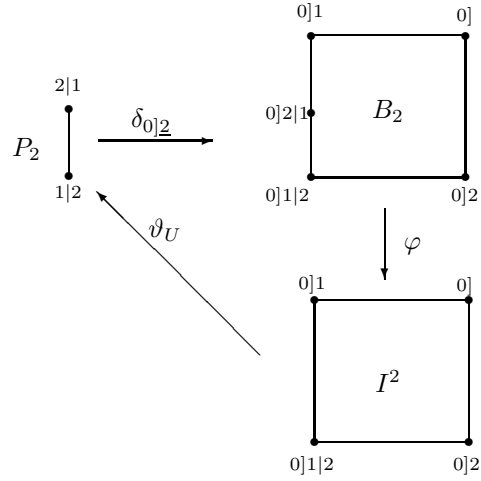
Figure 11: The universal truncating twisting function
$\vartheta_U:I^2\rightarrow P_2.$
\end{center}

\vspace{0.2in}

\unitlength=0.80mm \special{em:linewidth 0.4pt}
\linethickness{0.4pt}
\begin{picture}(133.67,112.67)
\put(117.00,6.33){\line(-1,0){28.67}}
\put(88.33,6.33){\line(-5,3){20.67}}
\put(67.66,18.67){\line(0,1){25.67}}
\put(67.66,44.33){\line(5,-3){21.00}}
\put(88.66,31.67){\line(1,0){28.00}}
\put(116.66,31.67){\line(-3,2){18.33}}
\put(98.66,44.00){\line(-1,0){30.67}}
\put(116.33,31.67){\line(0,-1){25.33}}
\put(88.33,6.33){\line(0,1){25.33}}
\put(98.66,44.00){\line(0,-1){11.33}}
\put(98.66,30.33){\line(0,-1){11.67}}
\put(67.66,18.67){\line(1,0){19.33}}
\put(89.66,18.67){\line(1,0){9.00}}
\put(67.66,44.00){\circle*{1.33}}
\put(98.66,44.33){\circle*{1.33}}
\put(116.33,31.67){\circle*{1.33}}
\put(88.33,31.67){\circle*{1.33}}
\put(116.33,6.33){\circle*{1.33}} \put(88.33,6.33){\circle*{1.33}}
\put(67.66,18.67){\circle*{1.33}}
\put(98.66,18.67){\circle*{1.33}}
\put(118.00,72.33){\line(-1,0){28.67}}
\put(89.33,72.33){\line(-5,3){20.67}}
\put(68.66,84.67){\line(0,1){25.67}}
\put(68.66,110.33){\line(5,-3){21.00}}
\put(89.66,97.67){\line(1,0){28.00}}
\put(117.66,97.67){\line(-3,2){18.33}}
\put(99.33,110.00){\line(-1,0){30.67}}
\put(117.33,97.67){\line(0,-1){25.33}}
\put(89.33,72.33){\line(0,1){25.33}}
\put(99.66,110.00){\line(0,-1){11.33}}
\put(90.66,84.67){\line(1,0){9.00}}
\put(68.66,110.00){\circle*{1.33}}
\put(99.66,110.00){\circle*{1.33}}
\put(117.33,97.67){\circle*{1.33}}
\put(89.33,97.67){\circle*{1.33}}
\put(117.33,72.33){\circle*{1.33}}
\put(89.33,72.33){\circle*{1.33}}
\put(68.66,84.67){\circle*{1.33}}
\put(99.66,84.67){\circle*{1.33}}
\put(79.66,90.67){\line(-5,3){11.00}}
\put(79.66,96.67){\line(5,-3){9.67}}
\put(68.66,84.67){\line(1,0){10.00}}
\put(80.66,84.67){\line(1,0){7.67}}
\put(79.66,77.67){\circle*{1.33}}
\put(79.66,91.00){\circle*{1.33}}
\put(79.99,96.66){\circle*{1.33}}
\put(79.66,103.67){\circle*{0.00}}
\put(79.99,103.34){\circle*{1.33}}
\put(68.66,97.00){\circle*{1.33}}
\put(89.33,91.33){\circle*{1.33}}
\put(30.99,97.00){\line(5,-3){9.67}}
\put(40.66,91.00){\line(0,-1){20.00}}
\put(30.66,77.33){\line(3,-2){10.00}}
\put(40.66,70.33){\circle*{1.33}}
\put(30.66,77.00){\circle*{1.33}}
\put(30.66,97.00){\circle*{1.33}}
\put(40.66,90.67){\circle*{1.33}}
\put(40.66,80.33){\circle*{1.33}}
\put(30.66,96.67){\line(0,-1){19.67}}
\put(30.66,87.00){\circle*{1.33}}
\put(83.00,66.00){\vector(0,-1){17.33}}
\put(62.00,48.33){\vector(-1,1){17.67}}
\put(46.00,74.33){\vector(1,0){34.33}} \
\put(87.66,68.67){\makebox(0,0)[cc]{$_{0 ]1|2|3}$}}
\put(120.33,69.00){\makebox(0,0)[cc]{$_{0 ]2|3}$}} \
\put(121.99,86.67){\makebox(0,0)[cc]{$_{_{0 ]3|2}}$}}
\put(121.66,97.67){\makebox(0,0)[cc]{$_{0]2}$}}

\put(84.99,86.03){\makebox(0,0)[cc]{$_{_{0]1|3|2}}$}} \

\put(95.33,91.00){\makebox(0,0)[cc]{$_{_{0 ]3|1|2}}$}} \
\put(75.39,89.33){\makebox(0,0)[cc]{$_{_{0 ]2|3|1}}$}}
\put(74.99,77.33){\makebox(0,0)[cc]{$_{_{0 ]2|1|3}}$}} \
\put(75.34,97.67){\makebox(0,0)[cc]{$_{_{0]3|2|1}}$}}
\put(63.34,84.66){\makebox(0,0)[cc]{$_{0 ]1|3}$}} \
\put(64.34,96.67){\makebox(0,0)[cc]{$_{_{0] 3|1}}$}}
\put(65.00,111.34){\makebox(0,0)[cc]{$_{0]1}$}}
\put(84.66,103.66){\makebox(0,0)[cc]{$_{_{0]2|1}}$}}
\put(99.66,112.67){\makebox(0,0)[cc]{$_{0 ]}$}}
\put(91.67,100.33){\makebox(0,0)[cc]{$_{0 ]1|2}$}}

\put(103.33,85.33){\makebox(0,0)[cc]{$_{0 ] 3}$}}
\put(80.00,103.33){\line(0,-1){25.67}}
\put(46.01,70.33){\makebox(0,0)[cc]{$_{1|2|3}$}}
\put(45.00,80.33){\makebox(0,0)[cc]{$_{_{1|3|2}}$}}
\put(46.00,90.67){\makebox(0,0)[cc]{$_{3|1|2}$}}
\put(25.33,97.33){\makebox(0,0)[cc]{$_{3|2|1}$}}
\put(27.00,87.00){\makebox(0,0)[cc]{$_{_{2|3|1}}$}}
\put(25.00,76.67){\makebox(0,0)[cc]{$_{2|1|3}$}}
\put(133.67,84.00){\makebox(0,0)[cc]{$B_3$}}
\put(133.33,22.67){\makebox(0,0)[cc]{$I^3$}}
\put(14.00,87.33){\makebox(0,0)[cc]{$P_3$}}
\put(88.33,2.34){\makebox(0,0)[cc]{$_{0 ]1|2|3}$}}
\put(119.33,2.00){\makebox(0,0)[cc]{$_{0 ]2|3}$}}
\put(103.00,20.33){\makebox(0,0)[cc]{$_{0] 3}$}}
\put(119.66,32.00){\makebox(0,0)[cc]{$_{0 ]2}$}}
\put(99.66,47.00){\makebox(0,0)[cc]{$_{0 ]}$}}
\put(64.00,44.01){\makebox(0,0)[cc]{$_{0] 1}$}}
\put(63.34,17.99){\makebox(0,0)[cc]{$_{0 ]1|3}$}}
\put(91.33,33.66){\makebox(0,0)[cc]{$_{0 ]1|2}$}}
\put(56.00,59.67){\makebox(0,0)[cc]{$\vartheta_U$}}
\put(88.33,57.33){\makebox(0,0)[cc]{$\varphi$}}
\put(63.00,69.90){\makebox(0,0)[cc]{$\delta_{0]\uu 3}$}}
\put(116.33,6.33){\line(-3,2){18.00}}
\put(117.67,72.00){\line(-3,2){18.00}}
\put(89.33,87.67){\line(1,0){28.00}}
\put(99.33,84.33){\line(1,0){0.33}}
\put(99.67,84.33){\line(0,1){3.00}}
\put(99.67,88.33){\line(0,1){8.67}}
\put(89.33,87.67){\circle*{1.49}}
\put(117.33,87.67){\circle*{1.33}}
\end{picture}

\vspace{0.2in}
\begin{center}
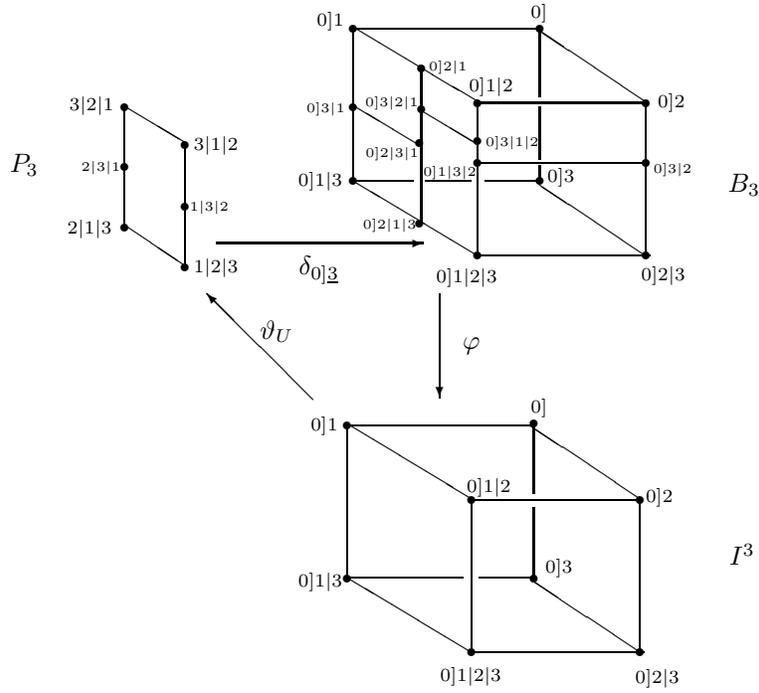
Figure 12: The universal truncating twisting function
$\vartheta_U:I^3\rightarrow P_3.$
\end{center}

\begin{definition}
Let $Q=(Q_n,d^0_i,d^1_i,\eta_i)$ be a 1-reduced cubical set and
${\mathcal Z}=({\mathcal Z}_{n},d_{M_1|M_2},\, \varrho_{M_1|M_2})$
be a monoidal permutahedral set. A sequence
$\vartheta=\{\vartheta_n\}_{n\geq 1} $ of degree $-1$ functions
$\vartheta_n :Q_n\to  {\mathcal Z}_n $
 is called a
truncating twisting function if it satisfies:
$$
\begin{array}{rllll}
\vartheta(a)
 & =&     e, & &
 a\in Q_1,
 $\newline$\vspace{1mm}\\

 d_{  M_1  | M_2  }  \vartheta (a) &
 =&
  \vartheta d^0_{M_2}  (a)
  \cdot
  \vartheta d^1_{M_1}(a),&    M_1|M_2\in   {\cal  P}_{*,* }( {n}), & a\in Q_n,
         $\newline$\vspace{1mm}
   \\
\varrho_{\uu{n}\sm i|i} \vartheta  (a)& =&         \vartheta
\eta_i(a),
& i\in \uu n.$\newline$\vspace{1mm}\\

\end{array}
$$
\end{definition}

Note that since the first condition above we in particular get

\[
 d_{i  | \uu n\setminus i} \vartheta  (a)   =  \vartheta
d^1_i(a) \ \ \text{and} \ \
 d_{\uu n\setminus i| i}\vartheta  (a)=
\vartheta d^0_i(a) \ \ \text{for} \ \  i\in \uu n   \ \ \text{and}
\ \    a\in Q_{n>0}.\]


\begin{remark}
By definition  a truncation twisting function commutes only with
the permutahedral degeneracy operator $\varrho_{\uu{n}\sm i|i},$
since it is in fact arisen  by the cubical degeneracy operator
$\eta_i$ (c.f. Remark \ref{lift}).
\end{remark}


We have the following

\begin{proposition}\label{universal}
Let $Q$ be a 1-reduced cubical set and ${\mathcal Z}$ be a
monoidal permutahedral set.
 A sequence $\vartheta=\{\vartheta_n\}_{n\geq
1} $ of degree $-1$ functions $\vartheta_n :Q_n\to  {\mathcal Z}_n
$
 is a
truncating twisting function  if and only if the  monoidal map
$f:{\bf \Omega} Q\to  {\mathcal Z}$ defined by $f(\bar{a}_1\cdots
\bar{a}_k)=\vartheta(a_1)\cdots \vartheta (a_k)$
 is a
map of   permutahedral sets.
\end{proposition}
\begin{proof}
Since $f$ is completely determined by its restriction to monoidal
generators, use the argument of verification of permutahedral
identities for a given single generator $\bar{\sigma}$ in
${\Omega} Q$ being equivalent to that of identities of the
universal truncating function $\vartheta_{U}:\sigma \rightarrow
\bar{\sigma}.$
\end{proof}
\begin{definition}\label{tmodule}
Let $Q=(Q_n,d^0_i,d^1_i,\eta_i)$ be a 1-reduced cubical set and
${\mathcal Z}=({\mathcal Z}_{n},d_{M_1|M_2},\,\varrho_{M_1|M_2})$
be a monoidal permutahedral set and ${\mathcal L}$ be a
permutahedral set with ${\mathcal Z}$-module structure. Let
 $\vartheta=\{\vartheta_n\}_{n\geq
1}, $\  $\vartheta_n :Q_n\to  {\mathcal Z}_n $ be a
 truncating twisting
function. The twisted Cartesian product
 $Q\times _{\vartheta}{\mathcal L}
 $
   is
the Cartesian product of sets $$Q\times {\mathcal L} =\{(Q\times
{\mathcal L})_{p,q} =\bigcup_{p\geq 0,q\geq 1} Q_{p} \times
{{\mathcal L}}_{q}\} $$ endowed with the face
 and degeneracy operators
${d}_i, d _{A]M}, d_{M_1|M_2},$
 $\eta_i,
 \varrho_{M_1|M_2}$  defined for $(a,b)\in Q_p\x {\mathcal L}_{q}$ by :
$$
\begin{array}{rlll}
{d}_i(a,b)&  = &  (d^1_i (a),\ b),   &   i\in \uu {p}, $\newline$\vspace{1mm}  \\

d_{ A  ] M  }(a,b)&
 =&  (d^0_{M}
(a),\  \vartheta d^1_{A}(a)\cdot b),  & A\sm 0\,|M\in   {\mathcal
P}^{0}_{*,* }(  {p}),
$\newline$\vspace{1mm}\\

d_{M_1|M_2}  (a,b)& =& (a,\ d_{M_1|M_2}(b)), &  M_1|M_2\in {\cal
P}_{*,* }({q})
$\newline$\vspace{2mm} \\

\eta_i(a,b)& =&  (\eta_i(a),\ b ),& i\in \uu {p+1}, $\newline$\vspace{1mm} \\

\varrho_{M_1|M_2}(a,b )&=& (a,\varrho_{M_1|M_2}(b) ),& M_1|M_2\in
{\cal P}_{*,* }({q+1}).  $\newline$\vspace{1mm} \\
\end{array}
$$

It is easy to check that $(Q\times _{\vartheta} {\mathcal L},\,
{d}_i, d _{A]M}, d_{M_1|M_2}, \eta_i,  \varrho_{M_1|M_2})$ is a
permutocubical set.

\end{definition}

\begin{remark}
Note that for the twisted Cartesian product $Q\times _{\vartheta}
{\mathcal L}$  we have  the following sequence of graded sets
$${\mathcal L}\xto{\iota} Q\times _{\vartheta} {\mathcal L}  \xto{\xi} Q$$
with $\iota (b)=(a_0,b)$ and $\xi(a,b)=a,\ a_0\in Q_0,\ (a,b)\in
Q\times {\mathcal L}.$
\end{remark}


\begin{example}
Let $M=\{  e_k \}_{k\geq 0}$ be the free minoid on a single
generator $e_1\in M_1$ with trivial permutahedral set structure
and let $\vartheta: Q\to  M$ be the sequence of constant maps
$\vartheta_n:Q_n\to  M_{n-1},\, n\geq 1.$ Then the twisted
Cartesian product $Q\times _{\vartheta} M$ can be thought of as a
permutocubical resolution of the 1-reduced cubical  set $Q.$
\end{example}


\subsection{The permutocubical set functor ${\bf P} Q$}\
The universal truncating twisting function $\vartheta_U$ defines a
special (acyclic)   twisted Cartesian product: Namely, we have

\begin{definition}
A functor from the category of 1-reduced cubical sets
 to the category of permutocubical sets
defined by $Q\to  Q\x_{\vartheta_U} {\bf \Omega} Q$ is the
universal
 permutocubical
functor and  denoted by  ${\bf P}.$
\end{definition}

\section{The  diagonal of permutocubes}\label{dsection}

Here we construct the explicit diagonal $\Delta_B:C_*(B_n)\to
C_*(B_n)\otimes C_*(B_n)$ for permutocubes which induces a
diagonal for a permutocubical set too.

\subsection{The orthogonal stream}\
Suppose that an $n$-dimensional  polytope $X$ is realized as a
subdivision of the cube $I^n$ so that each
   $m$-dimensional cell
$e\subset X,$ \   $0\leq m \leq n,$ is itself a subdivision of
$I^m$ ($I^m$ need not  be a face of $I^n;$ c.f. $B_n$). In
particular, we have an induced partial ordering on the set of all
vertices of $e$ defined by $x\leq y $ if there is an oriented
broken line from $x$ to $y.$
 For a cell $e^{\prime}\subset e,$ let $I^{m(e^{\prime})}\subset I^m$ be the
face of $I^m$ of the minimal dimension $m(e^{\prime})$ that
contains $e^{\prime}.$ Then we introduce the following
\begin{definition}
Let $ e\subset X$ be an $m$-cell and  $x\in e$ be a vertex. An
orthogonal stream $OS_x(e)$ of $x$ with the support $e$ is a pair
\[(U_x,V_x)= \left(\{ u_1,...,u_r \}\,,\{ v_1,...,v_s
\}\right)_{r,s\geq 1}\] of collections of
 faces of $e$  satisfying  the following conditions:

1. $\max u_r=x=\min v_1$  \ \  \text{and}\ \  $\dim u_r+ \dim
v_1=m; $

2. $    I^{m(u_i)}= I^{m(u_r)},\  \  \dim u_i=\dim u_r\ \
\text{and}\ \   \max u_i \leq x, \
     1\leq i\leq r;$

3. $  I^{m(v_j)}= I^{m(v_1)},\ \ \dim v_j=\dim v_1
 \  \ \ \text{and}\  \   \min v_j\geq x, \     1\leq j\leq s.   $

The union $\cup_{x\in e}SO_x(e)$ is denoted by $SO(e).$
\end{definition}

The pair
 $(u_r,v_1)\in OS_x(e)$
is referred to as
   {\it the strong
complementary pair} (SCP)  and  denoted by $(u_x,v_x);$ while a
pair $(u,v)\in OS_x(e)$ is referred to as {\it a complementary
pair }(CP) (compare, \cite{SU2}).

Clearly, any vertex $x\in e\subset B_n$
  uniquely
defines the SCP $(u_x,v_x)$ in  $OS_x(e),$ and, consequently, the
whole $OS_x(e)$ is uniquely determined by the vertex $x.$
 In particular, when $x$ coincides
with a vertex of $I^m$ then $\dim u_x=m(u_x)$ and $\dim
v_x=m(v_x),$ so that $U_x$ and $V_x$ actually lay on  orthogonal
faces of $I^m$ at the vertex $x.$

 For $B_n,$ an orthogonal stream $OS_x(e)
$ for each cell $e\subset B_n$ admits an explicit combinatorial
description. For example, for the top cell of $B_n$ we have the
following: Think of a vertex $x$=$0\,]\,x_1|...|x_k\in  B_n $ as
an ordered sequence of integers for  $1\leq k\leq n.$ Let $u_x$=$
A_0]A_1|...|A_p$ and $
 v_x$=$C_0]C_1|...|C_q$ be partitions   of  $\uu{n_0}$ in which $A_j$ for $1\leq j\leq p,$ and
 $C_i$ for $0\leq i\leq q,$ are the $(j+1)^{th}$ decreasing and the $(i+1)^{th}$
 increasing subsequence of maximal length of $x$ respectively (compare,
 \cite{SU2}); while for the vertex $x=0]\in B_n,$ let $(u_x,v_x)=\left(\uu{n_0}\,]\,,0]\right).$
 For
example,
  for $x $= $ 0]1|...|n\in B_n,$ one gets
$(u_x,v_x)$ = $\left(0]1|...|n\,, \uu{n_0}\,]\right);$  for $x $=
$ 0]2|1|3|6|5\in B_6,$ \, $(u_x,v_x)  $=$
\left(04]12|3|56\,,02]136|5\right).$

Next for a partition $a=A_0]A_1|...|A_{\ell}$ of $\uu{n_0},$ we
define the right-shift $R$ and the left-shift $L$ operators
respectively  as follows (compare, \cite{SU2}): For proper subsets
$M_i\subset A_i$ and $N_j\subset A_j,\,
 0  \leq i  < \ell, \,  0< j\leq \ell, $
 let
$$
\begin{array}{llll}

R_{M_i} ( a)  =A_0]A_{1}|\cdots|A_{i}\setminus M_i|A_{i+1}\cup
M_i|\cdots|A_{\ell} &  \ \ \text{for}\ \    \min M_i>\max A_{i+1}, \\

L_{N_j} ( a ) =A_{0}]A_1|\cdots|A_{j-1}\cup N_j|A_{j}\setminus
N_j|\cdots|A_{\ell} &  \ \ \text{for}\ \ \min N_{j}>\max A_{j-1}, \\
\end{array}
$$
where $R_{\varnothing}=Id =L_{\varnothing}.$ Then each CP
$(u,v)\in ( U_{x},V_{x})$ in the  orthogonal stream $OS_x(B_n)$
 can be  obtained from the SCP $(u_x,v_x)$ by
successive application of the above operators as
$$
(u ,v)=({R}_{M_{\ell-1}}\cdots R_{M_{1}} {R}_{M_0} (u_x) \,, \, {
L}_{N_1}\cdots L_{N_{\ell}}( v_x) )$$ for some $ \{M_i\}_{ 0 \leq
i< \ell}$ and $ \{N_j\}_{ 0<j\leq \ell}. $

For example, for the  vertex  $x=0]2|1|3|6|5,$  we obtain
\begin{multline*}OS_x(B_6)$=$(U_x,V_x)$=$\\
(\{0]12|34|56,\ 0]124|3|56,\ 04]12|3|56 \}\, , \{02]136|5],
\,023]16|5, 026]13|5, 0236]1|5 \} ).
\end{multline*}

 In particular, for the permutahedron $P_n$
the above description of CP's  agrees with that  by configuration
matrices in \cite{SU2}.

\subsection{The sign of $SO_x(B_n)$ }\
 The sign
 for a pair $(u,v)\in
OS_x(B_n)$ is deduced by motivation  that the cellular projection
$\varphi:B_n\to I^n$ preserves  diagonals.

First  for  a partition $a=A_0]A_1|...|A_p$ of $\uu{n_0}$ fix the
following signs:
\[
sgn_{1}(a)  =\left(  -1\right) ^{\epsilon_{1}}psgn (a) \text{\ \
and\ \ }\epsilon_{1}=\sum\nolimits_{i=1}^{p} i\cdot\#A_{p-i},
\]
\[
sgn_{2}(a)  =\left(  -1\right) ^{\epsilon_{2}}psgn(a) \ \
\text{and }\ \epsilon_{2}=\epsilon_{1}+\tbinom{p}{2},
\]
\[
rsgn (a)=(-1)^{\frac{1}{2}\left[ (\#A_{0}
)^{2}+\cdots+(\#A_{p})^{2}-(n+1)\right]}
\]
and $psgn (a)$ is the permutation sign $\{0,1,...,n\}\rightarrow
A_0]A_1|...|A_p .$

Now define
\[
sgn(u,v)=(-1)^{\binom{q+1}{2}}\text{ }rsgn (u)\cdot sgn_{1}
(v)\cdot sgn_{2}(u)\cdot sgn_{2}(u_x),
\]
where $q+1$ is the number of blocks in the partition $v.$

 Note that the above sign  agrees with  that of a CP from
$SO_x(P_n)$ established in \cite{SU2}.

\subsection{The diagonal of the permutocube}\
 By means of orthogonal streams we
 construct
an explicit diagonal for permutocubes as follows.

\begin{theorem}
The explicit diagonal of $B_n$
$$\Delta_B:C_*(B_n)\to  C_*(B_n)\otimes C_*(B_n)$$
is defined  for a cell $e \subset B_{n}$
 by
\begin{equation*}
 \Delta_B(e)=\sum_{ (e_1,e_2)\in
OS(e)}
 sgn(e_1,e_2)\, e_1\otimes e_2.
\end{equation*}

\end{theorem}

\vspace{5mm}
\begin{proof}
The proof is straightforward and analogous to that of Theorem 1 in
\cite{SU2}.
\end{proof}

In particular, in terms of orthogonal streams
 the diagonal $\Delta_P$ for permutahedra established in
\cite{SU2} can be formulated as follows.

\begin{theorem}\label{Pdiagonal}
The explicit diagonal of  $P_n$
$$\Delta_P:C_*(P_n)\to  C_*(P_n)\otimes C_*(P_n)$$
 is defined
for a cell $e  \subset P_n$
 by
\begin{equation*}
\Delta_P(e)=\sum _{(e_1,e_2)\in OS(e)}  sgn(e_1,e_2)\, e_1\otimes
e_2.
\end{equation*}
\end{theorem}
\vspace{5mm}

Below all components of $\Delta_B$ for the top cell of $B_n$ are
written  down for $n=1,2,3$ in which rows correspond to the
orthogonal streams.

\begin{example}
\end{example}
$$
\begin{array}{llllll}
\Delta_B\left(  1]\right)  =& \\

&\ \ \ 0]1  & \otimes   &  1] & & x=0]1 $\newline$\vspace{1mm} \\

&+\ 1] & \otimes  & 0]      &     & x=0].
\end{array}
$$

\begin{example}
\end{example}
$$
\begin{array}{llllll}
\Delta_B\left(  12]\right)  = &\\

& \ \  \ 0]1|2  & \otimes   &  12] & & x=0]1|2  $\newline$\vspace{1mm}\\

&- \ 0]12 & \otimes  & 2]1         &  & x=0]2|1 $\newline$\vspace{1mm}\\

&-\ (0]12 +2]1 )   &  \otimes   &   1]     &&     x=0]1 $\newline$\vspace{1mm}\\

&+   \ 1]2    & \otimes   &  2]   &&     x=0]2  $\newline$\vspace{1mm}\\

&+\ 12] & \otimes  &  0]       & &   x=0].

\end{array}
$$

\begin{example}
Up to sign, we have
\end{example}
$$
\begin{array}{lllll}
\Delta_B\left(  123]\right)  = &\\
 &\ \ \   0]1|2|3  & \otimes  &  123] &  x=0]1|2|3 \\

&+\  0]12|3   &\otimes  & 2]13           & x=0]2|1|3\\

&+\ 0]1|23  &  \otimes   &   13]2     &     x=0]1|3|2\\

&+  \  ( 0]12|3 + 0]1|23)  &  \otimes  &   3]12   &     x=0]3|1|2 \\

&+ \ 0]12|3  &  \otimes  &  (2]13+  23]1)     & x=0]2|1|3   \\

&+  \ 0]2|13 &  \otimes  & 23]1      & x= 0]2|3|1  \\

&+  \ ( 0]12|3+ 2]1|3)  &  \otimes  & 13]  &         x=0]1|3\\

&+   \  2]13   & \otimes  &  3]1          &     x=  0]3|1\\

&+    \ (  0]1|23+0]13|2+3]1|2)  & \otimes  &  12]             &         x=0]1|2\\

&+   \  (0]123+ 3]12)  &  \otimes     & 2]1                 & x=0]2|1\\

 &+   \  1]2|3  &  \otimes    &  23] &     x=0]2|3\\

& +\   1]23  &  \otimes   &   3]2 &     x=0]3|2\\

&+     \  (0]123+3]12+2]13+ 23]1) & \otimes  & 1]                &  x=0]1\\

&+ \  ( 1]23+13]2)   & \otimes    &  2] &     x=0]2\\

&+\    12]3   &   \otimes   &   3] &     x=0]3\\

&+    \  123]      &     \otimes  &  0]      &   x=0].

\end{array}
$$
\vspace{5mm}

\subsection{The diagonal on a permutocubical set}\
Now we use  the above combinatorial description  of an orthogonal
stream to define an explicit diagonal
 for a permutocubical set $
{\mathcal B}=\{ {\mathcal B}_{p,q}\}_{p\geq 0;q\geq 1}. $

A coproduct
$$\Delta:C_*({\mathcal B})\to  C_*({\mathcal B})\otimes C_*({\mathcal B})$$
is defined  for $a\in {\mathcal B}_{p,q}$ by
\begin{equation}\label{bformula}
\Delta(a)=\sum _{\substack{(u_1,u_2)\in OS(B_p)\\ (v_1,v_2)\, \in
OS(P_q)}} sgn(u_1,u_2)\cdot sgn(v_1,v_2)\cdot(-1)^{\epsilon}\,
  d_{u_1}d_{v_1}(a)\otimes d_{u_2}d_{v_2}(a),
\end{equation}
 $\epsilon=|d_{u_2}(a)||d_{v_1}(a)|.$

\vspace{0.2in}

\section{The permutocubical model for the path fibration}

 Let   $\Omega Y \xto{i} P Y\xto{\pi} Y$ be the Moore path
fibration on a topological space $Y.$ In \cite{Adams}
 Adams  constructed a dga map
$$
   \Omega C_*(Y)\to  C_*^{\Box}( \Omega Y)
$$
being a weak equivalence for a simply connected $Y,$ where $C_*$
denotes the singular simplicial chain complex, while in
\cite{Adams2} Adams and Hilton constructed a model for the path
fibration using the  singular cubical complex for each term of the
fibration. Here we obtain a natural combinatorial model for the
path fibration where for the base the singular cubical complex and
for the fibre the singular multipermutahedral complex are taken;
the total space in this case is modeled by the permutocubical set
being a twisted Cartesian product described in
 Section \ref{twistcart}. This model is naturally mapped into the
singular permutocubical complex of the total space. The chain
complex of the obtained model is a (comultiplicative) twisted
tensor product, while the Adams-Hilton model  is not. In
particular, the acyclic cobar construction $\Omega
\left(C^{\Box}_*(Y)\,;C^{\Box}_*(Y)\right)$ coincides with the
chain complex of the permutocubical set (compare, Theorem 5.1 in
\cite{KS1}).

For a space $Y,$ let  $\iota_0:\Sing ^{M}   Y \to \Sing ^{B}   Y $
be an inclusion of sets induced by the identification $P_{q}=B_0\x
P_{q}.$ Let denote  $\iota_*=\iota_0\circ i_*: \Sing ^{M}  \Omega
Y \to \Sing ^{B}  P Y  .$ Let
 $ (\varphi\x \rho)_{\ast}:\Sing ^{I}   Y \to \Sing ^{B}   Y$  be a natural map of graded sets
from Example \ref{Bsing}.
   Then we have the following theorem
(compare, \cite{Milgram}, \cite{CM}, \cite{Baues1}).

\begin{theorem}\label{cobar}
Let   $\Omega Y \xto{i} P Y\xto{\pi} Y$ be the Moore path
fibration.

(i)
 There are natural morphisms $\omega, p, (\varphi\x \rho)_{\ast} $
 such that
\begin{equation}\label{path}
\begin{CD}
\Sing ^{P}  \Omega Y @>\iota_*>>\Sing ^{B}  P Y     @>\pi _* >>\Sing ^B   Y    \\
@A \omega AA   @A  p AA       @A {(\varphi\x \rho)_{\ast}} AA  \\
{\bf \Omega} {\Sing^1} ^I Y @>\iota>>  {\bf P} {\Sing^{1}} ^I Y
@>\xi
>>{\Sing^{1}} ^I Y,
\end{CD}
\end{equation}
$(\varphi\x \rho)_{\ast}$ is a map of graded sets induced by
$\varphi\x \rho: B_{p}\x P_{n-p+1}\rightarrow I^n,$ while
 $p$ is a morphism  of   permutocubical  sets, and
 $\omega$ is a morphism  of monoidal  permutahedral  sets; $p$ and $\omega$ are
  homotopy equivalences whenever  $Y$ is
simply connected.

(ii)  The  chain complex $ C^{\diamondsuit}_*({\bf \Omega} {\Sing
^1}^I Y)$ coincides with the cobar construction $ \Omega
C^{\Box}_*(Y).$ \pagebreak

(iii) The  chain complex $C^{\boxminus}_*({\bf P}{\Sing ^1}^I Y)$
coincides with the acyclic cobar construction $\Omega \left(
C^{\Box}_*(Y)\,; C^{\Box}_*(Y)\right).$

\end{theorem}
\begin{proof}

(i). Morphisms  $p$ and $\omega$  are constructed simultaneously
 by induction on the dimension of
singular  cubes in ${\Sing ^{1}}^I Y.$ For $i=0,1$ and  $(\sigma
,e)\in {\bf P} {\Sing^{1}}^I Y,$  $ \sigma \in {\Sing^1}^I _iY,$
define $p\, (\sigma ,e)$ as
 the constant map $B_i \to  P Y$
to the base point $y,$  where $e$ denotes
 the unit of  the monoid ${{\bf \Omega}} {\Sing^1}^I Y$ (and of the monoid $\Sing
^P\Omega Y$ as well). Put $\omega (e)=e. $

Denote by ${\bf P}{\Sing ^1}^I _{(i,j)} Y$ the subset in ${\bf P}
{\Sing ^1}^I Y$
 consisting of the elements
$(\sigma,\tau )$
 with $|\sigma|\leq i $ and
$\tau \in {{\bf \Omega}} {\Sing ^1}^I _{(j)}Y,$ a
   submonoid in ${{\bf \Omega}} {\Sing ^1}^I Y$
 having (monoidal) generators $\bar \sigma$  with $|\bar{\sigma}|\leq j.$

Suppose by induction that we have constructed  $p$    and $\omega$
on  ${\bf P} {\Sing ^1}^I _{(n-1,n-2)} Y$ and
 ${{\bf \Omega}} {\Sing ^1}^I _{(n-2)} Y$ respectively such that
\begin{multline*}p\,(\sigma , \tau)= p\,(\sigma , e)\cdot \omega
(\tau)\ \ \text{and}\ \
  (\iota_*\circ \omega )(\bar \sigma)= p\,(
d_{0 ] \underline{r}} (\sigma ,e)), \,r=|\sigma|,\,1\leq
r<n,\end{multline*} where the $\cdot$ product is determined by the
action $PY \times \Omega Y\to  P Y.$ Let $\bar B_n \subset B_n$ be
the union of the all $(n-1)$-faces of $B_n$ except the
 $d_{0] \underline{n}}(B_n),$ and then for a singular
 cube
$\sigma: I ^n \to  Y$ define the map $\bar p: \bar B_n \to    P Y$
 by $$
 \bar p\, |_{ {d}_i(B_n)}= p\,({d}_i (\sigma , e) ),\ 1\leq
 i\leq n ,\
 \text{and}\ \
 \bar p\, |_{{d}_{A]M }(B_n)}= p\,({d}_{A]M }
  (\sigma , e) ),\ A]M\neq 0]\uu{n}.
  $$

Then
 the following  diagram commutes:
$$
\begin{CD}
\bar B_n @>\bar {p}_{\sigma}>>   P_{\sigma} Y   @> g_{\sigma} >>P   Y \\
@V \bar{i} VV          @V\pi _{\sigma} VV     @V\pi VV \\
B_n @>\varphi >>   I^n    @>\sigma >>   Y.
\end{CD}
$$
Clearly, $\bar{i}$ is a strong deformation retraction and we
define
 $ p\,(\sigma, e): B_n\to PY  $ as a lift of $\varphi.$
Define $ p\,(d_{0] \uu{n}}(\sigma ,e))=p(\sigma,e)|_{d_{0]
\uu{n}}(B_n)},$ and then $\omega (\bar \sigma)$ is determined by
$(\iota_* \circ \omega) (\bar \sigma)= p(\sigma, e)\circ\delta_{0]
\uu{n}}:P_n\to  B_n\to  PY.$

The proof of $p$ and $\omega$ being homotopy equivalences (after
the geometric realizations)
 immediately follows, for example,
from the  observation that $\xi$ induces a long exact homotopy
sequence. The last statement is a consequence of the following two
facts: (1) $|{\bf P} {\Sing ^{1}}^I Y|$ is contractible, (2) The
projection $\xi$ induces an isomorphism $\pi_* (|{\bf P} {\Sing
^1}^I Y |, |{{\bf \Omega}} {\Sing ^1}^I Y
|)\xto{\xi_*}\pi_*(|{\Sing^1}^I Y|).$

(ii)-(iii). This is straightforward.
\end{proof}

Thus, by passing on chain complexes in  diagram (\ref{path}) we
obtain the following comultiplicative model of $\pi$ formed by
dgc's (not necessarily coassociative ones).

\begin{corollary}
 For the path fibration   $\Omega Y \xto{i} P Y\xto{\pi} Y$
there is a comultiplicative model formed by  dgc's which  is
natural in $Y:$
\begin{equation*}\label{comultpath}
\begin{CD}
C^{\diamondsuit}_*( \Omega Y) @>\iota_*>>C^{\boxminus}_*(  P Y)
   @>\pi _* >>C^{\boxminus}_*(   Y)    \\
@A \omega_* AA   @A  p_* AA       @A {(\varphi\x \rho)_*} AA  \\
\Omega C^{\Box}_*( Y) @>>>   \Omega\left(C^{\Box}_*( Y)
;C^{\Box}_*( Y)\right)
 @>\xi_*>>C^{\Box}_*( Y).
\end{CD}
\end{equation*}
\end{corollary}


\section{Permutocubical  models for fibrations}

Here we prove the main result in this paper. Let $G$ be a
topological group, $F$ be a $G$-space $G\times F\rightarrow F$,
$G\rightarrow P\overset{\pi }{ \longrightarrow } Y$ be a principal
$G$-bundle and $F\rightarrow E\overset{ \zeta }{\longrightarrow }
Y$ be the associated fibration with the fiber $F$.
 Let $Q={\Sing^1}^I Y$, ${\mathcal Z}={\Sing}^M G$ and ${\mathcal
L}={\Sing}^M F$. The group operation $G\times G\to G$ induces  the
structure of a monoidal multipermutahedral set on ${\mathcal Z},$
and the action $G\times F\to F$ induces ${\mathcal Z}$-module
structure ${\mathcal Z}\times {\mathcal L}\to {\mathcal L}$ on
${\mathcal L}$ (c.f. Example \ref{monoidal}).

\begin{theorem}\label{percubmodel}
 The principal $G$-fibration $G\rightarrow P\overset{\pi }{
\longrightarrow }Y$ determines a truncating twisting function
$\vartheta :{\Sing^1}^I Y\to {\Sing}^M G$ such that twisted
Cartesian product ${\Sing^1}^I Y \times_{\vartheta } {\Sing}^M F$
models  the total space $E$ of the associated fibration
$F\rightarrow E\overset{\zeta }{ \longrightarrow }Y,$ that is,
there exists a permutocubical map
$$
{\Sing^1}^I Y\times_{\vartheta } {\Sing}^M F\to {\Sing^B E}
$$
inducing  homology isomorphism.
\end{theorem}

\begin{proof}
Let $\omega:{\bf \Omega} Q\to {\Sing}^M\Omega Y$ be the map of
monoidal multipermutahedral  sets  from Theorem \ref{cobar}.  By
Proposition \ref{universal} $\omega$ corresponds to  a truncating
twisting function $\vartheta':Q={\Sing^1}^I Y
\xto{\vartheta_U}{\bf\Omega} Q= {\bf \Omega} {\Sing^1}^I Y
\xto{\omega} {\Sing}^M\Omega Y.$
 Composing $\vartheta'$ with the map of monoidal multipermutahedral
sets ${\Sing}^M\Omega Y\to {\Sing}^M G={\mathcal Z} $ induced by
the canonical map $\Omega Y\to G$ of monoids we obtain a
truncating twisting function $\vartheta :Q\to {\mathcal Z}$. The
resulting twisted Cartesian product ${\Sing^1}^I Y
\times_{\vartheta }{\Sing}^M F$ is a permutocubical model of $E.$
Indeed, we have the canonical equality
$$Q \times_{\vartheta} {\mathcal L}=
 (Q\times_{\vartheta}{\mathcal Z}) \times {\mathcal L}/\sim,$$
where $(xg,y)\sim(x,gy).$ Next the argument of the proof of
Theorem \ref{cobar} gives a permutocubical map
$f':Q\times_{\vartheta_U} {\bf \Omega} Q\to \Sing^B P$ preserving
the actions of  ${\bf \Omega} Q$ and ${\mathcal Z}.$ Hence, this
map extents to a permutocubical map $f:Q\times_{\vartheta}
{\mathcal Z}\to \Sing^B P $ by $f(x,g)=f'(x,e)g.$ The map
\begin{multline*}
   (Q\times_{\vartheta} {\mathcal Z}) \times  {\mathcal L}\xto{f\times 1} \Sing^B
P\times {\mathcal L}\xto {\lambda} \Sing^B
 (P\times F),\\
 \lambda (h_1,h_2) =(h_1\times h_2)\circ (1_B\times \Delta_{r,s}),
 \end{multline*}
 induces the map of permutocubical sets
$$ {\Sing^1}  ^I Y \times_{\vartheta} \Sing^M F\to  \Sing^B E$$
as desired.
\end{proof}

\vspace{0.1in}

 For convenience, assume that $Q, {\mathcal Z}$ and ${\mathcal L}$ are as
in  Definition \ref{tmodule}. On the chain level  a truncating
twisting function $\vartheta$ induces  the twisting cochains
 $\vartheta_*:C^{\Box}_*(Q)\to C^{\diamondsuit}_{*-1}({\mathcal Z})$ and
$\vartheta^*:C_{\diamondsuit}^*({\mathcal Z})\to
C_{\Box}^{*+1}(Q)$  in the standard sense (\cite{Brown},
\cite{Berika3}, \cite{Gugenheim}).
 It is straightforward to verify that  the
following equality holds:
\begin{equation}\label{coalgebra}
C^{\boxminus}_*(Q\times _{\vartheta }{\mathcal
L})=C^{\Box}_*(Q)\otimes
_{\vartheta_*}C^{\diamondsuit}_*({\mathcal L}),
\end{equation}
and, consequently, the obvious injection
\begin{equation}
\label{algebra} C_{\boxminus}^*(Q\times _{\vartheta }{\mathcal
L})\supset C^*_{\Box}(Q)\otimes
_{\vartheta^*}C_{\diamondsuit}^*({\mathcal L})
\end{equation}
of dg modules (which is an equality  if the graded sets are of
finite type).

The permutocubical structure of $Q\times _{\vartheta }{\mathcal
L}$ induces a dgc sturcture on $C^{\boxminus}_*(Q\times
_{\vartheta }{\mathcal L}).$ Transporting this structure (diagonal
(\ref{bformula})) on the right-hand side of (\ref{coalgebra}) we
obtain a {\it comultiplicative} model of $C^{\Box}_*(Q)\ox
_{\vartheta }C_{\ast}^{\diamondsuit}({\mathcal L})$
  of our fibration. Dually,
$C_{\boxminus}^*(Q\times _{\vartheta }{\mathcal L})$ is a dga, so
a dga structure (a multiplication)  arises  on the right-hand side
of (\ref{algebra}) and  we obtain a {\it multiplicative} model
$C_{\Box}^*(Q)\ox _{\vartheta }C^{\ast}_{\diamondsuit}({\mathcal
L})$
 of
our fibration.

Below we describe these structures (the comultiplication on
$C^{\Box}_*(Q)\ox _{\vartheta }C_{\ast}^{\diamondsuit}({\mathcal
L})$
 and the
multiplication on $C_{\Box}^*(Q)\ox _{\vartheta
}C^{\ast}_{\diamondsuit}({\mathcal L})$
 in terms of certain (co)chain operations that form
a {\it Hirsch (co)algebra} structure on the (co)chain complex of
$Q$.

\subsection{The canonical  Hirsch  algebra
structure on $C_{\Box}^*(Q)$}\   Consider the equality
\begin{equation*}
C^{\diamondsuit}_*({\bf \Omega} Q)=\Omega C^{\Box}_*(Q)
\end{equation*}
from Theorem \ref{cobar}.
 As before, the permutahedral structure of ${\bf \Omega} Q$ induces a
coproduct  on $C_{\ast }^{\Box }({\bf \Omega} Q)$ (\cite{SU2});
consequently, this structure also appears on the right-hand side
of the above equality, so that the cobar construction $\Omega
C_{\ast }(Q)$ becomes a dg Hopf algebra.

To describe the above coproduct in terms of generators (singular
cubes) we need the following combinatorial analysis of the
diagonal $\Delta_P$ on permutahedra (compare \cite{Baues3},
\cite{KS1}).

Given an ordered subset $B\subset \mathbb{N}\cup 0$ and $a,b\in B$
with $a<b,$ let $[a\cdots b]=\{x\in B\,|\,a\leq x \leq b \}$ be a
{\it block}; for $A=(a_1<\cdots <a_k)\subset B,$ let
\[J_B(A)=
[a_1\cdots a_2]\cdot[a_2\cdots a_3]\cdots [a_{k-1}\cdots a_k]\]
 be a sequence of blocks; then
$\bar{J}_B(A) $ can be thought of as a generator of the monoid
${\bf \Omega} Q,$ i.e., $\bar{J}_B(A)\in \bar{Q}$
 (recall Proposition 3.1 from \cite{KS1} that we have the
correspondence between sequences of such blocks and compositions
of cubical face operators).
  On the other hand, $\Delta_P$  can be expressed on  $\uu{n}$
 as
\[\Delta_P(\uu{n})=\sum_{(u,v)\in OS(P_n)}\bar{J}_{\tilde{u}_1}(\tilde{u}_2)\cdots
 \bar{J}_{\tilde{u}_{p-1}}(\tilde{u}_p)
\ox \bar{J}_{\tilde{v}_1}(\tilde{v}_2)\cdots
\bar{J}_{\tilde{v}_{q-1}}(\tilde{v}_q),
\]
for
 $
(\tilde{u}_i\,,\tilde{v}_j)=(u_i\cup 0\cup (n+1)\,,v_j\cup 0\cup
(n+1)),$ \,  $( u_i,v_j)=(A_i\cup\cdots \cup A_p\,, C_j\cup\cdots
\cup C_q)_{ 1\leq i\leq p, 1\leq j\leq q},$\,
$(u,v)=(A_1|...|A_p\,,C_1|...|C_q).$ Consider
 the
identification ${J}_{w_{i}}(w_{i+1})=d^0_{w_{i+1}}d^1_{\uu{n}\sm
w_i}([01\cdots n+1]),\, w=u,v,$ to  obtain  the following formula
for the coproduct $\Delta : \Omega C^{\Box}_*(Q)\to \Omega
C^{\Box}_*(Q)\otimes \Omega C^{\Box}_*(Q)$: For a generator
$\sigma\in C^{\Box}_n(Q),$ let
\begin{equation}
\label{cohirsch}
 \Delta(\bar{\sigma})=\sum_{(u,v)\in
OS(P_n)}sgn(u,v) \left(\bigotimes_{i=1}^{p}\overline{
d^0_{u_{i+1}}d^1_{\uu{n}\sm u_i}(\sigma)}\right)\otimes
\left(\bigotimes_{i=1}^{q} \overline{d^0_{v_{i+1}}d^1_{\uu{n}\sm
v_i}(\sigma)}\right).
\end{equation}
 Note that since $Q$ is assumed to be 1-reduced, the
image $\overline{ d^0_{w_{i+1}}d^1_{\uu{n}\sm w_i}(\sigma)}$ of a
1-dimensional face ${ d^0_{w_{i+1}}d^1_{\uu{n}\sm w_i}(\sigma)}$
for $w=u,v,$ is the unit in $ \Omega C^{\Box}_{\ast }(Q)$ and
hence can be omitted.

Actually the diagonal consists of {\it components}
\begin{equation*}
E^{p,q}=pr\circ \Delta :C^{\Box}_{\ast }(Q)\rightarrow \Omega
C^{\Box}_{\ast }(Q)\otimes \Omega C^{\Box}_{\ast }(Q)\rightarrow
C^{\Box}_{\ast }(Q)^{\otimes p}\otimes C^{\Box}_{\ast
}(Q)^{\otimes q},\ p,q\geq 1,
\end{equation*}%
where $pr$ is the obvious projection.

The basic component $E^{1,1}$ is formed by those pairs $(u,v)\in
OS(P_n)$ in which all but one pair satisfy  $(\#A_i,\#C_j)=(1,1);$
this component  is a chain operation  dual to the cubical version
of Steenrod's $\smile _{1}$-product.

Dualizing the operations $E^{p,q},$ we obtain
 the sequence of cochain operations
$$
\{E_{p,q}: C^*_{\Box}(Q)^{\otimes p}\otimes C^{*}_{\Box}(Q
)^{\otimes q} \to  C^*_{\Box}(Q)\}_{p+q\geq 0},
$$
which define a
 multiplication on the bar construction
$BC_{\Box}^*(Q)\otimes BC_{\Box}^*(Q)\to\linebreak
BC_{\Box}^*(Q).$
 These cochain operations  form on $C^*_{\Box}(Q)$ the
structure of a {\it Hirsch algebra} (see the next section).

 The operations $E_{p,q}$ are
 restrictions of more general cochain  operations that arise on
$\bar{C}_{\Box}^*( Q)$ (the non-normalized chains) for a based
space $Y$ which is not necessarily 1-connected. In this case, for
$Q=\Sing ^{I}Y,$ we have the operations
$$
\{E_{p,q}: \bar {C}^*_{\Box}(Q)^{\otimes p}\otimes \bar
{C}^{*}_{\Box}(Q )^{\otimes q} \to  \bar
{C}^*_{\Box}(Q)\}_{p,q\geq 1}
$$
 given by the following explicit formulas:
 For $ a_i\in \bar{C}^{\geq 2}( Q)$, $b_j\in \bar{C}^{\geq 2}(
Q)$,\, $1\leq i\leq p$, $1\leq j\leq q$, let
$${E}_{p,q}(a_1,...,a_p;\, b_1,...,b_q)=
\sum_{s\geq p;\,t\geq q} \bar{E}_{s,t}(\epsilon^1 ,a_1,
\epsilon^1,...,
     \epsilon^1,a_p,  \epsilon^1;\,
\epsilon^1 ,b_1, \epsilon^1,...,
     \epsilon^1,b_q,  \epsilon^1
     ),$$
 $\epsilon^1\in \bar{C}^1( Q)$ is the generator represented by
 the constant singular 1-cube at the base point $I\rightarrow y\in Y$
 and the operations $\bar{E}_{s,t}$ are  defined
for  $ a_i\in \bar{C}^{k_i}_{\Box}(Q),\  b_j\in
\bar{C}^{r_j}_{\Box}(Q),\, \sigma \in Q_n,$
   by
$$\bar{E}_{s,t}(a_1,...,a_s\,;\, b_1,...,b_t)=c\in \bar{C}^{n}_{\Box } (Q), $$
\[
c(\sigma)=  \sum_{\substack{u\in {\mathcal
P}_{k_1,...,k_s}(n)\\v\in{\mathcal P}_{r_1,...,r_t}(n)\\(u,v)\in
OS(P_n) }}sgn(u,v)\,
 a_1(\sigma_1)\cdots a_s(\sigma_s)\cdot b_1(\sigma'_1)\cdots
 b_t(\sigma'_t),\]

\[ \sigma_i=d^0_{u_{i+1}}d^1_{\uu{n}\sm u_{i}}(\sigma),\, 1\leq
i\leq s,\ \ \
 \sigma'_j=d^0_{v_{j+1}}d^1_{\uu{n}\sm
v_{j}}(\sigma),\, 1\leq j\leq t,
\]
 where
$(u_i,v_j)$ is as in $\Delta_P(\uu{n})$ above,
 and where
$\bar{E}_{s,t}(a_1,...,a_{s}; b_1,...,b_{t})=0$ otherwise.

Thus, the above formula for $p,q=1$ defines ${E}_{1,1}$
 as  the cubical version of
 Steenrod's cochain $\smile_1$-operation without any restriction on $Y.$

\begin{remark}
 The  operations  $\{E^{p,q}\}$ on $C^{\Box}_{*}(Q)=\Omega
C_{*}(X),$
  $Q={\bf \Omega}
\Sing^2X,$ in fact have the form
$$
E^{p,q}=\sum \Delta^{p-1}_{E}\otimes \Delta^{q-1}_{E}
$$
where $\Delta^k_{E}:\Omega C_{*}(X)\to \Omega C_{*}(X)^{\otimes
k+1}$ is the $k$-th iteration of the comultiplication
$\Delta_E:\Omega C_{*}(X)\to \Omega C_{*}(X)\otimes \Omega
C_{*}(X)$ being itself induced by the homotopy G-coalgebra
structure $\{E^{k,1}\}$ on $C_{*}(X)$ (c.f. \cite{KS1}).
\end{remark}

\subsection{Twisted multiplicative model for a fibration}\
Next we further explore the twisted Cartesian product
 $Q\times _{\vartheta }{\mathcal L}$. To describe the corresponding
coproduct and product on the right-hand sides of (\ref{coalgebra})
and (\ref{algebra}) respectively, it is very convenient to express
the diagonal $\Delta_B$ in terms of
 combinatorics of the cubes and permutahedra.
  Namely, for the top cell $\uu{n_0}$ of $B_n,$
 let
\begin{multline}\label{deltaB}
\Delta_B(\uu{n_0})= [01\cdots n+1]\,\ox \,[0,n+1]+\\
\sum_{\substack{ \uu{n}\sm N|N\in {\mathcal P}^0_{n-s,s}(n)
\\(u,v)\in OS(P_{s})}}
{J}_{u_0}( u_1^{\prime})\cdot \bar{J}_{u^{\prime}_1}(u^{\prime}_2)
\cdots \bar{J}_{u^{\prime }_{p-1}}(u^{ \prime}_p)
 \ox
{J}_{v^{ \prime}_1}(v^{ \prime}_2)\cdot
 \bar{J}_{v^{
\prime}_2}(v^{ \prime}_3)\cdots \bar{J}_{v^{ \prime}_{q-1}}(v^{
\prime}_q),
\end{multline}
for
   $u_0=\{01\cdots n+1\},$ $( u^{\prime }_i,v^{
\prime}_j)=(I^{-1}_{N^{\prime}}(\tilde{u}_i)\, ,
I^{-1}_{N^{\prime}}(\tilde{v}_j)), $ $N^{\prime}=N\cup 0\cup
(n+1),$ where
 $(\tilde{u}_i,\tilde{v}_j)_{1\leq i\leq p, 1\leq j\leq q}$ is  as
in $\Delta_P(\uu{s})$ above. In particular,
 the summand $ [01\cdots n+1]\ox [0,n+1]$ is a primitive component of the
diagonal, while the second one is obtained by $N=\uu{n},\,
(u,v)=(1|2|...|n\,, \uu{n}),$ and is equal to
\begin{multline*}
[01][12]...[n,n+1]\cdot \overline{[012][23]...[n,n+1]}\cdot
\overline{[023][34]...[n,n+1]} \cdots \overline{[0,n,n+1]}\,
 \ox \\ \ox [01...n+1].
\end{multline*}

\begin{remark}
Note that we abuse the notation when we mean under ${[01\cdots
n+1]}$ an $n$-permutocube, since this notation was used for
combinatorial description of the $n$-cube in \cite{KS1}.
Accordingly here $\overline{[01\cdots n+1]}$ corresponds to the
$(n-1)$-permuta-\linebreak hedron.
\end{remark}

Furthermore, the action ${\mathcal Z}\times {\mathcal
L}\rightarrow {\mathcal L}$ induces a comodule structure
$\Delta_{{\mathcal {\mathcal L}}} :C_{\diamondsuit}^*({\mathcal
L})\to C_{\diamondsuit}^*({\mathcal Z})\otimes
C_{\diamondsuit}^*({\mathcal L})$ and it is not hard to see that
the permutocubical multiplication of (\ref{algebra}) can be
expressed by this comodule structure, the diagonal (\ref{deltaB}),
the twisting cochain $\vartheta ^{\ast }$, and the operations
$\{E_{p,q}\}_{p,q\geq 1}$ by the following formula:
 Let $a_1\otimes m_1,a_2\otimes
m_2\in C_{\Box}^*(Q)\otimes
_{\vartheta^*}C_{\diamondsuit}^*({\mathcal L})$ and
 $\Delta_L^k: C_{\diamondsuit}^*({\mathcal L})\rightarrow
  C_{\diamondsuit}^*({\mathcal Z})^{\otimes k}\otimes
C_{\diamondsuit}^*({\mathcal L}) $ be the iterated
$\Delta_{{\mathcal L}} $ with $\Delta_{{\mathcal
L}}^0=\Id:C_{\diamondsuit}^*({\mathcal L})\rightarrow
C_{\diamondsuit}^*({\mathcal L});$  let $\Delta_{{\mathcal L}}^p
(m_1)=\sum c_1^1\otimes \ldots \otimes c_1^p\otimes m_1^{p+1}$ and
$\Delta_{{\mathcal L}}^{q-1} (m_2)=\sum c_2^1\otimes \ldots
\otimes c_2^{q-1}\otimes m_2^{q};$ then
\begin{multline}\label{tformula}
 \mu((a_1\otimes m_1)\otimes (a_2\otimes m_2)) =\\
 \sum_{p\geq 0;\ q\geq 1}
 (-1)^{\epsilon}
a_1  E_{p,q }( \vartheta(c_1^1),... , \vartheta(c_1^p); a_2,
       \vartheta(c_2^{1}),...,\vartheta(c_2^{q-1})) \otimes
 m_1^{p+1}m_2^{q},
\end{multline}

$\epsilon= |m_1^{p+1}|(|a_2|+|c_2^{1}|+\cdots+|c_2^{q-1}|).$

\begin{corollary}\label{twisted}
Under the circumstances of Theorem \ref{percubmodel}, the twisted
differential $d_{\vartheta}$ and multiplication $\mu $ turn the
tensor product $C_{\Box}^*(Z)\otimes C_{\diamondsuit}^*(F)$
 into a dga
$(C_{\Box}^*(Z)\otimes
C_{\diamondsuit}^*(F),d_{\vartheta},\mu_{\vartheta})$
 weakly equivalent to the dga $C_{\diamondsuit }^{\ast
}(E)$.
\end{corollary}

\begin{corollary}

There exists on the acyclic bar construction
 $B(C^*_{\Box}(Z);C^*_{\Box}(Z))$
 the following  multiplication:
 For
 $a=a_0\otimes
[\bar{a}_1|\cdots |\bar{a}_n],\ \  b=b_0\otimes [\bar{b}_1|\cdots
|\bar{b}_m],\ \ a_i,b_j \in C_{\Box}^*(Z),\ 0\leq i\leq n,\ 0\leq
j\leq m, $ let \vspace{1mm}
\begin{equation}\label{tbformula}
 ab=  \sum _{p\geq 0;\,q\geq 1}
 (-1)^{\epsilon}
a_0  E_{p,q }( a_1, ... , a_p\,; b_0,   b_{1},...,b_{q-1}) \otimes
[\bar a_{p+1 }|\dotsb |\bar a_n]\circ [\bar b_{q}|\dotsb |\bar
b_m],
\end{equation}
$\epsilon= \left(|\bar{a}_{p+1}|+\cdots+|\bar{a}_{n}|\right)
\left(|b_0|+|\bar{b}_{1}|+\cdots+|\bar{b}_{q-1}| \right).$
\end{corollary}
\begin{proof}
Take ${\mathcal Z}={\mathcal L}={\bf \Omega} Q$. Then the
multiplication (\ref{tformula}) looks as (\ref{tbformula}).
\end{proof}

Using  the fact that $BC^*(Y)$ has an  associative multiplication
\cite{KS1} we canonically introduce on the acyclic bar
construction
 $B(BC^*(Y);BC^*(Y))$
  the multiplication by (\ref{tbformula}) that agrees with the one
 on the double bar construction  $BB C^*(Y)$ \cite{SU2}.

\section{Twisted tensor products for  Hirsch algebras}

The notion of  Hirsch (co)algebra naturally generalizes that of a
homotopy G-(co)algebra. We generalize the theory of {\it
multiplicative} twisted tensor products for homotopy G-algebras,
and, consequently, for commutative dga's \cite{KS1}. Namely we
define a twisted tensor product with both twisted differential and
{\it twisted multiplication } inspired by  formulas
(\ref{tformula}) and (\ref{tbformula}) established in the previous
section.

 Let $A$ be a dga and consider the dg module
$({Hom}(BA\otimes BA,A),\nabla )$ with differential $\nabla $. The
$\smile $-product induces a dga structure (the tensor product
$BA\otimes BA$ is a dgc with the standard coalgebra structure).

\begin{definition}
A Hirsch algebra is a  1-reduced associative dga $A$ eqwipped with
multilinear maps
$$
E_{p,q}:A^{\otimes p}\otimes A^{\otimes q}\to A, \ p,q\geq 0,\
p+q> 0,
$$
satisfying  the following conditions:
\begin{enumerate}
\item[(i)]
 $E_{p,q}\ \  \text{is of degree}\ \ 1-p-q$;

 \item [(ii)] $E_{1,0}=Id=
E_{0,1}\ \ \text{and}\ \ E_{p>0,0}= 0=E_{0,q>0};$

 \item [(iii)] The
homomorphism $E:BA\otimes BA \to A$ defined by
$$
E([\bar a_1|\dotsb | \bar a_p]\otimes [\bar b_1|\dotsb | \bar
b_q])= E_{p,q}(a_1, ... , a_p;b_1,...,b_q)
$$
is a twisting element in the dga $(Hom(BA\otimes BA,A),\nabla )$,
i.e., it satisfies $\nabla E =- E\smile E.$
\end{enumerate}
\end{definition}
 \noindent Condition (iii) implies that $\mu_E$ is a
chain map; thus $BA$ becomes a dg Hopf algebra with not
necessarily associative multiplication $\mu_E$ (c.f. \cite{GJ},
\cite{Voronov2}).
 Condition (iii)  can be rewritten in terms of components
$E_{p,q}.$  In particular, the operation $E_{1,1}$ satisfies
conditions similar to  Steenrod's $\smile_1$ product:
\[
dE_{1,1}(a;b)-E_{1,1}(da;b)+(-1)^{|a|}E_{1,1}(a;db)= (-1)^{|a|}ab
-(-1)^{|a|(|b|+1)}ba,
\]
so it measures the non-commutativity of the product of $A$ (thus,
a Hirsch algebra with $E_{p,q}=0$ for $p,q\geq 1$ is just a
commutative dga).

The dual notion is that of a Hirsch coalgebra. For a Hirsch
coalgebra  \linebreak $ \left(C, d,\Delta, \{E^{p,q}:C\rightarrow
C^{\otimes p}\otimes C^{\otimes q}\}\right),$ the cobar
construction $\Omega C$ is a dg Hopf algebra with a
comultiplication induced by $\{E^{p,q}\}$.

Main examples of Hirsch  (co)algebras are: $C_{\Box}^*(Q)$ (see
the previous section), in particular, Adams' cobar construction
$\Omega C_*(X)$ (\cite{SU2}), and the singular simplicial cochain
complex $C^*(X)$: In \cite{Khelaia} a twisting cochain
$E:BC^*(X)\otimes BC^*(X)\to C^*(X)$ satisfying (i)-(iii) is
constructed and  these conditions determined $E$ uniquely up to
the standard equivalence of twisting cochains.

\subsection{Multiplicative twisted tensor products}\
Let $A$ be a Hirsch algebra, $C$ be a dg Hopf algebra, and $M$  be
a dga being  a dg comodule over $C$.

\begin{definition}
A  twisting cochain $\vartheta : C\to  A$  in $Hom(C,A)$ is
multiplicative if the comultiplicative  extension $C \to  BA$ is
an algebra map.
\end{definition}

It is clear that if $\vartheta : C\to  A$ is a multiplicative
twisting element and if $g:B\to  C$ is a map of dg Hopf algebras
then the composition $\vartheta g: B\to  A $ is again a
multiplicative twisting cochain.
 The canonical projection
$BA\to  A$ provides an example of the universal multiplicative
cochain. The argument for the proof of formula (\ref{tformula})
immediately yields the following:
\begin{theorem}
Let $\vartheta : C\to  A$ be a multiplicative twisting cochain.
Then the tensor product $A\otimes M$ with the  twisting
differential $d_{\vartheta}= d\otimes Id +Id\otimes d +\vartheta\,
\cap _{-}$ becomes a dga $(A\otimes M,
d_{\vartheta},\mu_{\vartheta})$ with the twisted multiplication
$\mu_{\vartheta}$ determined by formula (\ref{tformula}).
\end{theorem}

The above theorem includes the twisted tensor product theory both
for homotopy G-algebras (\cite{KS1}) and for commutative algebras
(\cite{Proute}).

\begin{corollary} For a Hirsch algebra $A,$
the acyclic bar construction $B(A;A)$ endowed with the twisted
multiplication determined by formula (\ref{tbformula})
 asquires
a dga  structure.
\end{corollary}

\subsection{Examples}\label{bottsamelson}
 For simplicity we assume that the ground ring $R$ is a field, and all spaces
are path connected. In examples below we further explore the fact
that for a space being a suspension the corresponding homotopy
G-algebra structure is extremely simple: it consists just of
$E_{1,1}=\smile_1$ and all other operations $E_{k>1,1}$ are
trivial \cite{KS1}, and so does the corresponding Hirsch algebra
structure for the loop space on a double suspension.

  1. {\em Multiplicative models for $\Omega^2 S^2X$.}
   Given  a polyhedron $X,$
   consider the space $Y=\Omega S^2X.$
As in \cite{KS1} we regard a suspension $SX$  as the geometric
realization of   quotient simplicial set $C_{+}X\cup
C_{-}X/C_{-}X.$
   It is
immediate to check by (\ref{cohirsch}) that  $E^{p,q}=0$ for
$(p,q)\neq (1,1)$ on $C^{\Box}_{\ast}(Q),$ where $Q={\bf \Omega}
S^2X, $ \, $\bf \Omega$ is the cubical set functor constructed in
\cite{KS1} (if we had $Q={\bf \Omega} SX,$ then the Hircsh
coalgebra structure would be reduced to that of homotopy
$G$-coalgebra on $C^{\Box}_{\ast}(Q);$ verification of this fact
is left to the interested reader). Furthermore,
$E^{1,1}:C^{\Box}_{\ast}(Q)\rightarrow C^{\Box}_{\ast}(Q) \otimes
C^{\Box}_{\ast}(Q)$ becomes a coassociative chain map of degree 1.
Since the comultiplication on $C^{\Box}_{\ast}(Q)$ is cocomutative
(more precisely, $C^{\Box}_{\ast}(Q)$ is a primitively generated
Hopf algebra),
  $E^{1,1}$ also induces a binary cooperation of
degree 1 on the homology denoted by
$Sq^{1,1}:H_{\ast}(Q)\rightarrow H_{\ast}(Q)\otimes H_{\ast}(Q)$.

Notice that both $(C^{\Box}_{\ast}(Q),d,{\Delta}, E^{1,1})$ and
$(H_{\ast}(Q),d=0, {\Delta}_{\ast}, Sq^{1,1})$ are Hirsch
coalgebras, thus $\Omega C^{\Box}_{\ast}(Q)$ and $\Omega
H_{\ast}(Q)$ both are dg Hopf algebras.

Similarly to \cite{KS1} the cycle choosing homomorphism $\iota :
H_{\ast}(Q)\to C^{\Box}_{\ast}(Q)$ is a dg coalgebra map and
induces an isomorphism of dg Hopf algebras
\begin{equation*}
\label{Bott} H(\Omega T\tilde{H}_{\ast}(SX))
\overset{\approx}{\longrightarrow}
 H(\Omega H_{\ast}(Q))
\overset{(\Omega\iota)_{\ast}}{\longrightarrow}H_{\ast}(\Omega
C^{\Box}_{\ast}(Q))=H_{\ast}(\Omega Y).
\end{equation*}

2. Let
 $\Omega Y\rightarrow PY\overset{\pi}{\rightarrow} Y$
be the Moore path fibration with the base $Y=\Omega S^2X.$ Let
$f:Y\rightarrow Z$ be a map, $ \Omega Y\times \Omega Z\rightarrow
\Omega Z$ be the induced action via the composition
\[\Omega Y\times \Omega Z \overset{\Omega f\times \Id}{\longrightarrow}
\Omega Z\times \Omega Z\rightarrow \Omega Z,\] and  $\Omega
Z\rightarrow E_{f}\overset{\zeta}{\rightarrow} Y$ be the
associated fibration; for simplicity assume that
  $Z$ is
the suspension and  simply connected $CW$-complex of finite type,
as well. We present two multiplicative models for the fibration
$\zeta$ using the permutocubical  model $Q\times_{\vartheta}
{\bf\Omega} Z $ with the universal truncating twisting function
$\vartheta=\vartheta_{U}: Q\rightarrow  {\bf \Omega} Q$.

Notice that the twisted  differential of the cochain complex
$(C_{\boxminus}^{\ast}(Q\times_{\vartheta} {\bf\Omega}
Z),d)=(C_{\Box}^{\ast}(Q)\otimes C_{\diamondsuit}^{\ast}({\bf
\Omega} Z),d_{\vartheta^{\#}})=(C_{\Box}^{\ast}(Y)\otimes B
C_{\Box}^{\ast}( Z),d_{\vartheta^{\#}})$ with universal
$\vartheta^{\#}:BC_{\Box}^{\ast}(Q)\rightarrow C_{\Box}^{\ast}(Q)$
becomes the form
\begin{multline*}
d_{\vartheta^{\#}}(a\otimes [\bar{m}^1|...|\bar{m}^n]) = da\otimes
[\bar{m}^1|...|\bar{m}^n] +\sum_{k=1}^{n}\, a\otimes
[\bar{m}^1|...|d\bar{m}^k|...|\bar{m}^n]+
 \\
a\cdot m_1\otimes [\bar{m}^2|...|\bar{m}^n].
\end{multline*}
Since the simplified  structure of the Hirsch algebra
$(C_{\Box}^{\ast}(Q),d,\mu, E_{1,1})$  formula (\ref{tformula})
becomes the following form:
\begin{equation}
\label{multc} \mu_{\vartheta^{\#}}((a_1\otimes m_1)(a_2\otimes
m_2))=a_1a_2\otimes m_1m_2+a_1E_{1,1}(f^{\#}(m^1_1),a_2)\otimes
[\bar{m}^2_1|...|\bar{m}^n_1]\cdot m_2,
\end{equation}
where $f^{\#}:C_{\Box}^{\ast}( Z)\rightarrow C_{\Box}^{\ast}(Q)
,$\, $a_1,a_2\in C_{\Box}^{\ast}(Q), \,
m_1=[\bar{m}^1_1|...|\bar{m}^n_1],\,m_2\in B C_{\Box}^{\ast}(
Z),\\
n\geq 0.$

So that we  get that $H(C_{\Box}^{\ast}(Y)\otimes B
C_{\Box}^{\ast}( Z),d_{\vartheta^{\#}},\mu_{\vartheta^{\#}})$ and
$H^*(E_f)$ are isomorphic as algebras.

On the other hand, let us consider the following multiplicative
twisted tensor product $(H^{\ast}(Y)\otimes H^{\ast}({\bf \Omega}
Z),d_{\vartheta^{\ast}})=(H^{\ast}(Y)\otimes B H^{\ast}(
Z),d_{\vartheta^{\ast}})$ with universal
$\vartheta^{\ast}:BH^{\ast}(Y)\rightarrow H^{\ast}(Y)$. The
differential here is of the form:
$$
d_{\vartheta^{*}}(a\otimes [\bar{m}^1|...|\bar{m}^n]) = a\cdot
m_1\otimes [\bar{m}^2|...|\bar{m}^n].
$$
Again since the simplified structure of the Hirsch algebra
$(H^{\ast}(Y),d=0, \mu^{\ast}, Sq_{1,1})$ the formula
(\ref{tformula}) becomes the following form:
\begin{equation}
\label{multh} \mu_{\vartheta^{\ast}}((a_1\otimes m_1)(a_2\otimes
m_2))=a_1a_2\otimes m_1m_2+a_1Sq_{1,1}(f^{\ast}(m^1_1),a_2)\otimes
[\bar{m}^2_1|...|\bar{m}^n_1]\cdot m_2,
\end{equation}
where $ f^{\ast}:H^{\ast}(Z)\rightarrow H^{\ast}(Y),\, a_1,a_2\in
H^{\ast}(Y),\, $ $m_1=[\bar{m}^1_1|...|\bar{m}^n_1],\, m_2\in B
H^{\ast}( Z),\\ n\geq 0.$  Remark that for  an element $a\in
H^{\ast}(Y),$ one gets $Sq_{1,1}(a,a)=Sq_1(a),$ the Steenrod
square.

We claim that $(H^{\ast}(Y)\otimes B
H^{\ast}(Z),d_{\vartheta^{\ast}})$ is a "small" multiplicative
model of the fibration $\zeta,$ i.e., $H(H^{\ast}(Y)\otimes B
H^{\ast}(Z),d_{\vartheta^{\ast}})$ and $H^*(E_f)$ are isomorphic
as algebras. Indeed, since the explicit formulas (\ref{multc}) and
(\ref{multh}) it is straightforward to check that a "cocyle
choosing" homomorphism $(H^{\ast}(Y)\otimes B H^{\ast}(
Z),d_{\vartheta^{\ast}} )\rightarrow (C^{\ast}(Y)\otimes B
C^{\ast}( Z),d_{\vartheta^{\#}})$ induces an algebra isomorphism
\[
H\left(H^{\ast}(Y)\otimes B H^{\ast}(
Z),d_{\vartheta^{\ast}}\right)\overset{\approx}{\longrightarrow}
H\left(C^{\ast}(Y)\otimes B C^{\ast}( Z),d_{\vartheta^{\#}}
\right)\approx H^*(E_f) \]
 as required.

As a byproduct we obtain that the multiplicative structure of the
total space $E_f$ does not depend on a map $f$ in a sense that if
$f^{\ast}=g^{\ast}$ then $H^{\ast}(E_f)=H^{\ast}(E_g)$ as
algebras.  Note also that this multiplicative structure is purely
defined by the $\smile,$  $\smile_1$ and $\smile_2$  operations on
the simplicial cochain complex $C^{\ast}(S^2X).$

\vspace{0.2in}


\end{document}